\documentclass[12pt]{article}
\usepackage[utf8]{inputenc}

\usepackage{amsthm, amsfonts, amssymb, amsmath,enumitem, natbib, bbm, mathabx}
\usepackage{cases, amsthm}
\usepackage{fullpage}
\usepackage{wasysym} 
\usepackage{fancyhdr}
\usepackage{graphicx}
\usepackage{bbm}
\usepackage{caption}
\usepackage{subcaption}
\usepackage{booktabs}
\usepackage{graphicx}
\usepackage{authblk}
 \usepackage{setspace}
\usepackage[colorlinks=true, pdfstartview=FitV,
linkcolor=blue, citecolor=blue, urlcolor=blue]{hyperref}
\usepackage{comment}
\usepackage{mathtools}
\usepackage{algpseudocode}

\newtheorem{theorem}{Theorem}
\newtheorem{prop}{Proposition}
\newtheorem{lemma}{Lemma}
\newtheorem{assumption}{Assumption}

\newtheorem{remark}{Remark}

\usepackage[linesnumbered,ruled,vlined]{algorithm2e}

\usepackage{mathrsfs}

\newcommand{\rscp}{\mathsf{RSCP}}
\newcommand{\rsp}{\mathsf{RSP}}

\newcommand{\Pbb}{\mathbb P}

\newcommand{\1}{\mathbf 1}
\newcommand{\calC}{\mathcal C}

\newcommand{\Var}{\operatorname{Var}}

\usepackage{xcolor}
\usepackage{cleveref}
\newcommand{\E}{\mathbbm{E}}

\def\E{\mathbb{E}}

\def\E{\mathbb{E}}
\usepackage[top=0.5in, bottom=0.7in, left=0.5in, right=0.5in]{geometry}
\title{\large Multiply Robust Conformal Risk Control\\
with Coarsened Data}

\author[1]{Manit Paul \thanks{Corresponding author. Department of Statistics \& Data Science, University of Pennsylvania, Academic Research Building, 265 S 37th Street, Philadelphia 19104-1686, PA 19104, USA. Email: paulman@wharton.upenn.edu}}
\author[2]{Arun Kumar Kuchibhotla}
\author[3]{Eric J. Tchetgen Tchetgen}
\affil[1,3]{Department of Statistics \& Data Science, University of Pennsylvania, PA, US}
\affil[2]{Department of Statistics \& Data Science, Carnegie Mellon University, PA, US}
\date{}

\AtBeginDocument{\doublespacing}
\begin{document}
\maketitle
\begin{abstract}
Conformal Prediction (CP) provides prediction sets with formal coverage guarantees, but most existing methods assume fully observed training and calibration data. In many applications, however, outcomes, covariates, or both may be missing or coarsened. We develop a general framework for predictive inference with coarsened data by combining semiparametric theory with conformal methods.

We introduce two complementary procedures. $(i)$ Robust Split Prediction ($\rsp$) generalizes the doubly robust estimating-equation approach of \cite{yang2024doubly} beyond the classical missing outcome problem by constructing an efficient influence function based estimating equation and solving its empirical analogue. This yields multiply robust prediction sets under general coarsening patterns. $(ii)$ Robust Split Conformal Prediction ($\rscp$) instead combines semiparametric theory with conformal risk control. Because efficient influence functions are generally non-monotone in the threshold, we first estimate the target threshold on an independent sample and freeze the threshold-dependent nuisance functions at this preliminary estimate. This restores monotonicity and allows standard conformal risk control arguments to be applied. The resulting coverage guarantees inherit the higher-order mixed-bias structure of the influence function. We develop the theory for classical, multi-stage monotone, and non-monotone missingness and illustrate the performance of the methods through simulations and a longitudinal rheumatoid-arthritis application.
\end{abstract}

\section{Introduction and Motivation}
\label{sec:introduction}

Prediction algorithms are widely used to forecast future outcomes, inform
decision-making, and quantify uncertainty across domains such as healthcare,
finance, and the social sciences. Given a training data set
$S=\{(X_1,Y_1),\ldots,(X_n,Y_n)\}$ consisting of $n$ observations drawn
i.i.d.\ from a distribution $\mathcal P$ on $\mathcal X\times\mathbb R$, a
new observation $(X_{n+1},Y_{n+1})\sim\mathcal P$, and a nominal coverage
level $1-\alpha\in(0,1)$, the usual predictive-inference objective is to
construct a prediction set
$\widehat{\mathcal C}_{n,\alpha}(X_{n+1})$ satisfying,
\begin{equation}
\label{eq:gen_pred}
\mathbb P
\left(
Y_{n+1}\in
\widehat{\mathcal C}_{n,\alpha}(X_{n+1})
\right)
\geq
1-\alpha.
\end{equation}
Conformal prediction, originating with \cite{vovk2005algorithmic}, provides
a general approach for obtaining such finite-sample coverage guarantees
under exchangeability without imposing a parametric model on the joint
distribution of $(X,Y)$.

Since the development of the inductive conformal method of
\cite{papadopoulos2002inductive}, also known as split conformal prediction
\citep{lei2014distribution}, and the transductive conformal method of
\cite{vovk2005algorithmic}, the conformal-prediction literature has
developed in several important directions. One line of research considers
methods such as jackknife$+$, CV$+$, and bootstrap-based variants that use
the available data more efficiently while retaining computational
tractability; see, for example, \cite{barber2021predictive},
\cite{romano2020classification}, and \cite{kim2020predictive}. We refer the
reader to \cite{kuchibhotla2020exchangeability} for a detailed discussion.
While these methods exploit exchangeability to provide marginal
coverage, another line of work studies guarantees conditional on the
training data, typically under i.i.d.\ sampling and additional regularity
conditions; see, for example, \cite{gyorfi2019nearest} and
\cite{yang2021finite}.

A related literature seeks prediction guarantees that are more local than
marginal coverage. Exact distribution-free conditional coverage is
impossible without substantial loss of efficiency in sufficiently rich
distribution classes \citep{vovk2012conditional}. This has motivated
intermediate notions such as group-conditional or approximately
conditional coverage; see, among others, \cite{vovk2003mondrian},
\cite{foygel2021limits}, \cite{jung2022batch}, and
\cite{gibbs2023conformal}. Other methods, including
\cite{romano2019conformalized}, \cite{sesia2021conformal}, and
\cite{chernozhukov2021distributional}, provide methods with asymptotic conditional validity guarantees.

A separate challenge arises when the observations used for calibration are
not distributed in the same way as a fully observed future observation.
The covariate-shift problem was studied from a conformal perspective by
\cite{tibshirani2019conformal}; subsequent work includes
\cite{lei2021conformal}, \cite{qiu2023prediction},
\cite{cauchois2024robust}, and \cite{yang2024doubly}. More general forms of
non-exchangeability have also been studied, for example for dependent time
series \citep{stankeviciute2021conformal} and through weighted conformal
methods \citep{barber2023conformal}.

The focus of the present paper is different but closely related. We allow
the data used for training and calibration to be \emph{coarsened}: the
outcome, auxiliary covariates, or both may be missing for some calibration
units. The goal is to obtain a valid prediction set for the outcome when only the baseline covariates are observed for the test point. A canonical example is a longitudinal cohort study in which
baseline covariates are observed for all subjects, additional time-varying
covariates are collected during follow-up, and subjects may drop out before
the final outcome is recorded. The time-varying variables may have little
direct predictive interest for a future patient, but they can be essential
for correcting selection bias caused by dropout in the calibration data.

The central methodological difficulty is that the distribution of the
target outcome or non-conformity score is not directly observed in the
calibration sample. Semiparametric theory provides observed-data
representations of the target distribution through efficient influence
functions. Such representations are attractive because the bias in population moments under nuisance-function estimation often has a higher-order
mixed-bias structure. However, using an efficient influence function inside
a conformal procedure creates an additional obstacle: although its {expectation} as a function of the threshold can be identified as a monotone
distribution function, the influence function itself generally need not be
monotone in the threshold. Consequently, the standard conformal risk
control argument cannot in general be applied to the estimated influence
function directly.

This observation leads to the two complementary procedures studied in this
paper. The first procedure directly solves the empirical estimating
equation obtained from the efficient influence function. This approach is
the natural extension of the doubly robust calibration method of
\cite{yang2024doubly}, and we call its general version Robust Split
Prediction ($\rsp$). Because the threshold is obtained as a root of an
empirical estimating equation, its coverage analysis involves using techniques from empirical process theory. Thus, even when all the nuisance functions are estimated correctly, a finite-sample stochastic remainder remains.

Our second procedure instead modifies the influence-function estimating
equation so that conformal risk control approach can be used. We first construct an
independent preliminary estimate $\widehat r_{\alpha,0}$ of the desired
quantile. We then evaluate, or ``freeze,'' all threshold-dependent nuisance
functions at $\widehat r_{\alpha,0}$. In the classical and monotone missing-data problems studied in detail in the main text, the remaining dependence on the final threshold is through a monotone term such as $\mathbf 1\{R\leq\theta\}$, allowing us to apply the conformal finite-sample correction to the modified estimating function. We call the
resulting procedure Robust Split Conformal Prediction ($\rscp$). The same idea also extends to the non-monotone MAR setting by freezing all threshold-dependent nuisance components at the preliminary estimate; for ease of exposition, we primarily use $\rsp$ when discussing and empirically studying the non-monotone setting. This
construction is the principal contribution of the second approach. A
formal review of conformal risk control and its relation to our construction
is provided in \Cref{app:crc_background}.


\subsection{Notation and semiparametric ingredients}
\label{subsec:intro_notation}

We use uppercase letters for random variables or random vectors and
lowercase letters for their realizations. The full-data vector is denoted by $\widetilde O=(X,Z,Y)$, where $X$ denotes baseline covariates available for prediction, $Z$ denotes
auxiliary variables that may be useful for adjusting for coarsening, and
$Y$ denotes the target outcome. The observed data are written as $O=(C,\mathfrak G_C(\widetilde O))$, where $C$ denotes the coarsening level and
$\mathfrak G_C(\widetilde O)$ denotes the portion of the full data that is
observed under that coarsening pattern. We reserve $C=0$ for the pattern in
which only the baseline information required for prediction is observed and
$C=\infty$ for complete observation of the full data. These notations are standard in the missing data literature (see \cite{tsiatis2006semiparametric} for example).

In the general coarsening formulation, $L$ denotes the number of
intermediate coarsening levels. In the monotone longitudinal setting, $D$
denotes the number of intermediate auxiliary-variable stages. In the
non-monotone setting, $K$ denotes the number of distinct nontrivial
coarsening patterns. The generic collection $\eta^*=\{\eta_i^*\}_{i=1}^{K_\eta}$ denotes nuisance functions appearing in an identification equation of the coverage probability using the observed-data distribution. These nuisance functions are features of the observed-data law, such as
propensity scores, coarsening probabilities, conditional distribution
functions, and iterated outcome regressions involving $X$ and $Z$. Their
precise form depends on the missing-data model under consideration.

For completeness, we briefly recall the role of an influence function.
Let $\Psi(P)$ be a pathwise differentiable functional of an observed-data
distribution $P$. For a regular parametric submodel
$\{P_\varepsilon:\varepsilon\in(-\epsilon,\epsilon)\}$ through $P$ with
score $S(O)$, an influence function $\varphi_P(O)$ satisfies,
\[
\left.
\frac{d}{d\varepsilon}
\Psi(P_\varepsilon)
\right|_{\varepsilon=0}
=
\mathbb E_P[\varphi_P(O)S(O)],
\qquad
\mathbb E_P[\varphi_P(O)]=0.
\]
The efficient influence function is the influence function with the smallest variance among all influence functions representing the first-order behavior of the target functional. In the problems considered here,
$\Psi(P)$ is a distribution function or quantile-related functional for
the target population. If
$\mathrm{IF}(O,\theta,\eta^*)$ denotes the corresponding efficient
influence function, then
\[
f(O,\theta,\eta^*)
=
\mathrm{IF}(O,\theta,\eta^*)
+
\mathbb P(R(X,Y)\leq\theta)
\]
provides an observed-data representation satisfying $\mathbb E[f(O,\theta,\eta^*)]
=
\mathbb P(R(X,Y)\leq\theta).$ This identity is the bridge between semiparametric theory and the
prediction problems studied below.

We use the term \emph{multiple robustness} in its standard semiparametric
sense. An estimating function is multiply robust if its population bias
vanishes under several distinct configurations of correctly specified
nuisance components. For example, in a $(D+1)$-stage monotone
missing-data problem, the bias of the efficient influence function has the schematic form $\sum_{j=1}^{D+1} A_jE_j$, where $A_j$ and $E_j$ denote errors in estimating propensity scores and
the outcome-regression functions respectively. Hence the leading bias vanishes whenever, at every stage, either $A_j=0$ or $E_j=0$, producing $2^{D+1}$ robust configurations. 

\subsection{Related literature on coarsened data}
\label{sub_sec:related_lit}

It is well known that in most applied settings, a nontrivial portion of
the training or calibration sample may be missing or coarsened regardless
of the overall sample size; see, for example,
\cite{little2019statistical} and \cite{josse2018introduction}. Coarsened
data arise across many substantive areas: patient dropout in randomized
trials and observational studies, unanswered income questions in surveys,
equipment malfunction in environmental sensing, trading-day absences in
financial data, or connectivity-related losses in digital telemetry.
Such missingness will rarely be completely at random. It is therefore
important to account for potential bias due to missingness, a
task that standard machine-learning algorithms do not generally address; see, for example,
\cite{josse2024consistency},
\cite{le2020neumiss,le2020linear,le2021sa},
\cite{ayme2022near}, and \cite{van2023missing}.

Missing and coarsened data have been extensively studied in the
semiparametric literature; see, among others, \cite{robins1994estimation}
and \cite{tsiatis2006semiparametric}. The principal focus here, however, is estimation and inference for population
functionals rather than predictive sets for future outcomes. Our goal is
to use the same tools from the semiparametric literature
to construct prediction sets whose coverage error can be characterized
under estimation of the true nuisance functions.


Uncertainty quantification using conformal prediction with missing
covariates is studied by \cite{zaffran2023conformal}, who use an
impute--then-predict approach and establish marginal validity under
symmetrical imputation of calibration and test observations.
\cite{zaffran2024predictive} study prediction sets which are conditionally valid
given missingness patterns under a more restrictive setting in which the
missingness pattern is independent of both the covariates and the response.
These works do not directly address missing outcomes, which are central to
our setting.

Conformal methods for coarsened outcomes are also considered by
\cite{cauchois2024predictive}. Their observed response identifies a set of
possible full outcomes (called the weak set), and the target is to construct a prediction set intersecting
that weak set. This is distinct from our objective, which is prediction of
the underlying target response itself after accounting for the coarsening
mechanism.

The classical missing-outcome problem has received significant attention in recent years.
In this problem, the outcome is observed only for one subpopulation and the goal is
to construct a prediction set for a target unit whose outcome is
unobserved. Weighted conformal prediction \citep{tibshirani2019conformal}
gives finite-sample validity when the relevant density ratio or propensity
score is known, while \cite{lei2021conformal} studies the effect of
estimating this ratio on the coverage probability. Conformal risk control
\citep{angelopoulos2022conformal} instead exploits a correctly specified
conditional outcome distribution to construct valid prediction set. Both approaches are singly robust in
the sense that their validity relies on the accurate estimation of one particular nuisance component. \cite{yang2024doubly} use the efficient influence function of the target
quantile to obtain a doubly robust estimating equation under covariate
shift/classical missingness. Their procedure remains asymptotically valid
when either the propensity score or the outcome-regression component is
consistently estimated. 

We further point to two recent works \citep{angelopoulos2026nonmonotonic,aldirawi2026nonmonotone} arxived after the first version of our work was released on arxiv. These two works relax the monotonicity requirement of conformal risk control and therefore provide, in principle, an alternative way of handling the non-monotonicity of an efficient influence function based estimating function. Conditional on the sample used to estimate the nuisance functions, and under bounded inverse-probability weights so that the resulting influence function based loss can be treated as bounded, one could apply these methods directly to the estimated influence function. However, doing so treats the influence function as an essentially generic non-monotone loss and therefore does not exploit the additional structure available in our problem. For general bounded non-monotone losses, \citet{angelopoulos2026nonmonotonic} obtains a stability-based excess risk term of order $\widetilde O(n^{-1/2})$; the sharper $O(n^{-1})$ result in that work requires continuity and Lipschitz smoothness of the loss together with a strong crossing condition, whereas the influence functions considered in our work contain the discontinuous threshold term $\mathbf 1\{Y\leq\theta\}$. Similarly, \citet{aldirawi2026nonmonotone} obtain, over a finite grid of $m$ candidate thresholds, a generic correction of order $O\{\sqrt{\log(m)/n}\}$, and show that this rate is minimax-optimal over arbitrary bounded non-monotone losses. Thus, to guarantee risk at the nominal level, their procedure must calibrate at the more conservative level $\alpha-D(m,n)$, where $D(m,n)=O\{\sqrt{\log(m)/n}\}$. In contrast, $\rsp$ exploits the bounded-variation and the second-moment structure of the influence-function and, when the nuisance estimation errors vanish, has coverage slack $O(\log n/n^{3/4})$. The $\rscp$ construction goes further by freezing the threshold-dependent nuisance functions at an initial estimate of the target quantile, thereby restoring pointwise monotonicity and allowing the standard CRC argument to be applied without any penalty coming from using a generic non-monotone-loss function. The coverage slack of $\rscp$ vanishes in finite samples when the nuisance estimation errors disappear. Hence, although generic non-monotone CRC methods could be applied to the estimated influence function, the specialized $\rsp$ and $\rscp$ constructions exploit the structure of the underlying efficient influence functions to obtain substantially sharper guarantees in the coarsened-data settings considered in this work.


\subsection{Contributions and Organization}

In this work we study general missing/coarsened data, that is, missing or coarsened covariates or outcomes in the training and calibration samples, with the goal of establishing how conformal methods can be used to obtain valid prediction sets for the target outcome. Throughout the paper, we consider the setting in which a designated subset of the covariates (termed baseline covariates) is observed for all units in the calibration sample, and we anticipate making predictions for the unobserved outcome in a new sample based solely on their observed baseline covariates. However, in the calibration sample the outcome of interest may be missing for some units. Fortunately, we suppose that we are in the common setting in which a rich set of covariates are also potentially observed in the calibration sample (auxiliary data that one does not necessarily expect to observe for future units). The value of these auxiliary variables is that, although they are not necessarily of primary scientific interest, they can be essential for accounting for complex dependent coarsening. Our formal development focuses primarily on missing-at-random and coarsening-at-random models under which the relevant full-data functionals are identified. Under coarsening not at
random, the full-data distribution is generally not identified from the
observed data without additional restrictions, sensitivity parameters, or
an explicit selection model. Therefore CNAR
models are not studied in this paper. The contributions of this paper are as follows:

\begin{enumerate}
\item We develop the Robust Split Prediction ($\rsp$) framework, which generalizes the doubly robust prediction method of \cite{yang2024doubly} from the classical missing-outcome/covariate-shift problem to substantially more general coarsened-data settings. The key idea is to derive an observed-data estimating function from the efficient influence function (EIF) of the target quantile and to estimate the prediction threshold by directly solving the corresponding empirical moment equation. In the classical missing-data problem, $\rsp$ reduces to the doubly robust construction of \cite{yang2024doubly}. We sharpen the coverage guarantee available for this method in that setting and, more importantly, show how the same influence-function-based construction can be systematically extended to considerably more complex missing-data structures.

\item We introduce a second method, Robust Split Conformal Prediction ($\rscp$), which combines efficient influence functions with conformal risk control \citep{angelopoulos2022conformal}. A central difficulty is that, although the expectation of the EIF identifies a monotone distribution function of the threshold $\theta$, the EIF itself is generally not monotone in $\theta$. Therefore, it cannot in general be inserted directly into the standard conformal risk control framework. To overcome this difficulty, $\rscp$ first computes an initial estimate $\widehat r_{\alpha,0}$ of the target threshold using an independent data split and then freezes the threshold-dependent nuisance functions at $\widehat r_{\alpha,0}$. After freezing the threshold-dependent nuisance functions, the resulting estimating function is monotone in the final threshold in the classical and monotone missing-data problems developed in detail below, allowing us to apply a conformal finite-sample correction. The same construction is also applicable under non-monotone MAR when all threshold-dependent outcome-regression and operator-derived nuisance functions are frozen at the preliminary threshold. This distinction between $\rsp$ and $\rscp$ is important: $\rsp$ directly solves a non-monotone empirical estimating equation and consequently retains a vanishing empirical-process contribution to the coverage error, whereas $\rscp$ can provide finite-sample coverage validity when the relevant nuisance components are correctly specified. 

\item The efficient influence functions used in both constructions typically depend on unknown nuisance functions that must be estimated from the observed data. Because these influence functions have higher-order mixed-bias structure \citep{rotnitzky2021characterization}, the corresponding prediction procedures naturally inherit double or multiple robustness. We explicitly quantify the resulting coverage error in terms of nuisance-estimation errors. In the classical missing-data problem (\Cref{sec:cov_shift}), the leading bias is a product of the propensity-score and outcome-regression errors. In the multi-stage monotone missingness problem (\Cref{sec:monotone}, \Cref{app:multi_stage_monotone}), the leading bias becomes a sum of stagewise products: at every stage it is sufficient for either the corresponding propensity-score component or the outcome-regression component to be correctly specified for that stage's leading bias term to vanish. This yields the $2^{D+1}$-robust structure described in \Cref{rem:multi_stage_multiple_robustness}.

\item We show how these ideas yield systematic constructions of valid prediction sets under missing-data structures for which, to the best of our knowledge, no general predictive inference procedure is currently available. In particular, we develop prediction methods for both multi-stage monotone (\Cref{app:multi_stage_monotone}) and non-monotone missingness (\Cref{sec:non_monotone}). In the multi-stage monotone setting, the calibration data may arise from a longitudinal cohort study in which each fully observed data vector has the form $(X,Z^1,\ldots,Z^D,Y)$. Here, $X$ denotes baseline covariates that are always observed, $Z^j$ denotes auxiliary covariates collected at the $j$-th follow-up visit, and the target outcome $Y$ is recorded at the final stage. Under monotone dropout, once an individual leaves the study, all subsequent auxiliary variables and the final outcome are missing. We derive the corresponding multi-stage efficient influence function and show how $\rsp, \rscp$ produce multiply robust prediction sets in this setting.

We also consider non-monotone missingness, where the observed-data patterns need not be nested, for example, a participant may miss an intermediate visit and subsequently re-enter the study. Existing conformal prediction methods do not provide a systematic way of constructing valid prediction sets under such general coarsening patterns. Using the efficient influence function available under non-monotone MAR, the general $\rsp$ framework provides doubly robust prediction sets. The $\rscp$ construction is also applicable in this setting: by freezing all threshold-dependent outcome-regression and operator-derived nuisance functions at an independent preliminary root, the remaining dependence on the final threshold is through the positively weighted indicator term, and the resulting estimating function is monotone. For ease of exposition, we develop and empirically illustrate the $\rsp$ construction for non-monotone missingness in this paper.

\item We demonstrate the practical behavior of the proposed methods in both simulated and real-data examples. Our simulations in \Cref{subsec:simulation} and \Cref{app:add_simul} compare $\rsp$ and $\rscp$ across different calibration sample sizes and under model misspecification, using different relevant non-conformity score functions. We additionally compare the proposed procedures with existing approaches in the classical missing-data problem in \Cref{tab:cov_compare_simul} in \Cref{sec:cov_shift}. Finally, we apply the methods to a longitudinal cohort study of patients with rheumatoid arthritis \citep{choi2002methotrexate}, where the goal is to obtain valid prediction sets for total tender joint count at a later follow-up visit using baseline information while accounting for incomplete intermediate covariates and outcomes. We consider both monotone and non-monotone missingness structures in this application.

\end{enumerate}

\paragraph{{\bf Organization}}
We use the classical missing-data problem in
\Cref{sec:cov_shift} as a motivating example that makes the distinction
between $\rsp$ and $\rscp$ transparent; additional background on the connection with the
covariate shift problem is in \Cref{app:cov_shift}. In
\Cref{sec:gen_missingness}, we study the construction of prediction set for general coarsened-data and establish coverage guarantees for $\rsp$ and $\rscp$. In
\Cref{sec:monotone}, we specialize the theory developed for general missingness problems to the two-stage monotone missingness problem. The results for the multi-stage monotone missingness problem have been deferred to \Cref{app:multi_stage_monotone}. The construction of valid prediction set for non-monotone missingness 
problem is studied in \Cref{sec:non_monotone}. Empirical results (for both synthetic and real data experiments) are presented in \Cref{sec:simulation}. \Cref{sec:extension} concludes the main text with a discussion of limitations of the methods and
future directions. A self-contained review of conformal risk control is
provided in \Cref{app:crc_background}; proofs and additional technical
results are collected in the remaining appendix sections.

\section{Motivating example: Classical missing data problem}
\label{sec:cov_shift}
In this section, we discuss the construction of valid prediction intervals in the context of the classical missing data problem. The observed dataset $\mathcal D$ is partitioned into two groups: the training $\mathcal D^{tr}$ and the calibration  dataset $\mathcal D^{cal}$. The training dataset comprises $n_0$ i.i.d.\ data points $O_i = (C_i^t, \mathfrak G_{C_i^t}(\tilde O_i^t))$ (for all $i \in [n_0]$). The calibration dataset comprises $n$ i.i.d.\ data points $O_i = (C_i, \mathfrak G_{C_i}(\tilde O_i))$ (for all $i \in [n]$). Here $\tilde O_i $ (similarly $\tilde O_i^t$) denotes the complete observation $ (X_i, R_i)$, $C_i$ (similarly $C_i^t$) is a random coarsening level that takes values in $\{0, \infty\}$ and $\mathfrak G_0(\tilde O_i) = X_i, \mathfrak G_{\infty}(\tilde O_i) = (X_i, R_i)$ where $R_i = R(X_i, Y_i)$ is the non-conformity score for the $i$-th data point computed using a pre-trained non-conformity score function $R(\cdot, \cdot)$. We make the missing at random (MAR) assumption,
\begin{equation}
    \label{eq:mar_cmdp}
    R_i \perp C_i|X_i \quad \mbox{for}\quad i \in [n].
\end{equation}
For identification we further assume that the distribution of $X_i|C_i = 0 $ is absolutely continuous with respect to the distribution of $X_i|C_i = \infty $. Given a new test point $(0, X_{n+1})$ with missing outcome $Y_{n+1}$, our objective is to construct a valid $(1-\alpha)$-prediction set $\widehat C_{n,\alpha}(X_{n+1}) = \{y: R(X_{n+1}, y) \leq \widehat r_{\alpha} \}$ for $Y_{n+1}$. In other words, we want to find $\widehat r_{\alpha}$ such that,
\begin{equation}
    \label{eq:cmdp}
    \mathbb{P}\left( R(X_{n+1}, Y_{n+1}) \leq \widehat r_{\alpha} |C_{n+1} = 0\right) \geq 1 - \alpha.
\end{equation}
Since $R_i$ is not observed when $C_i = 0$, it is not possible to find any non-trivial $\widehat r_{\alpha}$ such that \eqref{eq:cmdp} holds in finite samples (refer to Theorem $1$ of \cite{yang2024doubly} for a rigorous proof of this impossibility result) without any additional assumptions. Therefore instead of targeting finite sample validity, we focus on the construction of asymptotically valid prediction sets in this context. We denote the propensity score function $\mathbb{P}(C = 0| X = x)/\mathbb{P}(C = \infty|X = x)$ by $\pi^*(x)$ and the outcome regression function $\mathbb P(R \leq r|X = x)$ by $m_R^*(r, x)$. The Weighted Conformal Prediction (WCP) introduced in \cite{tibshirani2019conformal} computes $\widehat r_{\alpha}$ by using the following equivalent representation of the coverage probability, 
\[
 \mathbb{P}\left( R(X_{n+1}, Y_{n+1}) \leq \widehat r_{\alpha} |C_{n+1} = 0\right) = 1- \alpha +  \mathbb E[\textbf{1}\{C_{n+1} = \infty\} \pi^*(X_{n+1})[\textbf{1}\{R(X_{n+1}, Y_{n+1}) \leq \widehat r_{\alpha}\} - (1-\alpha)] ],
\]
and shows that the resulting estimate guarantees a coverage of $1 - \alpha - \mathfrak C\|\widehat \pi -  \pi^*\|_2$ where $\mathfrak C$ is used to denote a fixed constant and $\widehat \pi$ is an estimate of $\pi^*$ obtained from the training data. The conformal risk control approach proposed in \cite{angelopoulos2022conformal} can be used to compute $\widehat r_{\alpha}$ by using a different representation of the coverage probability,
\[
\mathbb{P}\left( R(X_{n+1}, Y_{n+1}) \leq \widehat r_{\alpha} |C_{n+1} = 0\right) = 1- \alpha +  \mathbb E[\textbf{1}\{C_{n+1} = 0\} [m_R^*(\widehat r_{\alpha}, X_{n+1}) - (1-\alpha) ]].
\]
The resulting prediction set maintains a coverage of at least $1-\alpha - \mathfrak C \sup_{\theta}\| \widehat m_R(\theta, \cdot) - m_R^*(\theta, \cdot) \|_2 $ where $\widehat m_R(\theta, \cdot)$ is an estimate of $m_R^*(\theta, \cdot)$ obtained from the training data. Refer to \Cref{app:cov_shift} for further details on the prior work. 
Note that both of these approaches rely on the correct specification of a nuisance function ($\pi^*(\cdot)$ for WCP and $m_R^*(\cdot, \cdot)$ for conformal risk control) for the required coverage guarantee of $1-\alpha$; i.e.\ the prediction sets computed using these methods result in a first order coverage slack. In order to ensure that the coverage guarantee holds if either of the two nuisance functions are correctly specified, \cite{yang2024doubly} computes $\widehat r_{\alpha}$ using the following alternate representation of the coverage probability,
\begin{equation}
\begin{split}
    \label{coverage_link_cmdp}
     \mathbb{P}\left( R(X_{n+1}, Y_{n+1}) \leq \widehat r_{\alpha} |C_{n+1} = 0\right) = & 1 - \alpha + \frac{\mathbb{E}\left[\mathrm{IF}(\widehat r_{\alpha},C, \mathfrak G_C((X, R)), \pi^*, m_R^*) \right]}{\mathbb{P}(C = 0)}
\end{split}
\end{equation},
where $\mathrm{IF}(r_{\alpha}, c, \mathfrak G_c((x, r)), \pi^*, m_R^*)$ is the efficient influence function of the $(1-\alpha)$-th quantile of $R|C = 0$,
\begin{equation}
    \label{eq:eff_influ_cmdp}
 \mathrm{IF}(r_{\alpha},c, \mathfrak G_c((x, r)), \pi^*, m_R^*) = \textbf{1}\{c = \infty\}\pi^*(x)\left[\textbf{1}\{r \leq r_{\alpha}\}  - m_R^*(r_{\alpha},x)  \right] + \textbf{1}\{c = 0\}\left[m_R^*(r_{\alpha},x) - (1 - \alpha)\right].
\end{equation}
Here $r_{\alpha}$ is the $(1-\alpha)$-th quantile of $R|C = 0$. \cite{yang2024doubly} shows that if $\widehat r_{\alpha}$ is defined using the efficient influence function \eqref{eq:eff_influ_cmdp}, 
\begin{equation}
    \label{eq:dr_cmdp}
    \widehat r_{\alpha} = \inf \left\{\theta \in \mathbb{R}: \frac{1}{n}\sum_{i = 1}^n \mathrm{IF}(\theta, X_i, O_i, \widehat \pi, \widehat m_R)  \geq 0 \right\},
\end{equation}
then the resulting prediction set has the following double robust coverage slack (see Theorem $4$ of \cite{yang2024doubly}), 
\[
\mathbb{P}\left( R(X_{n+1}, Y_{n+1}) \leq \widehat r_{\alpha} |C_{n+1} = 0\right) \geq 1 - \alpha - \mathfrak{C}_1\|\widehat \pi - \pi^*\|_2 \sup_{\theta \in \mathbb{R}} \|\widehat m_R(\theta, \cdot) - m_R^*(\theta, \cdot) \|_2 - \frac{\mathfrak{C}_2}{\sqrt{n}},
\]
where $\mathfrak C_1, \mathfrak C_2$ are constants. The key phenomenon behind this improved coverage slack is the double robustness of the underlying efficient influence function (see \cite{robins1994estimation} and Theorem $2$ of \cite{yang2024doubly}). We will see in the following sections that this method can be generalized to more complex missing data problems. To maintain naming consistency throughout the paper, we shall refer this method as the robust split prediction ($\rsp$) algorithm. The following theorem shows that the lower bound of \cite{yang2024doubly} on the coverage of the $\rsp$ prediction set constructed using $\widehat r_{\alpha}$ (defined in \eqref{eq:dr_cmdp}) can be made sharper. 
\begin{theorem}
   \label{thm:influence_old_cmdp}
Suppose $\theta \mapsto \widehat m_R(\theta, x)$ is a monotone continuous function of $\theta$ for each $x$. We also assume that \eqref{eq:mar_cmdp} holds, $R(X,Y)$ is a continuous random variable and $\widehat \pi(X), \pi^*(X) \leq \overline \pi < \infty$. Then the $\rsp$ prediction set $\widehat C_{n,\alpha}(X_{n+1}) = \{y: R(X_{n+1}, y) \leq \widehat r_{\alpha}\}$, with $\widehat r_{\alpha}$ as defined in \eqref{eq:dr_cmdp}, satisfies the following conditional on the training data, 
\begin{equation*}
    \begin{split}
        &\left| \mathbb{P}( R(X_{n+1},Y_{n+1}) \leq \widehat r_{\alpha}| C_{n+1} = 0) - (1-\alpha) \right| \\
         \leq & \frac{a_\pi a_m}{\mathbb P(C = 0)}+  \frac{\mathfrak C_n}{\mathbb P(C = 0)}\left\{\frac{1}{n^{3/4}} +  \frac{\{ \sqrt{a_\pi a_m} + a_m \}}{n^{1/4}} \left[\frac{1}{n^{1/4}}+ \sqrt{a_\pi a_m} \right]\right\}, 
    \end{split}
\end{equation*}
where $a_\pi = \|\widehat \pi - \pi^*\|_2$ and $a_m = \sup_{\theta}\| \widehat m_R(\theta, \cdot) - m_R^*(\theta, \cdot)\|_2$ respectively denote the errors in estimating the nuisance functions and $\mathfrak C_n = \mathcal C(\log n)$ for a universal constant $\mathcal C$. 
\end{theorem}
The proof of \Cref{thm:influence_old_cmdp} is provided in \Cref{app:thm_old_influ_cmdp}. \Cref{thm:influence_old_cmdp} provides a tighter coverage guarantee than Theorem $4$ of \cite{yang2024doubly} for the $\rsp$ prediction set $\widehat C_{n,\alpha}(\cdot)$ introduced in \cite{yang2024doubly}. In particular, if both the nuisance functions $\pi^*, m_R^*(\cdot, \cdot)$ are correctly specified, Theorem $4$ of \cite{yang2024doubly} shows that the coverage is at least $1 - \alpha - O(n^{-1/2})$. However, we see from \Cref{thm:influence_old_cmdp} that if both the nuisance functions are correctly specified, the slack in the coverage of the prediction set $\widehat C_{n,\alpha}(\cdot)$ instead decays at a faster rate of $O(n^{-3/4})$. 

We now introduce the robust split conformal prediction ($\rscp$) algorithm for producing prediction sets with improved coverage guarantees. The algorithm comprises the following steps.
\begin{enumerate}
    \item Split the training dataset into two parts: $\mathcal D^{tr}_1 = \{O_i^{t_1}\}_{i = 1}^{n_1}$ comprising of $n_1$ points and $\mathcal D^{tr}_2 = \{O_i^{t_2}\}_{i = 1}^{n_2}$ comprising of $n_2$ points where $n_0 = n_1 + n_2$. 
    \item Obtain estimates $\widehat \pi(x), \widehat m_R(\theta, x)$ of the nuisance functions $\pi^*(x), m_R^*(\theta, x)$ from the sub-sample $\mathcal{D}^{tr}_1$.
    \item Compute an initial estimate $\widehat r_{\alpha,0}$ of $r_\alpha$ using the sub-sample $\mathcal{D}^{tr}_2$,
   \begin{equation}
   \label{eq:rscp_cmdp_init}
      \begin{split}
      \widehat  r_{\alpha,0} = &\inf\Bigg\{\theta \in \mathbb{R}: \frac{\sum_{i = 1}^{n_2} \textbf{1}\{C_i^{t_2} = 0\} [\widehat m_R(\theta, X_i^{t_2}) - (1 - \alpha)] }{n_2+1} -  \frac{(1 - \alpha)}{n_2+1}  \geq 0\Bigg\}.
      \end{split}
  \end{equation}
  \item Return the prediction set $\widehat C_{n,\alpha}(x)  = \{y \in \mathbb{R}: R(x,y) \leq \widehat  r_{\alpha}\}$ where $\widehat r_{\alpha}$ is computed using the calibration set $\mathcal{D}^{cal}$,  
  \begin{equation}
  \label{eq:rscp_cmdp}
      \begin{split}
      \widehat  r_{\alpha} = &\inf\Bigg\{\theta \in \mathbb{R}: \frac{\sum_{i = 1}^n f(O_i, \theta, \widehat \pi, \widehat m_R(\widehat r_{\alpha, 0}))}{n+1} -  \frac{\overline \pi}{n+1}  \geq 0\Bigg\},
      \end{split}
  \end{equation}
  where $\widehat \pi(x), \pi^*(x) \leq \overline \pi$ and for any $\theta, \theta' \in \mathbb R$,
  \begin{equation*}
      \begin{split}
     f(O, \theta, \pi, m_R(\theta')) = & \textbf{1}\{C = \infty\}\pi(X) [\textbf{1}\{R(X, Y) \leq \theta\} - m_R(\theta', X)] + \textbf{1}\{C = 0\} [m_R(\theta', X) - (1 - \alpha)]  .
      \end{split}
  \end{equation*}
\end{enumerate}
The robust split conformal prediction ($\rscp$) algorithm re-writes the estimating equation \eqref{eq:dr_cmdp} by adding an additional term of the order of $O(1/n)$ for the test point (see step $4$ of the algorithm). This additional term can not be identified by semi-parametric theory as the estimating equation does not change asymptotically. However, it turns out that the inclusion of the additional $O(1/n)$ term in the estimating equation enhances the finite-sample coverage of the prediction set. The next theorem characterizes the coverage slack in terms of the quality of the nuisance estimate $\widehat \pi(x)$ and the initial estimate $\widehat r_{\alpha,0}$.  
\begin{theorem}
    \label{thm:validity_proof}
Under the assumption that $\theta \mapsto \widehat m_R(\theta, x)$ is a monotonic right continuous function of $\theta$ for each $x$ and $\widehat \pi(X), \pi^*(X) \leq \overline \pi < \infty$, the $\widehat r_{\alpha}$ defined in \Cref{eq:rscp_cmdp} using the $\rscp$ algorithm satisfies the following conditional on both the training splits, 
\begin{equation*}
\begin{split}
    & \mathbb{P}( R(X_{n+1},Y_{n+1}) \leq \widehat r_{\alpha}| C_{n+1} = 0) 
       \geq  1 - \alpha - \frac{a_\pi a_m}{\mathbb P(C = 0)} - \frac{a_\pi \| m_R^*(\widehat r_{\alpha, 0}, \cdot) - m_R^*(r_\alpha, \cdot)\|_2}{\mathbb P(C = 0)} ,
\end{split}
\end{equation*}
where $a_\pi, a_m$ have the same definitions as in \Cref{thm:influence_old_cmdp}. Additionally if $R(X, Y)$ is a continuous random variable then the following holds,
\begin{equation*}
    \begin{split}
       & \mathbb{P}( R(X_{n+1},Y_{n+1}) \leq \widehat r_{\alpha}| C_{n+1} = 0) 
       \leq  1 - \alpha + \frac{3 \overline \pi}{\mathbb P(C = 0)(n+1) } +  \frac{a_\pi a_m}{\mathbb P(C = 0)} + \frac{a_\pi \| m_R^*(\widehat r_{\alpha, 0}, \cdot) - m_R^*(r_\alpha, \cdot)\|_2}{\mathbb P(C = 0)}.
    \end{split}
\end{equation*}
\end{theorem}
The proof of \Cref{thm:validity_proof} is provided in \Cref{app:proof_valid_cmdp}. We note that the bounds in \Cref{thm:validity_proof} hold irrespective of the definition of the initial estimate $\widehat r_{\alpha, 0}$. Therefore \Cref{thm:validity_proof} shows that the coverage slack of $\rscp$ prediction set is small if either $\pi^*(\cdot)$ is estimated well or if the estimates $\widehat m_R(\cdot, \cdot)$ and $\widehat r_{\alpha, 0}$ are close to the oracle quantities. The next proposition shows that under suitable regularity conditions the term $\| m_R^*(\widehat r_{\alpha, 0}, \cdot) - m_R^*(r_\alpha, \cdot)\|_2$ appearing in the coverage slack is small if we use equation~\eqref{eq:rscp_cmdp_init} for computing the initial estimate $\widehat r_{\alpha, 0}$.   
\begin{prop}
    \label{prop:rscp_init_error}
Under assumptions that $(1)$ $\theta \mapsto \widehat m_R(\theta, X)$ is a right continuous function, $(2)$ $|P[\textbf 1\{C = \infty\} \widehat \pi(X) [\textbf{1}\{R\leq \theta\} - \textbf{1}\{R\leq \theta'\}] \}]| \geq \mathfrak C_1| \theta - \theta'|$ for $\theta, \theta' \in [r_{\alpha }- \eta, r_{\alpha}+ \eta]$, $(3)$ $\|m_R^*(\theta, X) - m_R^*(\theta', X)\|_2 \leq \mathfrak C_2 | \theta - \theta'| $ for $\theta, \theta' \in [r_{\alpha }- \eta, r_{\alpha}+ \eta]$, we have the following upper bound with probability at least $ 1- 2 \delta$ over the training split used to construct $\widehat r_{\alpha, 0}$ (computed using \eqref{eq:rscp_cmdp_init}),
\begin{equation*}
    \begin{split}
     \| m_R^*(\widehat r_{\alpha, 0}, \cdot) - m_R^*(r_\alpha, \cdot)\|_2 \leq & \frac{\mathfrak C_{n,\delta}}{\sqrt{n}} + \frac{\mathcal{C}}{n_2} + \mathcal{C} a_m,
    \end{split}
\end{equation*}
where $a_m$ has the same definition as in \Cref{thm:influence_old_cmdp}, $\mathfrak C_{n,\delta} = \mathcal{C} \max\{\log n/\sqrt n, \log(1/\delta)/\sqrt{n}, 1/\sqrt{\delta}\}$, and $\mathcal C$ denotes a universal constant.
\end{prop}
The proof of Proposition~\ref{prop:rscp_init_error} is provided in \Cref{app:proof_bound_init_cmdp}. In table~\ref{tab:method_compare_cmdp} we show the comparison of the coverage guarantees derived in this paper with that of the earlier methods under the assumptions and regularity conditions of the theorems and propositions proved in this section. For ease of understanding we let $n_2 = O(n)$ in the $\rscp$ algorithm. Three different scenarios are considered: $(i)$ correctly specified $\pi^*(\cdot)$ and misspecified $m_R^*(\cdot, \cdot)$ ($a_\pi = 0, a_m >0$) $(ii)$ misspecified $\pi^*(\cdot)$ and correctly specified $m_R^*(\cdot, \cdot)$ ($a_\pi > 0, a_m = 0$) $(iii)$ correctly specified $\pi^*(\cdot)$ and correctly specified $m_R^*(\cdot, \cdot)$ ($a_\pi = 0, a_m = 0$). Since the nuisance functions are generally non-parametrically estimated, the nuisance estimation errors $a_\pi, a_m$ are often much larger than $O(n^{-1/2})$. We observe from table~\ref{tab:method_compare_cmdp} that $\rsp$ and $\rscp$ provide more robust coverage guarantees than WCP and conformal risk control across different scenarios and the slacks in coverage of $\rsp$ (and $\rscp$) do not explode if one of the nuisance functions is not correctly estimated. We also observe that the proposed $\rscp$ algorithm removes the additional $O(n^{-3/4})$ slack of the $\rsp$ algorithm when both of the nuisance functions are correctly specified. 
\begin{table}[h]
\centering
\resizebox{\textwidth}{!}{  
\begin{tabular}{l|c|c|l}
\toprule
\textbf{Method} & $a_\pi = 0, a_m >0$ & $a_\pi > 0, a_m = 0$ & $a_\pi = 0, a_m = 0$ \\
\midrule
$\begin{aligned}
  &\mbox{WCP} \\
  &(\mbox{\cite{tibshirani2019conformal}})
\end{aligned}$ & $0$ & $O(a_\pi)$ & 0\\
\hline 
$\begin{aligned}
    &\mbox{Conformal risk control (CRC)} \\
    &\mbox{\cite{angelopoulos2022conformal}}
\end{aligned}$& $O(a_m)$  & $0$ & $0$  \\
\hline 
$\begin{aligned}
    &\rsp (\mbox{ Coverage guarantee } \\
    &\mbox{proved in \cite{yang2024doubly}})
\end{aligned}$ & $O(n^{-1/2})$ & $O(n^{-1/2})$ & $O(n^{-1/2})$ \\
\hline 
$\begin{aligned}
    &\rsp (\mbox{ Coverage guarantee } \\
    &\mbox{proved in this work})
\end{aligned}$ & $O(n^{-3/4} + a_mn^{-1/2})$ & $O(n^{-3/4})$ & $O(n^{-3/4})$ \\
\hline 
$\begin{aligned}
    &\rscp (\mbox{ Coverage guarantee } \\
    &\mbox{proved in this work})
\end{aligned}$ & $0$ & $O(a_\pi n^{-1/2})$ & $0$ \\
\bottomrule
\end{tabular}
}  
\caption{Comparison of the coverage slack (departure from the target coverage level of $1-\alpha$) of $\rsp$ and $\rscp$ with prior works in the classical missing data problem. Here $a_\pi = \|\widehat \pi - \pi^*\|_2$ and $a_m = \sup_{\theta}\| \widehat m_R(\theta, \cdot) - m_R^*(\theta, \cdot)\|_2$ respectively denote the errors in estimating the nuisance functions.}
\label{tab:method_compare_cmdp}
\end{table}

It is useful to clarify the interpretation of this comparison. In the
classical missing-data problem, the $\rsp$ procedure is precisely the
efficient-influence-function estimating-equation construction underlying
the doubly robust method of \cite{yang2024doubly}; our contribution here
is to place that construction inside the general $\rsp$ framework and to
sharpen its coverage analysis. The $\rscp$ procedure, by contrast,
changes the calibration mechanism: after freezing the
threshold-dependent nuisance function, it uses a conformal risk control
correction. The substantive gain is that this construction can remove
the empirical-process remainder appearing in $\rsp$ and can yield exact
finite-sample lower validity when the relevant nuisance components are
correctly specified.

The advantage of $\rscp$ is most pronounced when the
multiply robust coverage slack is small, for example, when one of the nuisance components is accurately estimated. The numerical comparisons below include WCP, ordinary CRC,
$\rsp$ (the Yang et al.\ construction in this special case), and
$\rscp$. To illustrate the performance of these methods, we consider a simple simulation setting. We generate samples from the DGP,
\[
X_1, X_2, \varepsilon \stackrel{iid}{\sim} \mathcal N(0,1), \qquad
R = a_0 + a_1 X_1 + a_2 X_2 + \varepsilon,
\]
with \((a_0,a_1,a_2)=(20,70,-50)\). The missingness is generated as,
\[
\mathbb P(C = 0 \mid X) =p^*(X), \qquad
p^*(X)=0.2+0.6\,\mathrm{logit}^{-1}(b_0+b_1X_1+b_2X_2),
\]
with \((b_0,b_1,b_2)=(0.4,0.6,-0.6)\), so \(p^*(X)\in[0.2,0.8]\). Hence $\pi^*(X)={p^*(X)}/({1-p^*(X)}) \in [0.25,4]$. In order to demonstrate the advantage of the doubly robust methods $\rsp$ and $\rscp$, we consider the following adversarial estimates of $\pi^*$ and $m_R^*$.
\begin{equation*}
    \begin{split}
       \widehat\pi(X)={\hat p(X)}/({1-\hat p(X)}) \quad & \mbox{where} \quad \hat p(X)=0.2+0.6\,\mathrm{logit}^{-1}(b_0-b_1X_1-b_2X_2), \\
       \widehat m_R(\theta, X) = \Phi(\theta) \quad & \mbox{for all $\theta \in \mathbb R$ where $\Phi(\cdot)$ is the cdf of normal distribution}.
    \end{split}
\end{equation*}
As in \Cref{tab:method_compare_cmdp}, we consider three different scenarios, depending on whether the nuisance functions $\pi^*(\cdot)$ and $m_R^*(\cdot, \cdot)$ are estimated correctly: $(1)$ $a_\pi = 0, a_m > 0$ (we use the true nuisance function $\pi^*$ and the estimated nuisance function $\widehat m_R$), $(2)$ $a_\pi > 0, a_m = 0$ (we use the estimated nuisance function $\widehat \pi$ and the true nuisance function $m_R^*$), and $(3)$ $a_\pi = 0, a_m = 0$ (we use the true nuisance function $ \pi^*$ and the true nuisance function $m_R^*$). We use \(\alpha=0.1\), \(n=100\), and compute the coverage over $500$ $\{C = 0\}$ test points generated from the above DGP. The experiment is repeated \(100\) times to compute the mean and standard deviation of the coverage of the prediction sets. The initial estimate $\widehat r_{\alpha, 0}$ required for the $\rscp$ algorithm is computed using a seperate training data set of size $n_2 = 100$.

\begin{table}[h]
\centering
\begin{tabular}{lccc}
\hline
Method & \(a_\pi = 0, a_m > 0\) & \(a_\pi > 0, a_m = 0\) & \(a_\pi = 0, a_m = 0\) \\
\hline
WCP  & \(0.879 \pm 0.070\) & \(0.658 \pm 0.075\) & \(0.879 \pm 0.070\) \\
CRC  & \(0.337 \pm 0.023\) & \(0.896 \pm 0.045\) & \(0.896 \pm 0.045\) \\
$\rsp$  & \(0.876 \pm 0.079\) & \(0.895 \pm 0.046\) & \(0.886 \pm 0.049\) \\
$\rscp$ & \(0.932 \pm 0.066\) & \(0.955 \pm 0.066\) & \(0.948 \pm 0.046\) \\
\hline
\end{tabular}
\caption{Coverage (\(\text{mean} \pm \text{sd}\)) for WCP, CRC, $\rsp$, and $\rscp$ across three different settings.}
\label{tab:cov_compare_simul}
\end{table}

We observe from table~\ref{tab:cov_compare_simul} that the performance of WCP degrades strongly when (\(\pi \neq \hat\pi, m_R=m_R^*\)). This aligns with the sensitivity of WCP to \(\pi\)-misspecification. CRC behaves as expected; poor when \(m_R\) is misspecified and near target when \(m_R\) is correct. Both $\rsp$ and $\rscp$ are comparatively robust across the three different settings. We observe that although $\rscp$ maintains valid coverage ($\geq 0.9$) across all the settings, it is conservative (primarily because of the small size of the calibration data set).

\section{Conformal prediction under general missingness patterns}
\label{sec:gen_missingness}     
Having introduced $\rsp$ and $\rscp$ for the classical missing data problem, we now describe these methods in the general problem set-up. We denote the full data vector by $\tilde O = (X, Z, Y) \sim P_{\mathrm{full}}$ where $X$ denotes baseline covariates of primary scientific interest, $Y$ denotes the target outcome of interest, while $Z$ are auxiliary variables, which are not necessarily of scientific interest; however, they can be leveraged to account for missing outcome. Most often, the outcome of interest $Y$ would be real-valued, but our framework allows for more general outcomes of interest; see Remark~\ref{rem:nc_score}. 

The observed data vector is denoted by $O = (C, \mathfrak{G}_C(\tilde{O})) \sim P_{\mathrm{obs}}$ where $C$ is a random variable that denotes the coarsening level and takes value in $\mathcal K \coloneqq\{0, 1, \cdots, L, \infty\}$ and $\mathfrak G_k(\tilde{O})$ denotes a subset of the full data vector $\tilde{O}$ for all $k \in \mathcal K$. We designate $\{C = 0\}$ and $\{C = \infty\}$ to denote observing baseline covariates $X$ and observing full data $\tilde{O}$, respectively. The baseline covariates $X$ are assumed to be fully observed throughout i.e.,\ $X \subset \mathfrak{G}_k(\tilde{O})$ for $k \in \mathcal K$. The auxiliary covariates $Z$ and the target outcome $Y$ may be fully observed, ``partially" observed, or missing. \Cref{tab:missing_data_type} shows how this setting includes the classical, monotone and non-monotone missing data framework. From this, it is also clear that for a given coarsening $\mathfrak G_k(\tilde{O})$ ($k \in \mathcal{K}\setminus \{\infty\}$) of the data, there can be multiple values of $\tilde{O}$. In other words, $\tilde{O}\mapsto \mathfrak{G}_k(\tilde{O})$ is a many-to-one map, except when $k = \infty$. 

\begin{table}[h]
\centering
\resizebox{\textwidth}{!}{  
\begin{tabular}{l|c|c|l}
\toprule
\textbf{Missing Data} & \textbf{Full Data ($\tilde O$)} & \textbf{Coarsening Levels} ($\mathcal{K}$) & \textbf{Coarsened Data ($\{\mathfrak G_k(\tilde{O})\}_{k \in \mathcal K}$)} \\
\midrule
$\begin{aligned}
  &\mbox{Classical} \\
  &(\mbox{Covariate-Shift})
\end{aligned}$ & $(X,Y)$ & $\{0, \infty\}$ & $
   \begin{cases}
     X& \mbox{ if } k = 0,\\
  (X, Y) & \mbox{ if } k = \infty.
\end{cases} $ \\
\hline 
Monotone & $(X, Z_1,\cdots,Z_D, Y)$  & $\{0,1,\cdots,D, \infty\}$ & $\begin{cases}
     X& \mbox{ if } k = 0,\\
  (X, Z^1,\cdots,Z^k) & \mbox{ if } k \in \{1,\cdots, D\}, \\
  (X, Z^1,\cdots,Z^D, Y) & \mbox{ if } k = \infty.
\end{cases} $  \\
\hline 
Non-Monotone\footnotemark[2] & $(X, Z_1,\cdots,Z_D, Y)$ & $\begin{aligned}
    & \{0,1,\cdots,K, \infty\} \\
    (1 &\leq K \leq 2^{D+1} -2)
\end{aligned}$ & $\begin{aligned}
& \begin{cases}
     X& \mbox{ if } k = 0,\\
     \mathscr S_k & \mbox{ if } k \in \{1, \cdots, K\}, \\
  (X, Z^1,\cdots,Z^D, Y) & \mbox{ if } k = \infty.
\end{cases}\\
&\mbox{where the sets $\mathscr S_1,\cdots,\mathscr S_K$ satisfy}, \\
&X \subset \mathscr S_k \subset (X, Z^1,\cdots,Z^D, Y) \mbox{ $\forall$ } k \in [K].
\end{aligned}$ \\
\bottomrule
\end{tabular}
}  
\caption{Different types of missing data in the calibration sample and the corresponding coarsening notations. Here $X$ denotes baseline covariates, $Z^1,\cdots,Z^D$ are auxiliary variables, and $Y$ denotes the response variable.}
\label{tab:missing_data_type}
\end{table}
\footnotetext[2]{Note that Non-Monotone missing data framework is a more general regime and comprises Monotone Missing Data as one of the special cases. In particular, if we set $K= D$, $\mathscr S_k = (X, Z_1, \cdots, Z_k)$ for all $k \in [D]$, the setup reduces to the monotone missing data framework.}
Our overarching objective is to predict $Y$ for a future sample for which $X$ is observed using a \emph{well-calibrated} prediction algorithm, trained on coarsened data $(C, \mathfrak G_C(\tilde O))$ of $\tilde O$. The observed data set $\mathcal{D} =  \mathcal{D}^{tr} \cup\mathcal{D}^{cal}$ is a disjoint union of two subsets: the training subset $ \mathcal{D}^{tr} = \{O_1^t, \cdots, O_m^t \}$ and the calibration subset $ \mathcal{D}^{cal} = \{O_1^{cal}, \cdots, O_n^{cal} \} = \{O_1, \cdots, O_n\}$ (we suppress the superscript $cal$ wherever it is clear from the context). All the $(m+n)$ data points $ \{O_1^t, \cdots, O_m^t \} \cup \{O_1, \cdots, O_n\}$ in $\mathcal{D}$ are iid from the distribution $\mathcal{P}_O$. Given a new coarsened observation $O_{n+1} \sim \mathcal{P}_O$, our objective is to obtain a valid $1 - \alpha$ prediction interval for the corresponding $Y_{n+1}$. We describe the Robust Split Prediction ($\rsp$) method in \Cref{alg:split_robust} and the Robust Split Conformal Prediction ($\rscp$) method in \Cref{algg:split_robust_three_split}. 
\begin{algorithm}
    \caption{Robust Split Prediction ($\rsp$)}
    \label{alg:split_robust}
    \KwIn{Observed data set: $\mathcal{D}  = \{O_1^t, \cdots, O_m^t, O_1, \cdots, O_n\} = \mathcal{D}^{tr} \cup \mathcal{D}^{cal}$, Confidence level: $1-\alpha$}
    \KwOut{A valid prediction set $\widehat{\mathcal C}_{n,\alpha}$}
    \begin{algorithmic}[1]
   \State Obtain a function $f(O, \theta, \{\eta_i^*\}_{i = 1}^K)$ of the observed data with nuisance functions $\{\eta_i^*\}_{i = 1}^K$ such that $\mathbb{P}(Y \leq \theta) = \mathbb{E}[f(O, \theta, \{\eta_i^*\}_{i = 1}^K) ]$ for all $\theta \in \mathbb{R}$.
   \State Obtain estimates $\{\widehat \eta_i \}_{i = 1}^K$ of the nuisance functions $\{\eta_i^*\}_{i = 1}^K$ from the sub-sample $\mathcal{D}^{tr}$. 
  \State  Define $\widehat r_{\alpha}$ using the calibration set $\mathcal{D}^{cal}$,  
  \[
   \widehat  r_{\alpha} = \inf\left\{\theta \in \mathbb{R}: \frac{\sum_{i = 1}^n f(O_i, \theta,  \{\widehat \eta_i \}_{i = 1}^K)}{n}  \geq 1 - \alpha\right\}.
  \]
  \State Return the prediction set, 
  \[
  \widehat{\mathcal C}_{n,\alpha}  = \{y \in \mathbb{R}: y \leq \widehat r_{\alpha}\}. 
  \]
 \end{algorithmic}
\end{algorithm}

\begin{algorithm}
    \caption{Robust Split Conformal Prediction ($\rscp$)}
    \label{algg:split_robust_three_split}
    \KwIn{Observed data set: $\mathcal{D}  = \{O_1^{t_1}, \cdots, O_{n_1}^{t_1},O_1^{t_2}, \cdots, O_{n_2}^{t_2} , O_1, \cdots, O_n\} = \mathcal{D}^{tr}_1 \cup \mathcal{D}^{tr}_2  \cup \mathcal{D}^{cal}$, Confidence level: $1-\alpha$}
    \KwOut{A valid prediction set $\widehat{\mathcal C}_{n,\alpha}$}
    \begin{algorithmic}[1]
  \State Obtain a function $f(O, \theta, \{\eta_{1i}^*(\theta)\}_{i = 1}^{K_1}, \{\eta_{2i}^*\}_{i = 1}^{K_2})$ of the observed data with nuisance functions $\{\eta_{1i}^*(\theta)\}_{i = 1}^{K_1}, \{\eta_{2i}^*\}_{i = 1}^{K_2}$ such that $\mathbb{P}(Y \leq \theta) = \mathbb{E}[f(O, \theta, \{\eta_{1i}^*(\theta)\}_{i = 1}^{K_1}, \{\eta_{2i}^*\}_{i = 1}^{K_2}) ]$ for all $\theta \in \mathbb{R}$.
   \State Obtain estimates $\{\widehat \eta_{1i}(\theta)\}_{i = 1}^{K_1}, \{\widehat \eta_{2i}\}_{i = 1}^{K_2}$ of the nuisance functions $\{\eta_{1i}^*(\theta)\}_{i = 1}^{K_1}, \{\eta_{2i}^*\}_{i = 1}^{K_2}$ from the sub-sample $\mathcal{D}^{tr}_1$. 
   \State Define $\widehat r_{\alpha,0}$ using the sub-sample $\mathcal{D}^{tr}_2$,
   \begin{equation*}
      \begin{split}
      \widehat  r_{\alpha,0} = & \inf\left\{\theta \in \mathbb{R}: \frac{\sum_{i = 1}^{n_2} f(O_i^{t_2}, \theta,  \{\widehat \eta_{1i}(\theta)\}_{i = 1}^{K_1}, \{\widehat \eta_{2i}\}_{i = 1}^{K_2})}{n_2}  \geq 1 - \alpha\right\}.
      \end{split}
  \end{equation*}
  \State  Define $\widehat r_{\alpha}$ using the calibration set $\mathcal{D}^{cal}$,  
  \begin{equation*}
      \begin{split}
      \widehat  r_{\alpha} = &\inf\Bigg\{\theta \in \mathbb{R}: \frac{\sum_{i = 1}^n f(O_i, \theta, \{\widehat \eta_{1i}(\widehat r_{\alpha, 0})\}_{i = 1}^{K_1}, \{\widehat \eta_{2i}\}_{i = 1}^{K_2})}{n+1} +  \frac{\inf_{o} f(o, \theta,  \{\widehat \eta_{1i}(\widehat r_{\alpha, 0})\}_{i = 1}^{K_1}, \{\widehat \eta_{2i}\}_{i = 1}^{K_2})}{n+1}  \geq 1 - \alpha \Bigg\}.
      \end{split}
  \end{equation*}
  \State Return the prediction set, 
  \[
  \widehat{\mathcal C}_{n,\alpha}  = \{y \in \mathbb{R}: y \leq \widehat r_{\alpha}\}. 
  \]
 \end{algorithmic}
\end{algorithm}
\begin{remark}[Choice of non-conformity score]
\label{rem:nc_score}
We use the non-conformity score $R(X,Y)=Y$ in the general algorithms
only to simplify notation. The theory applies to any real-valued
non-conformity score $R(X,Y)$ for which the coverage probability can be identified using observed-data distribution i.e.\ there exists an estimating function $f_R(\cdots)$ such that $\mathbb P
(
R(X,Y)\leq\theta
)
=
\mathbb E
[
f_R(O,\theta,\eta_R^*)
]$. The corresponding prediction set is $\widehat{\mathcal C}_{n,\alpha}(x)
=
\{
y:
R(x,y)\leq\widehat r_\alpha
\}.$ Examples include the absolute-residual score $R(X,Y)
=
|Y-\widehat\mu(X)|$ \citep{lei2018distribution}, the conformalized-quantile score $R(X,Y)
=
\max
\{
\widehat q_{\alpha/2}(X)-Y,\,
Y-\widehat q_{1-\alpha/2}(X)
\}$ \citep{romano2019conformalized}, and the inverse-quantile scores used for discrete outcomes in \cite{romano2020classification}.

If the score function itself is learned from a data set independent of the calibration samples and of the data set where the nuisance functions are estimated, then throughout the
analysis we may condition on this sample from which the score function is learned. Conditional on
this sample, $R$ is fixed, and all of the coverage arguments apply
unchanged to the distribution of $R(X,Y)$. In particular, the double or
multiple robustness of the resulting prediction set is inherited from
the efficient influence function for the distribution of
$R(X,Y)$: changing the non-conformity score changes the corresponding
outcome-regression nuisance functions, but does not alter the
stagewise mixed-bias structure. If the same observations are used both to learn $R$ and to calibrate
the prediction set, then we must account for this dependency.
Sample splitting, as used in our theoretical development and
experiments, avoids this issue. Finally, in \Cref{alg:split_robust} and \Cref{algg:split_robust_three_split}, the outcome $Y$ need not be
univariate. The same algorithms apply to a multivariate outcome
provided that the chosen map $(X,Y)\mapsto R(X,Y)$ is real valued.
\end{remark}


Before stating the general coverage guarantees, we introduce the necessary assumptions.  Throughout this section we use $R(X,Y)=Y$ as the non-conformity score, as in
\Cref{alg:split_robust} and \Cref{algg:split_robust_three_split}. We write $Ph=\mathbb E[h(O)]$ (where the expectation is taken conditional on the training data set $\mathcal D^{\mathrm{tr}}$) and $P_nh = (1/n)  \sum_{i=1}^n h(O_i) $. We also define the following,
\[
    F_Y(\theta)=\mathbb P(Y\leq\theta),
    \qquad
    r_\alpha=\inf\{\theta:F_Y(\theta)\geq1-\alpha\}.
\]

\begin{assumption}[Identification]
\label{assump:gen_identification} 
We assume that $r_\alpha$ satisfies $ F_Y(r_\alpha)=1-\alpha.$ Moreover there exist nuisance functions
$\eta^*=\{\eta_i^*\}_{i=1}^K$ such that, for every
$\theta\in\mathbb R$, $P f(O,\theta,\eta^*)=F_Y(\theta).$
\end{assumption}

The first condition ensures that
$r_\alpha$ is an exact $(1-\alpha)$-quantile.  It is satisfied, for
example, whenever the distribution of $Y$ is continuous at
$r_\alpha$.
The second condition is the identifying relation that connects the distribution of the target
outcome to the distribution of the observed data. 

\begin{assumption}[Error induced by nuisance estimation]
\label{assump:gen_nuisance}
For any collection of nuisance functions
$\widetilde\eta=\{\widetilde\eta_i\}_{i=1}^K$,
\[
    \sup_{\theta\in\mathbb R}
    \left|
        P f(O,\theta,\eta^*)
        -
        P f(O,\theta,\widetilde\eta)
    \right|
    \leq
    d(\eta^*,\widetilde\eta),
\]
where $d(\eta^*,\widetilde\eta)$ is a nonnegative measure of the
difference between the two collections of nuisance functions.
\end{assumption}

The quantity $d(\eta^*,\widetilde\eta)$ measures the bias in the
population estimating equation caused by replacing the true nuisance
functions by estimated ones.  In the applications considered in this
paper, $f$ is typically multiply robust and $d(\cdot, \cdot)$ is therefore a
higher-order quantity. For example, in the classical missing-data
problem, $d(\cdot, \cdot)$ is proportional to the product of the errors in estimating
the propensity-score and outcome-regression functions.

\begin{assumption}[Bounded variation in the threshold]
\label{assump:gen_bv}
There exist constants $B<\infty$ and $V<\infty$ such that, for every
collection of nuisance functions
$\widetilde\eta$ under consideration,
\[
    \sup_{\theta\in\mathbb R}
    |f(O,\theta,\widetilde\eta)|
    \leq B
    \qquad\mbox{almost surely},
\]
and the map $\theta\mapsto f(O,\theta,\widetilde\eta)$ is right-continuous and has total variation uniformly bounded by $V$, 
\[
    \sup_{M \geq 1}\sup_{\theta_0<\cdots<\theta_M}
    \sum_{j=1}^{M}
    \left|
        f(O,\theta_j,\widetilde\eta)
        -
        f(O,\theta_{j-1},\widetilde\eta)
    \right|
    \leq V
    \qquad\mbox{almost surely}.
\]
In addition, for every fixed $\theta\in\mathbb R$, we have $ \mathbb P(
        f(O,\cdot,\widetilde\eta)
        \mbox{ is discontinuous at }\theta
    )=0.$
\end{assumption}

The assumption of bounded variation is automatically satisfied for several standard missing data regimes. In particular, the function $f(O,\cdot,\widetilde\eta)$ for the classical and monotone missing data problems considered in this work can be written as a difference of two bounded monotone functions, and hence has bounded total variation.  


\begin{assumption}[Second-moment bound]
\label{assump:gen_local_variance}
For any collection of nuisance functions
$\widetilde\eta$, we define the centered function,
\[
    g(O,\theta,\widetilde\eta)
    =
    f(O,\theta,\widetilde\eta)
    -
    f(O,r_\alpha,\widetilde\eta).
\]
There exists a constant $K<\infty$ and a nonnegative quantity
$d'(\eta^*,\widetilde\eta)$ such that, for every
$\theta\in\mathbb R$,
\[
    P g^2(O,\theta,\widetilde\eta)
    \leq
    K
    \left|
        P g(O,\theta,\widetilde\eta)
    \right|
    +
    d'(\eta^*,\widetilde\eta).
\]
\end{assumption}

The above bound is a Bernstein-type localization condition, the
second moment of the increment becomes small when its population mean
is small.  The additional variance inflation term $d'(\cdot, \cdot)$ need not be equal to
$d(\cdot, \cdot)$.  In the classical missing-data problem, for example, the
corresponding calculation gives
$d'(\eta^*,\widehat \eta) =  a_\pi a_m+a_m^2$, while $d(\eta^*,\widehat \eta) = a_\pi a_m$ (where $a_\pi, a_m$ are the nuisance estimation errors defined in \Cref{thm:influence_old_cmdp}). The following theorem provides the coverage guarantee for $\rsp$ under general missingness patterns.

\begin{theorem}[Coverage guarantee for $\rsp$]
\label{thm:gen_influ_func}
Suppose Assumptions~\ref{assump:gen_identification}--\ref{assump:gen_local_variance}
hold. Let $\widehat r_\alpha$ be the threshold returned by \Cref{alg:split_robust}. Then, conditionally on the training set $\mathcal D^{\mathrm{tr}}$,
\[
    \left|
        \mathbb P
        \left(
            Y_{n+1}\leq\widehat r_\alpha
        \right)
        -(1-\alpha)
    \right|
    \leq
    \Delta_{\rsp,n}, \qquad \mbox{where $\Delta_{\rsp,n}$ is defined as,}
\]
\begin{equation}
\label{eq:rsp_general_slack}
\begin{split}
    \Delta_{\rsp,n}
    =
    d_{\rsp}
    +
    \mathfrak C_n
    \left(
        \sqrt{d_{\rsp}}+n^{-1/4}
    \right)
    \left(
        n^{-1/2}
        +
        n^{-1/4}\sqrt{d'_{\rsp}}
    \right),
\end{split}
\end{equation}
where $\mathfrak C_n=\mathfrak C\log(en)$ and $\mathfrak C$ is a constant depending only on the constants in
Assumptions~\ref{assump:gen_bv} and
\ref{assump:gen_local_variance}, $d_{\rsp}
    =
    d(\eta^*,\widehat\eta)$, $ d'_{\rsp}
    =
    d'(\eta^*,\widehat\eta)$. 

\end{theorem}

\begin{remark}[Interpretation of the coverage slack of $\rsp$ algorithm]
\label{rem:rsp_general_slack}
The leading term $d_{\rsp}$ is the population bias induced by nuisance
estimation.  The remaining term is the finite-sample error arising from using the sample moment equation $P_n f(O, \theta, \widehat \eta) = 0$ to compute the threshold $\widehat r_\alpha$.  In particular, if
$d_{\rsp}=d'_{\rsp}=0$, then $\Delta_{\rsp,n}
    =
    O(
        {\log n}/{n^{3/4}}
    ).$
More generally, if $d_{\rsp}$ and $d'_{\rsp}$ are second-order
multiply robust quantities, then the entire coverage error inherits the corresponding robustness property.
\end{remark}

\begin{remark}[Conditioning on a coarsening pattern]
\label{rem:general_pattern_conditional}
Theorem~\ref{thm:gen_influ_func} targets the marginal distribution of the outcome $Y$.  Coverage guarantee conditional on a particular missingness or coarsening pattern follows
from exactly the same argument.  For example, for an event
$\mathcal A\subseteq\mathcal K$ with
$\mathbb P(C\in\mathcal A)>0$, suppose one can construct a function $f_{\mathcal A}$ of observed data satisfying,
\[
    \mathbb P(Y\leq\theta\mid C\in\mathcal A) = \mathfrak C  P f_{\mathcal A}(O,\theta,\eta^*),
\]
where $\mathfrak C$ is an arbitrary constant. Then Theorem~\ref{thm:gen_influ_func} provides the lower and upper bounds on $\mathbb P
    (
        Y_{n+1}\leq\widehat r_\alpha
        \mid
        C_{n+1}\in\mathcal A
    )$ provided the required assumptions hold for the function $f_{\mathcal A}$. Therefore no separate
proof is required for coverage conditional on a missingness pattern;
only the corresponding identification formula changes.
\end{remark}



We next state the general coverage guarantee for $\rscp$. Recall that
the nuisance functions in \Cref{algg:split_robust_three_split} are
partitioned into $\{\eta_{1i}^*(\theta)\}_{i=1}^{K_1}$ and $\{\eta_{2i}^*\}_{i=1}^{K_2}$, where the first collection depends on the threshold $\theta$, whereas
the second collection does not.

Throughout the following assumptions, theorem, and proof, all
probabilities and expectations are conditional on the two training
samples $\mathcal T
    :=
    \mathcal D_1^{\mathrm{tr}}
    \cup
    \mathcal D_2^{\mathrm{tr}}$. Consequently, the nuisance estimates
$\{\widehat\eta_{1i}(\theta)\}_{i=1}^{K_1}$,
$\{\widehat\eta_{2i}\}_{i=1}^{K_2}$, and the initial estimate
$\widehat r_{\alpha,0}$ are treated as fixed. For $\theta,s\in\mathbb R$, define,
\[
\begin{split}
    f^*_\theta(O)
    &:=
    f\left(
        O,\theta,
        \{\eta_{1i}^*(\theta)\}_{i=1}^{K_1},
        \{\eta_{2i}^*\}_{i=1}^{K_2}
    \right),\qquad
    \widehat f_{\theta,s}(O)
    :=
    f\left(
        O,\theta,
        \{\widehat\eta_{1i}(s)\}_{i=1}^{K_1},
        \{\widehat\eta_{2i}\}_{i=1}^{K_2}
    \right),
\end{split}
\]

Using these notations, the final calibration step of $\rscp$ can be re-written as
\begin{equation}
\label{eq:rscp_general_centered}
    \widehat r_\alpha
    =
    \inf\left\{
        \theta\in\mathbb R:
        \frac{
            \sum_{i=1}^n
            \phi_{\theta,\widehat r_{\alpha,0}}(O_i)
            +
            \inf_o
            \phi_{\theta,\widehat r_{\alpha,0}}(o)
        }{n+1}
        \geq0
    \right\}, \quad  \phi_{\theta,s}(O)
    :=
    \widehat f_{\theta,s}(O)-(1-\alpha).
\end{equation}

We require the following additional assumptions.

\begin{assumption}[Error induced by nuisance estimation and the initial threshold estimate]
\label{assump:rscp_anchored_bias}
There exists a nonnegative function $s \mapsto d_{\mathrm{init}}(s,r_\alpha)$ satisfying $d_{\mathrm{init}}(r_\alpha,r_\alpha)=0$ such that, for every $s,\theta\in\mathbb R$,
\begin{equation}
\label{eq:rscp_anchored_bias_bound}
\left|
    P f^*_\theta
    -
    P\widehat f_{\theta,s}
\right|
\leq
d_{\rscp}
+
d_{\mathrm{init}}(s,r_\alpha), \qquad \mbox{where $d_{\rscp}$ is defined as,}
\end{equation}
\begin{equation}
\label{eq:rscp_population_bias}
\begin{split}
    d_{\rscp}
    :=
    \sup_{\theta\in\mathbb R}
    d\Big(
        &\{\eta_{1i}^*(\theta)\}_{i=1}^{K_1},
        \{\eta_{2i}^*\}_{i=1}^{K_2},\{\widehat\eta_{1i}(\theta)\}_{i=1}^{K_1},
        \{\widehat\eta_{2i}\}_{i=1}^{K_2}
    \Big),
\end{split}
\end{equation}
where $d(\cdot, \cdot)$ is the same measure as defined in Assumption~\ref{assump:gen_nuisance}.
\end{assumption}
Assumption~\ref{assump:rscp_anchored_bias} separates two fundamentally
different sources of error. The quantity $d_{\rscp}$ quantifies the usual difference in the population moments because of nuisance estimation.  In contrast,
$d_{\mathrm{init}}(s,r_\alpha)$ measures the additional error caused
by evaluating the threshold-dependent nuisance functions at $s$ instead of the oracle threshold $r_\alpha$.  

\begin{assumption}[Monotonicity after freezing the nuisance functions]
\label{assump:rscp_monotone}
For every fixed $s\in\mathbb R$ and every observed data value $o$, the
map $\theta
    \longmapsto
    \widehat f_{\theta,s}(o)$ is non-decreasing and right-continuous.
\end{assumption}
The full estimating function need not be monotone when all nuisance
functions are evaluated at the varying threshold $\theta$.  The
purpose of the preliminary estimate $\widehat r_{\alpha,0}$ is to
freeze the threshold-dependent nuisance functions.  Once they are
frozen, the remaining dependence on $\theta$ is typically through a
term such as $\mathbf 1\{Y\leq\theta\}$, making the estimating
function monotone.  This monotonicity is what allows the conformal
exchangeability argument to go through and establish the validity of the method.

\begin{assumption}[Range and jump control]
\label{assump:rscp_range}
There exist constants $B_{\rscp}<\infty$ and
$J_{\rscp}<\infty$ such that, for every $\theta,s\in\mathbb R$,
\[
    \sup_o \widehat f_{\theta,s}(o)
    -
    \inf_o \widehat f_{\theta,s}(o)
    \leq
    B_{\rscp}.
\]
Moreover, every jump of
$\theta\mapsto\widehat f_{\theta,s}(o)$ has magnitude at most
$J_{\rscp}$, and, almost surely, two independent observations do not
have a common discontinuity.  
\end{assumption}
The range bound guarantees that the additional conformal term
$\inf_o\widehat f_{\theta,s}(o)$ is sufficiently conservative to
account for the unobserved $(n+1)$-st observation.  The jump bound is
only required for the upper coverage bound and controls the amount by
which the full-sample estimating equation can overshoot its target at
the first crossing.  For most missing-data settings,
these conditions follow from bounded inverse-probability weights and
continuity of the outcome. An immediate consequence of Assumption~\ref{assump:rscp_range} is the following bound,
\begin{equation}
\label{eq:rscp_empirical_jump}
    \sup_{\theta}
    J\left(
        \frac{1}{n+1}
        \sum_{i=1}^{n+1}
        \widehat f_{\theta,s}(O_i),
        \theta
    \right)
    \leq
    \frac{J_{\rscp}}{n+1}
    \qquad\mbox{almost surely}.
\end{equation}
\begin{theorem}[Coverage guarantee for $\rscp$]
\label{thm:gen_three_split}
Suppose Assumptions~\ref{assump:gen_identification},
\ref{assump:rscp_anchored_bias},
\ref{assump:rscp_monotone}, and
\ref{assump:rscp_range} hold. Let the threshold $\widehat r_\alpha$, computed by $\rscp$ algorithm, be defined by
\eqref{eq:rscp_general_centered}. Then we have the following coverage guarantee conditional on both the training splits,
\begin{equation}
\label{eq:rscp_lower_general}
    \mathbb P
    \left(
        Y_{n+1}\leq\widehat r_\alpha
    \right)
    \geq
    1-\alpha
    -
    d_{\rscp}
    -
    \Delta_{\mathrm{init}},
\end{equation}
where $ \Delta_{\mathrm{init}}
    :=
    d_{\mathrm{init}}
    (\widehat r_{\alpha,0},r_\alpha)$. Moreover, we have the following upper bound on the coverage probability,
\begin{equation}
\label{eq:rscp_upper_general}
    \mathbb P
    \left(
        Y_{n+1}\leq\widehat r_\alpha
    \right)
    \leq
    1-\alpha
    +
    d_{\rscp}
    +
    \Delta_{\mathrm{init}}
    +
    \frac{B_{\rscp}+J_{\rscp}}{n+1}.
\end{equation}
\end{theorem}

Theorem~\ref{thm:gen_three_split} shows that the coverage slack of
$\rscp$ has two principal components, $ d_{\rscp}$ and $d_{\mathrm{init}}
    (\widehat r_{\alpha,0},r_\alpha)$. The first term is the usual bias in population moment resulting from nuisance
estimation. The second term measures only how much the threshold-dependent nuisance functions change between the preliminary
threshold $\widehat r_{\alpha,0}$ and the oracle threshold $r_\alpha$. For the classical missing-data problem, after normalizing by
$\mathbb P(C=0)$, these two quantities reduce to,
\[
    d_{\rscp}
    =
    \frac{
        a_\pi a_m
    }{
        \mathbb P(C=0)
    }, \qquad  d_{\mathrm{init}}
    (\widehat r_{\alpha,0},r_\alpha)
    =
    \frac{
        a_\pi
        \left\|
            m_R^*(\widehat r_{\alpha,0},\cdot)
            -
            m_R^*(r_\alpha,\cdot)
        \right\|_2
    }{
        \mathbb P(C=0)
    }.
\]
Thus the coverage guarantee obtained in \Cref{thm:validity_proof} for the classical missing data problem becomes a special case of Theorem~\ref{thm:gen_three_split}. In the following sections, we show that the same phenomenon occurs under monotone missingness. The bias in the population moment for the monotone missing data setting generally has the following form,
\[
    d_{\rscp}
    \lesssim
    \sum_{j=1}^{D+1}
    A_j E_j,
\]
where $A_j$ is the error in estimating propensity score associated with stage
$j$ and $E_j
    =
    \sup_{\theta}
    \|
        \widehat m_j(\theta,\cdot)
        -
        m_j^*(\theta,\cdot)
    \|_2$ is the corresponding error in estimating outcome-regression function.  The error term $ d_{\mathrm{init}}
    (\widehat r_{\alpha,0},r_\alpha)$ has the form,
\[
    d_{\mathrm{init}}
    (\widehat r_{\alpha,0},r_\alpha)
    \lesssim
    \sum_{j=1}^{D+1}
    A_j
    \left\|
        m_j^*(\widehat r_{\alpha,0},\cdot)
        -
        m_j^*(r_\alpha,\cdot)
    \right\|.
\]
If, locally around $r_\alpha$, $ \|
        m_j^*(s,\cdot)
        -
        m_j^*(r_\alpha,\cdot)
    \|_2
    \leq
    L_j |s-r_\alpha|$, then $d_{\mathrm{init}}
    (\widehat r_{\alpha,0},r_\alpha)
    \lesssim
    |\widehat r_{\alpha,0}-r_\alpha|
    \sum_{j=1}^{D+1} A_jL_j.$ Thus any consistency or rate result for the preliminary estimator
$\widehat r_{\alpha,0}$ translates directly into a coverage guarantee
for $\rscp$.

\begin{remark}[Comparison of $\rsp$ and $\rscp$]
\label{rem:rsp_rscp_comparison}
The coverage
slack for $\rsp$ contains, in addition to the usual nuisance estimation error $d_{\rsp}$,
the empirical-process term,
\[
    \mathfrak C_n
    \left(
        \sqrt{d_{\rsp}}+n^{-1/4}
    \right)
    \left(
        n^{-1/2}
        +
        n^{-1/4}\sqrt{d'_{\rsp}}
    \right).
\]
This term arises because $\rsp$ estimates the target quantile by
directly solving an empirical estimating equation. In contrast, $\rscp$ uses an independent preliminary estimate
$\widehat r_{\alpha,0}$ to freeze the threshold-dependent nuisance
functions in order to make the estimating function monotonic in the threshold and then adds the conformal correction in
\eqref{eq:rscp_general_centered}.  This allows us to use the usual exchangeability arguments directly.  As a result, the empirical-process term disappears entirely from the lower bound on the coverage. Its role is replaced by $ d_{\mathrm{init}}
    (\widehat r_{\alpha,0},r_\alpha),$ which is controlled by the accuracy of an independent preliminary
estimator. In classical and monotone missing data problems, $\rscp$ generally provides improved coverage guarantee (see for example, \Cref{tab:method_compare_cmdp}) as the initialization
error $ d_{\mathrm{init}}
    (\widehat r_{\alpha,0},r_\alpha)$ based on an independent data split can be controlled more sharply than the empirical-process remainder term appearing in the $\rsp$ bound. 

\end{remark}

\paragraph{How to choose the function $f$}The pivotal step in \Cref{alg:split_robust} is to find a function $f(O, \theta, \{\eta_i^*\}_{i = 1}^K)$ of the observed data such that $\mathbb{P}(Y \leq \theta) = \mathbb{E}[f(O, \theta, \{\eta_i^*\}_{i = 1}^K) ]$ for all $\theta \in \mathbb{R}$. Therefore a natural question that arises here is how to systematically choose this function $f$. We detail a specific method for construction below, 
\begin{enumerate}
    \item Derive the efficient influence function $\mathrm{IF}(O, \theta, \{ \eta_i^* \}_{i = 1}^K)$ corresponding to the functional $\mathbb{P}(Y \leq \theta)$ under a given semi-parametric model for which the identification is established. By definition of influence function, $\mathbb{E}[\mathrm{IF}(O, \theta, \{ \eta_i^* \}_{i = 1}^K)] = 0$. 
    \item Since $\mathrm{IF}(O, \theta, \{ \eta_i^* \}_{i = 1}^K) + \mathbb{P}(Y \leq \theta)$ is only a function of $\theta$, the observed data $O$, and the nuisance functions $\{ \eta_i^*\}_{i = 1}^K$, we set $f(O, \theta, \{ \eta_i^* \}_{i = 1}^K) = \mathrm{IF}(O, \theta, \{ \eta_i^* \}_{i = 1}^K) + \mathbb{P}(Y \leq \theta)$. Clearly we have $\mathbb{E}[f(O, \theta, \{\eta_i^*\}_{i = 1}^K) ] = \mathbb{P}(Y \leq \theta)$. 
\end{enumerate}

To illustrate the application of the general methods discussed in this section to specific problems, we demonstrate the application of $\rscp$ (\Cref{algg:split_robust_three_split}) to the monotone missing data problem in the following section. The coverage guarantee for $\rsp$ can be obtained analogously. 

\section{Application: The monotone missing data problem}
\label{sec:monotone}
Monotone missing data is a common occurrence in empirical studies of social and health sciences (see \cite{robins1994estimation, laan2003unified, tsiatis2006semiparametric, tchetgen2006statistical, sun2018inverse}). For example, suppose the complete data-vector is $\tilde O = (X, Z, Y)$, where the baseline covariates $X$ are always observed, the intermediate covariates $Z$ are observed only for a subset of individuals, and the target outcome $Y$ is observed only for a further subset of individuals who have complete data for both $X$ and $Z$. This is a two-stage monotone missing data problem which commonly occurs in longitudinal cohort studies when participants never return after dropping out from the study. We introduce some definitions and formalize the general multi-stage monotone missing data problem. We assume that $X$ denotes the baseline covariates, $Z^j$ denotes the covariates at the $j$-th stage for $1 \leq j \leq D$, and $Y$ denotes the potential outcome. Therefore $\tilde O = (X, Z^1, \cdots, Z^D, Y)$ denotes the complete information (full data) about a randomly drawn data point. Following the coarsening notation in \Cref{tab:missing_data_type}, the observed data is of the form $(C, \mathfrak G_C(\tilde O))$ where $C$ is a random variable that denotes the coarsening level and takes values in the set $\mathcal{K} = \{0, 1, \cdots, D, \infty\}$. As discussed in prior sections, $\mathfrak G_0(\tilde O) =X, \mathfrak G_k(\tilde O) = (X, Z_1,\cdots,Z_k)$ for all $k \in [D]$, and $\mathfrak G_{\infty}(\tilde O) = (X, Z_1, \cdots, Z_D, Y)$.
For example, we will later consider a study where the goal is to predict the future outcome at time $D+1$, for a new sample with observed baseline characteristics $X$.  The predictive algorithm is  calibrated using data from a longitudinal study where $j$ indexes study visits of subjects in the training sample at which covariates $Z^j$ are collected. Censoring due to participants' dropout is not uncommon in such a longitudinal setting, inducing the monotone coarsening process whereby $C = j-1$ for an individual who drops out by time $j$, in which case $Z^j \ldots Z^D,Y$ are missing. In monotone missing data problem we have $\mathfrak G_k(\tilde O) \subset \mathfrak G_{k+1}(\tilde O)$ for all $k \neq \infty$ i.e.\ $\mathfrak G_{k+1}(\tilde O) \mapsto \mathfrak G_k(\tilde O)$ is a many to one function for all $k \neq \infty$. In other words, $\mathfrak G_0(\tilde O)$ is the most coarsened data, followed by $\mathfrak G_1(\tilde O), \cdots$, and $\mathfrak G_{\infty}(\tilde O)$ is the complete data (no coarsening). Therefore, the observed calibration dataset $\mathcal{D}^{cal}$ comprises the $n$ i.i.d.\ data points, $\mathcal{D}^{cal} = \{O_1, \cdots, O_n\}$ where $O_i = (C_i, \mathfrak G_{C_i}(\tilde O_i))$ (here $C_i$ is the random variable denoting the coarsening level, and $\tilde O_i = (X_i, Z^1_i,\cdots,Z^D_i, Y_i)$ denotes the full information about the $i$-data point) for $i = 1,\cdots,n$. As in the covariate shift problem (the classical missing data setting), we make the missing at random (MAR) assumption. 
\begin{equation}
\label{eq:mar_monotone}
    \begin{split}
        \mathbb{P}(C = k | X, Z^1, \cdots,Z^D, Y) =  \omega(k, \mathfrak G_k((X, Z^1, \cdots,Z^D, Y))) \quad \mbox{for all } k \in \{0, 1, \cdots, D, \infty\}. 
    \end{split}
\end{equation}
This notation is due to \cite{robins1994estimation} and then adopted by \cite{tsiatis2006semiparametric} in Chapter $7.1$. Although prior works have addressed identification and inference under monotone missingness, the problem of constructing valid outcome prediction sets in the multi-stage monotone missing data setting has not been investigated.  It is important to note that the time-varying data $Z^j$ are themselves not of scientific interest in building the predictive algorithm, as they typically will not be available in a new sample for which we wish to predict the outcome given observed baseline characteristics $X$, e.g. a new patient being seen by a cardiologist wishing to predict cardiovascular disease risk. The role of $Z^j$ is primarily to account for informative dropout in the calibration sample, which may be a cohort study available to the analyst. To facilitate understanding, we apply $\rscp$ and analyze the coverage of the resulting prediction set in the two-stage monotone missing data setting, and defer the results for the general monotone missing data problem to \Cref{app:multi_stage_monotone}. For the two stage monotone missingness problem, the calibration data is $\mathcal{D}^{cal} = \{O_1, \cdots, O_n\}$ where $O_i = (C_i, \mathfrak G_{C_i}(\tilde O_i))$ for $i \in [n]$. For the $i$-th data point, the complete data vector is $\tilde O_i = (X_i, Z^1_i, Y_i) = (X_i, Z_i, Y_i)$ (we suppress the superscript $1$ for the two stage monotone missingness problem to enhance readability) and the coarsening random variable $C_i$ takes value in the set $\{0, 1, \infty\}$ where $C_i = 0,1, \infty$ corresponds to the coarsened data $X, (X, Z), (X, Z, Y)$ respectively. The MAR assumption implies that,
\[
 \mathbb{P}(C = k | X, Z, Y) =  \omega(k, \mathfrak G_k((X, Z, Y))) \quad \mbox{for all } k \in \{0, 1, \infty\}.
\]
We consider the following definitions, 
\begin{equation*}
    \begin{split}
        m_2^*(\theta, x, z) = \mathbb{P}[Y \leq \theta | X = x, Z= z, C = \infty],& \quad  m_1^*(\theta, x) = \int m_2^*(\theta, x, z) dF(z|C \geq 1, x),\\
         \pi_2^*(x, z) = \mathbb{P}(C = \infty| X = x, Z = z, C \geq 1),& \quad \pi_1^*(x) = \frac{\mathbb{P}(C = 0|X = x)}{\mathbb{P}(C \geq 1|X = x)}.
    \end{split}
\end{equation*}
Note that $m_2^*(\theta, x, z)$ and $\pi_2^*(x, z)$ are respectively the outcome regression and the propensity score functions computed at the second stage (conditional on $\{C \geq 1\}$) of the missing data problem. Note that we have the following alternate representation of $\pi_2^*(x, z)$ in terms of $\omega(\cdot)$ in the MAR assumption,
\begin{equation*}
    \begin{split}
    \pi_2^*(x, z) = &  \mathbb{P}(C = \infty| X = x, Z = z, C \geq 1) =  1 - \frac{ \mathbb{P}(C = 1| X = x, Z = z)}{ \mathbb{P}(C \geq 1| X = x, Z = z)}=  1 - \frac{\omega(1, (x, z))}{1 - \omega(0, x)}. 
    \end{split}
\end{equation*}
The definitions of the first stage outcome regression $m_1^*(\theta, x)$ and propensity score function $\pi_1^*(x)$ are analogous to the corresponding definitions for the covariate shift problem discussed in earlier sub-sections. Note that in terms of $\omega(\cdot)$ in the MAR assumption, propensity score function $\pi_1^*(x)$ has the alternate representation $\pi_1^*(x) = \omega(0, x)/(1 - \omega(0, x))$. 
Using the general semiparametric theory developed by \cite{robins1994estimation}, one can show that $\mathrm{IF}$ (defined below) is the efficient influence function for the $(1 - \alpha)$-th quantile of $Y| C = 0$; see \Cref{app:proof_twostage_missingness}  for a stand-alone derivation of this result,
\begin{equation}
\label{eq:influ_mmdp}
    \begin{split}
       \mathrm{IF}(\theta,c, \mathfrak G_c((x,z, y)), m_2, m_1, \pi_2, \pi_1)  =& \frac{\textbf{1} \{ c = \infty\}\pi_1(x)}{\pi_2(x, z)} \left( \textbf{1}\{y \leq \theta \} - (1 - \alpha)  \right) \\
        & - \textbf{1} \{ c \geq 1\}\pi_1(x)\left( \frac{ \textbf{1} \{c = \infty\}}{\pi_2(x, z)} - 1 \right)  \left[m_2(\theta, x, z) - (1 - \alpha) \right] \\
        &- \left(  \textbf{1} \{ c \geq 1\}\pi_1(x) - \textbf{1} \{ c = 0\} \right)  \left[ m_1(\theta, x) - (1 - \alpha)   \right].
    \end{split}
\end{equation}
As before, we construct prediction sets of the form $\widehat{\mathcal C}_{n,\alpha}(x) = \{y \in \mathbb{R}: y \leq \widehat r_{\alpha}( x)\}$. More generally, the condition $y \leq \widehat r_{\alpha}( x)$ can be replaced by $R(x,y) \leq \widehat r_{\alpha}( x)$ where $R(x,y)$ is any non-conformity score function (see \Cref{rem:nc_score}). The construction of outcome prediction sets by $\rscp$ in the two-stage monotone missing data problem follows \Cref{algg:split_robust_mmdp} in
\Cref{app:algo}.  We summarize the construction here. Define the following for $\theta,\theta'\in\mathbb R$,
\begin{equation}
\label{eq:gprime_2mmdp}
\begin{split}
g'(O,\theta,\pi_1,m_1,m_2)
={}&
\textbf{1}\{C\geq 1\}\pi_1(X)
\left[m_2(\theta,X,Z)-(1-\alpha)\right]
\\
&-
\left[
\textbf{1}\{C\geq 1\}\pi_1(X)
-
\textbf{1}\{C=0\}
\right]
\left[m_1(\theta,X)-(1-\alpha)\right],
\end{split}
\end{equation}
\begin{equation}
\label{eq:g_2mmdp}
\begin{split}
g(O,\theta,\theta',m_2,m_1,\pi_2,\pi_1)
={}&
\frac{\textbf{1}\{C=\infty\}\pi_1(X)}
{\pi_2(X,Z)}
\left[
\textbf{1}\{Y\leq\theta\}
-
m_2(\theta',X,Z)
\right]
\\
&+
g'(O,\theta',\pi_1,m_1,m_2).
\end{split}
\end{equation}
When $\theta'=\theta$, the function in
\eqref{eq:g_2mmdp} coincides with the influence function in \eqref{eq:influ_mmdp}. The $\rscp$ prediction set is constructed as follows.
\begin{enumerate}
    \item Using the first training sample $\mathcal D^{tr}_1$, obtain
    estimates $\widehat\pi_1(x),
        \widehat\pi_2(x,z),
        \widehat m_1(\theta,x),
        \widehat m_2(\theta,x,z)$ of $\pi_1^*,\pi_2^*,m_1^*,m_2^*$, respectively.

    \item Using the independent second training sample
    $\mathcal D^{tr}_2=\{O_i^{t_2}\}_{i=1}^{n_2}$, compute the
    preliminary threshold,
    \begin{equation}
    \label{eq:rscp_init_2mmdp}
        \widehat r_{\alpha,0}
        =
        \inf\left\{
        \theta\in\mathbb R:
        \frac{1}{n_2}
        \sum_{i=1}^{n_2}
        g'\left(
        O_i^{t_2},
        \theta,
        \widehat\pi_1,
        \widehat m_1,
        \widehat m_2
        \right)
        \geq0
        \right\}.
    \end{equation}

    \item Let $\overline\pi\geq1$ be such that, almost surely, $0\leq\widehat\pi_1(X)\leq\overline\pi$ and $ 0<
        {\widehat\pi_1(X)}/
        {\widehat\pi_2(X,Z)}
        \leq\overline\pi$. Using the calibration sample
    $\mathcal D^{cal}=\{O_i\}_{i=1}^n$, compute,
    \begin{equation}
    \label{alg:split_robust_mmdp}
    \begin{split}
        \widehat r_\alpha
        =
        \inf\Bigg\{
        \theta\in\mathbb R:
        \frac{1}{n+1}
        \sum_{i=1}^n
        g\left(
        O_i,
        \theta,
        \widehat r_{\alpha,0},
        \widehat m_2,
        \widehat m_1,
        \widehat\pi_2,
        \widehat\pi_1
        \right)
        -
        \frac{2\overline\pi}{n+1}
        \geq0
        \Bigg\}.
    \end{split}
    \end{equation}

    \item Return the prediction set $\widehat{\mathcal C}_{n,\alpha}
        =
        \{
        y\in\mathbb R:
        y\leq\widehat r_\alpha
        \}.$
\end{enumerate}

We consider the following definitions of nuisance estimation errors before stating the coverage guarantee,
\begin{equation}
\label{eq:two_stage_error_notation}
\begin{split}
A_2
&:=
\|\widehat\pi_1\|_4
\left\|
\frac{\pi_2^*}{\widehat\pi_2}-1
\right\|_4,
\qquad
A_1
:=
\|\widehat\pi_1-\pi_1^*\|_2,
\\
E_2
&:=
\sup_{\theta\in\mathbb R}
\left\|
m_2^*(\theta,X,Z)
-
\widehat m_2(\theta,X,Z)
\right\|_2, \qquad
E_1
:=
\sup_{\theta\in\mathbb R}
\left\|
m_1^*(\theta,X)
-
\widehat m_1(\theta,X)
\right\|_2,
\\
D_2
&:=
\left\|
m_2^*(\widehat r_\alpha,X,Z)
-
m_2^*(\widehat r_{\alpha,0},X,Z)
\right\|_2,\qquad D_1
:=
\left\|
m_1^*(\widehat r_\alpha,X)
-
m_1^*(\widehat r_{\alpha,0},X)
\right\|_2.
\end{split}
\end{equation}



\begin{theorem}[Coverage guarantee for $\rscp$ under two-stage monotone missingness]
\label{thm:2mmdp}
Suppose that the MAR condition \eqref{eq:mar_monotone} holds and
$\mathbb P(C=0)>0$.  Assume that
$\widehat m_1(\theta,x),\widehat m_2(\theta,x,z)\in[0,1]$,
that $\widehat\pi_1(X)\geq0$, $\widehat\pi_2(X,Z)>0$, and that, for
some $\overline\pi\geq1$, we have $\widehat\pi_1(X)\leq\overline\pi$ and ${\widehat\pi_1(X)}/
    {\widehat\pi_2(X,Z)}
    \leq\overline\pi$ almost surely. Suppose further that $Y$ is continuously distributed and that the maps
$\theta\mapsto\widehat m_1(\theta,x)$ and
$\theta\mapsto\widehat m_2(\theta,x,z)$ are right-continuous for every
$x,z$.  Then the threshold $\widehat r_\alpha$ defined in
\eqref{alg:split_robust_mmdp} satisfies the following conditional on both the training splits,
\begin{equation}
\label{eq:two_stage_simple_coverage}
\begin{split}
&
\left|
\mathbb P
\left(
Y_{n+1}\leq\widehat r_\alpha
\,\middle|\,
C_{n+1}=0
\right)
-(1-\alpha)
\right| \leq
\frac{
A_2(E_2+D_2)
+
A_1(E_1+D_1)
}{
\mathbb P(C=0)
}
+
\frac{4\overline\pi}
{\mathbb P(C=0)(n+1)}.
\end{split}
\end{equation}
Suppose, in addition, that there exist $\eta>0$ and constants
$\mathfrak C_0,\mathfrak C_1,\mathfrak C_2,\mathfrak C_3>0$ such that
the following conditions hold on $\mathcal N_\alpha
    =
    [r_\alpha-\eta,r_\alpha+\eta]$: First, for all $\theta,\theta'\in\mathcal N_\alpha$,
\begin{equation}
\label{eq:two_stage_lipschitz_gprime}
\left\|
g'\left(
O,\theta,\widehat\pi_1,\widehat m_1,\widehat m_2
\right)
-
g'\left(
O,\theta',\widehat\pi_1,\widehat m_1,\widehat m_2
\right)
\right\|_2
\leq
\mathfrak C_0|\theta-\theta'|.
\end{equation}
Second, the following lipschitz conditions hold,
\begin{equation}
\label{eq:two_stage_lipschitz_m}
\begin{split}
\left\|
m_1^*(\theta,X)-m_1^*(\theta',X)
\right\|_2
&\leq
\mathfrak C_2|\theta-\theta'|, \qquad 
\left\|
m_2^*(\theta,X,Z)-m_2^*(\theta',X,Z)
\right\|_2
\leq
\mathfrak C_3|\theta-\theta'|.
\end{split}
\end{equation}
Third, the weighted distribution function is locally well separated,
\begin{equation}
\label{eq:two_stage_margin}
\left|
P\left[
\textbf{1}\{C=\infty\}
\frac{\widehat\pi_1(X)}{\widehat\pi_2(X,Z)}
\left\{
\textbf{1}\{Y\leq\theta\}
-
\textbf{1}\{Y\leq\theta'\}
\right\}
\right]
\right|
\geq
\mathfrak C_1|\theta-\theta'|
\end{equation}
for all $\theta,\theta'\in\mathcal N_\alpha$.  We also assume that,
for every fixed $\theta$, $\mathbb P(
J(\widehat m_1(\cdot,X),\theta)>0
)
=
\mathbb P(
J(\widehat m_2(\cdot,X,Z),\theta)>0
)
=0$. Then, for every $\delta\in(0,1)$, there exists a finite constant
$C_\delta$, depending only on $\delta$, $\overline\pi$, $\eta$ and the
constants in \eqref{eq:two_stage_lipschitz_gprime}--\eqref{eq:two_stage_margin},
such that, with probability at least $1-\delta$ over the training split used to construct $\widehat r_{\alpha, 0}$,
\begin{equation}
\label{eq:two_stage_high_prob_coverage}
\begin{split}
&
\left|
\mathbb P
\left(
Y_{n+1}\leq\widehat r_\alpha
\,\middle|\,
C_{n+1}=0
\right)
-(1-\alpha)
\right| \leq
\frac{4\overline\pi}
{\mathbb P(C=0)(n+1)} +
\frac{C_\delta}{\mathbb P(C=0)}
\Bigg[
A_2
\left\{
E_2
+
\frac1{\sqrt n}
+
\frac1{\sqrt{n_2}}
\right\}
\\
&\hspace{35mm}
+
A_1
\left\{
E_1
+
\sup_{\theta\in\mathbb R}
\left\|
\widehat m_1^\dagger(\theta,X)
-
\widehat m_1(\theta,X)
\right\|_2
+
\frac1{\sqrt n}
+
\frac1{\sqrt{n_2}}
\right\}
\Bigg],
\end{split}
\end{equation}
where $\widehat m_1^\dagger(\cdot,\cdot)$ is defined as,
\begin{equation}
\label{eq:m1dagger_2mmdp}
\widehat m_1^\dagger(\theta,X)
:=
\mathbb E
\left[
\widehat m_2(\theta,X,Z)
\,\middle|\,
X,C\geq1
\right].
\end{equation}
\end{theorem}

As in the classical missing data problem and in the general missing data setting, the first bound in \Cref{thm:2mmdp} separates the two sources of error: the terms
$A_2E_2$ and $A_1E_1$ are the usual second-order, multiply robust terms quantifying the bias in the population moment arising from nuisance estimation. The terms $A_2D_2$ and $A_1D_1$ quantify the
additional error introduced because of freezing the threshold-dependent outcome
regressions at the initial estimate $\widehat r_{\alpha,0}$. The second part of the theorem shows that this additional error
term can be removed under mild local regularity conditions. If
$\widehat m_1$ is obtained by a regression method for which $\sup_{\theta}
\|
\widehat m_1^\dagger(\theta,X)
-
\widehat m_1(\theta,X)
\|_2
=
O(n_2^{-1/2})$, as is the case under the conditions for the Forster Warmuth estimator
discussed below, then \eqref{eq:two_stage_high_prob_coverage} simplifies
to the following,
\[
\begin{split}
&
\left|
\mathbb P
\left(
Y_{n+1}\leq\widehat r_\alpha
\,\middle|\,
C_{n+1}=0
\right)
-(1-\alpha)
\right|
\\
&\qquad\leq
\frac{4\overline\pi}
{\mathbb P(C=0)(n+1)}
+
\frac{C_\delta}{\mathbb P(C=0)}
\left[
A_2
\left(
E_2+n^{-1/2}+n_2^{-1/2}
\right)
+
A_1
\left(
E_1+n^{-1/2}+n_2^{-1/2}
\right)
\right].
\end{split}
\]
\begin{remark}[Multiple robustness and finite-sample validity]
\label{rem:two_stage_multiple_robustness}
The proof of \Cref{thm:2mmdp} in fact establishes the following sharper lower bound on coverage,
\[
\mathbb P
\left(
Y_{n+1}\leq\widehat r_\alpha
\,\middle|\,
C_{n+1}=0
\right)
\geq
1-\alpha
-
\frac{
A_2(E_2+D_2)
+
A_1(E_1+D_1)
}{
\mathbb P(C=0)
},
\]
This display makes clear that $\rscp$ enjoys an exact finite-sample
validity whenever $A_1=A_2=0$, irrespective of the quality of the estimated
outcome-regression functions or of the preliminary threshold
$\widehat r_{\alpha,0}$. In other words, when both propensity
score functions are correctly specified, $\rscp$ provides finite-sample validity even if both outcome-regression
functions are misspecified.

More generally, the high-probability bound in
\eqref{eq:two_stage_high_prob_coverage} reveals the $2^2$-robust
structure of the two-stage monotone missingness problem. The leading
nuisance-estimation bias term, $A_1E_1+A_2E_2$, vanishes whenever $(A_1=0\ \mbox{or}\ E_1=0)$ and $(A_2=0\ \mbox{or}\ E_2=0)$. Therefore at each stage, for the leading nuisance-estimation bias term to vanish, it is sufficient to correctly specify either the propensity-score component or the corresponding outcome-regression component. Under the additional regularity conditions, the coverage slack in any of these four cases is, up to a
constant, $(A_1+A_2)
(
n^{-1/2}+n_2^{-1/2}
)
+
O(n^{-1})$.
\end{remark}

To implement \Cref{alg:split_robust_mmdp}, we must estimate the associated nuisance functions. The propensity models $\pi_1^*(x), \pi_2^*(x, z)$ can be learned via logistic regression or through non-parametric approaches such as random forests. The second stage outcome regression function $m_2^*(\theta, x, z)$ can be estimated by running a quantile random forest on the subset of complete $(C = \infty)$ observations. For estimating the first stage outcome regression function $m_1^*(\theta, x)$ we first define the pseudo outcome $H^*(\theta; O)$,
\[
H^*(\theta; O) = \frac{\textbf{1}\{ C = \infty\}}{\pi_2^*(X, Z)}\left[\textbf{1}\{Y \leq \theta\} - m_2^*(\theta, X, Z) \right] + m_2^*(\theta, X, Z).
\]
We observe that the outcome regression function $m_1^*(\theta, x)$ is the conditional mean of the pseudo outcome $H^*(\theta; O)$ given $X, C \geq 1$ i.e.\ $ m_1^*(\theta, X) =  \mathbb{E}[H^*(\theta; O)|X, C \geq 1]$.
Using the estimated nuisance functions $\widehat m_2(\theta, x, z)$ and $\widehat \pi_2(x, z)$ we construct the pseudo outcomes $\widehat H(\theta, O)$ for all the individuals in the $\{C \geq 1\}$ set,
\begin{equation}
\label{eq:estimated_pseudo_outcome}
    \widehat H(\theta; O_i) = \frac{\textbf{1}\{ C_i = \infty\}}{\widehat \pi_2(X_i, Z_i)}\left[\textbf{1}\{Y_i \leq \theta\} - \widehat m_2(\theta, X_i, Z_i) \right] + \widehat m_2(\theta, X_i, Z_i) \quad \mbox{for all} \quad i \in \{C \geq 1\}.
\end{equation}
To estimate the outcome regression function $m_1^*(\theta, x)$, we fit a non-parametric method (for instance random forest) on the $\{C \geq 1\}$ sub-sample using the pseudo outcomes $\widehat H(\theta; O_i)$ as the response variable i.e.\ the estimated outcome regression function satisfies $\widehat m_1(\theta, x) = \widehat{\mathbb{E}}[\widehat H(\theta, O)| X, C \geq 1]$. To derive theoretical results on the estimation error, one can apply results from \cite{yang2023forster}.

For instance, if we use the Forster-Warmuth [\cite{forster2002relative}] series regression estimator based on the first $J$ elements of a fundamental sequence $\Psi \coloneqq \{\phi_1(\cdot), \phi_2(\cdot), \cdots\}$ of functions in $L_2$ , then corollary-$1$ in \cite{yang2023forster} guarantees the following upper bound on the second term, 
\begin{equation*}
\begin{split}
  &\|\widehat{\mathbb{E}}[\widehat H(\theta, O)| X, C \geq 1]  - \mathbb{E}[\widehat H(\theta; O)|X, C \geq 1]\|_2 \\
  \leq & \sqrt{\frac{2 \sigma^2(\theta) J}{|\{C \geq 1\}|}} + \sqrt{2 \kappa} E_J^{\Psi}(m_1^*) + \sqrt{2} \left\lVert \frac{\pi_2^*(X, Z)}{\widehat \pi_2(X, Z)} - 1 
 \right\rVert_4  \| m_2^*(\theta, X, Z) - \widehat m_2(\theta, X, Z)\|_4, 
\end{split}
\end{equation*}
where $\sigma^2(\theta)$ is an upper bound on $\mathbb{E}[\widehat H(\theta, O)^2|X, C \geq 1]$, $\kappa$ is an upper bound on the density of $X$ w.r.t.\ the Lebesgue measure, $|\{C \geq 1\}|$ denotes the number of units in the set $\{C \geq 1 \}$, and $E_J^{\Psi}(m_1^*)$ denotes the $J$-th degree approximation error of the function $m_1^*$ by the first $J$ functions in the fundamental sequence $\Psi$. 
\begin{remark}[Coverage guarantee for $\rsp$ under two-stage monotone missingness]
\label{rem:rsp_two_stage_monotone}
The general $\rsp$ result in \Cref{thm:gen_influ_func} applies directly to the two-stage monotone missingness problem under mild regularity conditions. To see this, let $p_0=\mathbb P(C=0)$ and define,
\[
f_{\mathrm{2m}}
\left(
O,\theta;
m_2,m_1,\pi_2,\pi_1
\right)
:=
(1-\alpha)
+
\frac{1}{p_0}
\mathrm{IF}
\left(
\theta,O,m_2,m_1,\pi_2,\pi_1
\right),
\]
where $\mathrm{IF}$ is defined in \eqref{eq:influ_mmdp}. It is easy to check that $P f_{\mathrm{2m}}
(
O,\theta;
m_2^*,m_1^*,\pi_2^*,\pi_1^*
)
=
\mathbb P(Y\leq\theta\mid C=0)$, so Assumption~\ref{assump:gen_identification} holds. Moreover, the usual stagewise expansion of the influence function gives
\[
\sup_{\theta\in\mathbb R}
\left|
P f_{\mathrm{2m}}
\left(
O,\theta;
\widehat m_2,\widehat m_1,\widehat\pi_2,\widehat\pi_1
\right)
-
P f_{\mathrm{2m}}
\left(
O,\theta;
m_2^*,m_1^*,\pi_2^*,\pi_1^*
\right)
\right|
\leq
d_{\mathrm{2m}}, \quad d_{\mathrm{2m}}
:=
\frac{A_2E_2+A_1E_1}{\mathbb P(C=0)}.
\]
Thus Assumption~\ref{assump:gen_nuisance} follows from the same multiple-robust bias decomposition underlying \Cref{thm:2mmdp}. Under bounded inverse-probability weights and monotone, right-continuous outcome-regression functions, the map
$\theta\mapsto f_{\mathrm{2m}}(O,\theta;\widehat m_2,\widehat m_1,\widehat\pi_2,\widehat\pi_1)$
is a finite linear combination of bounded monotone functions and therefore has uniformly bounded total variation. If, in addition, $Y$ is continuously distributed and the estimated outcome-regression functions have no fixed jumps, Assumption~\ref{assump:gen_bv} follows. Finally, the same computation as in the classical missing data case yields the Bernstein-type second-moment bound in Assumption~\ref{assump:gen_local_variance},
\[
P g_{\mathrm{2m}}^2(O,\theta)
\leq
K_{\mathrm{2m}}
\left|
P g_{\mathrm{2m}}(O,\theta)
\right|
+
d'_{\mathrm{2m}}, \quad d'_{\mathrm{2m}}
=
\frac{\mathfrak C_{\mathrm{2m}}}{\mathbb P(C=0)^2}
\left\{
A_2E_2+A_1E_1+E_2^2+E_1^2
\right\}.
\]
where $g_{2m}(\cdots)$ is the centered version of the function $f_{2m}(\cdots)$ and $\mathfrak C_{\mathrm{2m}}<\infty$ is a constant. Consequently, if $\widehat r_\alpha^{\rsp}$ denotes the threshold obtained by applying \Cref{alg:split_robust} to $f_{\mathrm{2m}}$, then Theorem~\ref{thm:gen_influ_func} gives the following conditional on the training data,
\[
\begin{split}
\mathbb P
\left(
Y_{n+1}\leq\widehat r_\alpha^{\rsp}
\,\middle|\,
C_{n+1}=0
\right)
\geq{}&
1-\alpha-d_{\mathrm{2m}} -
\mathfrak C_n
\left(
\sqrt{d_{\mathrm{2m}}}+n^{-1/4}
\right)
\left(
n^{-1/2}
+
n^{-1/4}\sqrt{d'_{\mathrm{2m}}}
\right),
\end{split}
\]
where $\mathfrak C_n=\mathfrak C\log(en)$. In particular, if all nuisance functions are correctly specified, the remaining coverage slack is the $O(\log n/n^{3/4})$ empirical-process remainder term characterized in \Cref{thm:gen_influ_func}. We omit the detailed proof for brevity.
\end{remark}

To keep the main text focused, we defer the corresponding developments for the general multi-stage monotone missingness problem and the non-monotone missingness problem to \Cref{app:multi_stage_monotone} and \Cref{sec:non_monotone}, respectively. The broad principle remains unchanged: we first use semiparametric theory to construct an observed-data estimating function for the target distribution and then apply $\rsp$ by directly solving its empirical moment equation. The $\rscp$ construction can likewise be applied by freezing all threshold-dependent nuisance functions at an independent preliminary root whenever this leaves a monotone estimating function. This includes both the multi-stage monotone and non-monotone MAR settings considered in this paper. For ease of exposition, we focus on $\rsp$ when discussing the non-monotone missingness problem and when evaluating its empirical performance.

\section{Simulation studies}
\label{sec:simulation}
In this section we demonstrate the performance of $\rsp$ and $\rscp$ on a synthetic and a real data set, under multi-stage monotone and non-monotone missingness structures. Unless the non-conformity score function $R(\cdot, \cdot)$ that we use in our method is independent of the training data $\mathcal{D}^{tr}$, we further split the training data $\mathcal{D}^{tr}$ into two parts: we use the first part for learning the non-conformity score function $R(\cdot, \cdot)$, and we use the second part to estimate the various underlying nuisance functions. In the following analysis, we evaluate the performance of $\rsp$ (under multi-stage monotone and non-monotone missingness) and $\rscp$ (under multi-stage monotone missingness) with respect to both coverage and interval width, considering different nonconformity score functions, including the absolute residual [\cite{lei2018distribution}], conformalized quantile regression (CQR) [\cite{romano2019conformalized}], and the inverse quantile score [\cite{romano2020classification}].
\subsection{Synthetic data}
\label{subsec:simulation}
We use the following simulation setting to compare $\rsp$ and $\rscp$. For each unit $i \in \mathcal{D}$, we generate $(X_i,Z_{1i},Z_{2i},Y_i,T_i,S_i)$ independently of the other units as follows.

\textbf{Setting $1$:}
\begin{equation}
\label{eq:sim_dgp}
    \begin{split}
        X_i,Z_{1i},\epsilon_{1i}',\epsilon_{2i}'
        \stackrel{\mathrm{iid}}{\sim}{}& \mathcal{N}(0,1),\\
        Z_{2i}={}&2X_i+Z_{1i}+0.2\epsilon_{1i}',\\
        Y_i={}&2X_i+2Z_{1i}+0.6Z_{2i}+0.2\epsilon_{2i}',\\
        T_i\mid X_i \sim{}&\operatorname{Bernoulli}
        \!\left\{\operatorname{expit}(0.2X_i)\right\},\\
        S_i\mid X_i,Z_{1i},Z_{2i},T_i \sim{}&
        \begin{cases}
            \operatorname{Bernoulli}\!\left\{
            \operatorname{expit}(0.1X_i+0.1Z_{1i}+0.05Z_{2i})
            \right\}, & T_i=0,\\
            0, & T_i=1,
        \end{cases}
    \end{split}
\end{equation}
where $\operatorname{expit}(u)=1/(1+\exp(-u))$. For units with $T_i=1$, we set $C_i=0$, so only the baseline covariate $X_i$ is observed. For units with $T_i=0$, we observe $(X_i,Z_{1i},Z_{2i})$, so $C_i\geq1$. Finally, the complete data vector $(X_i,Z_{1i},Z_{2i},Y_i)$ is observed only when $T_i=0$ and $S_i=1$, which corresponds to $C_i=\infty$; when $T_i=0$ and $S_i=0$, we set $C_i=1$. This is the two-stage monotone-missingness structure studied in \Cref{sec:monotone}.
For both methods, the absolute-residual score is $R_{\mathrm{AR}}(X,Y)=|Y-\widehat\mu(X)|$ \citep{lei2018distribution}, while the CQR score is $R_{\mathrm{CQR}}(X,Y)
=\max\{\widehat q_{\alpha/2}(X)-Y,
Y-\widehat q_{1-\alpha/2}(X)\}$ \citep{romano2019conformalized}. We use independent samples to learn the score functions and nuisance functions. The baseline mean and quantile functions are fitted using natural-spline regression and quantile regression, respectively. The $\rsp$ algorithm solves the empirical efficient-influence-function equation with the nuisance functions evaluated at each candidate threshold. The $\rscp$ construction first computes $\widehat r_{\alpha,0}$ on an independent initialization sample, freezes the threshold-dependent nuisance functions at $\widehat r_{\alpha,0}$, and then applies the two-stage finite-sample correction in \Cref{alg:split_robust_mmdp} on the calibration sample.

\begin{figure}[ht]
  \centering
  \includegraphics[width=\textwidth]{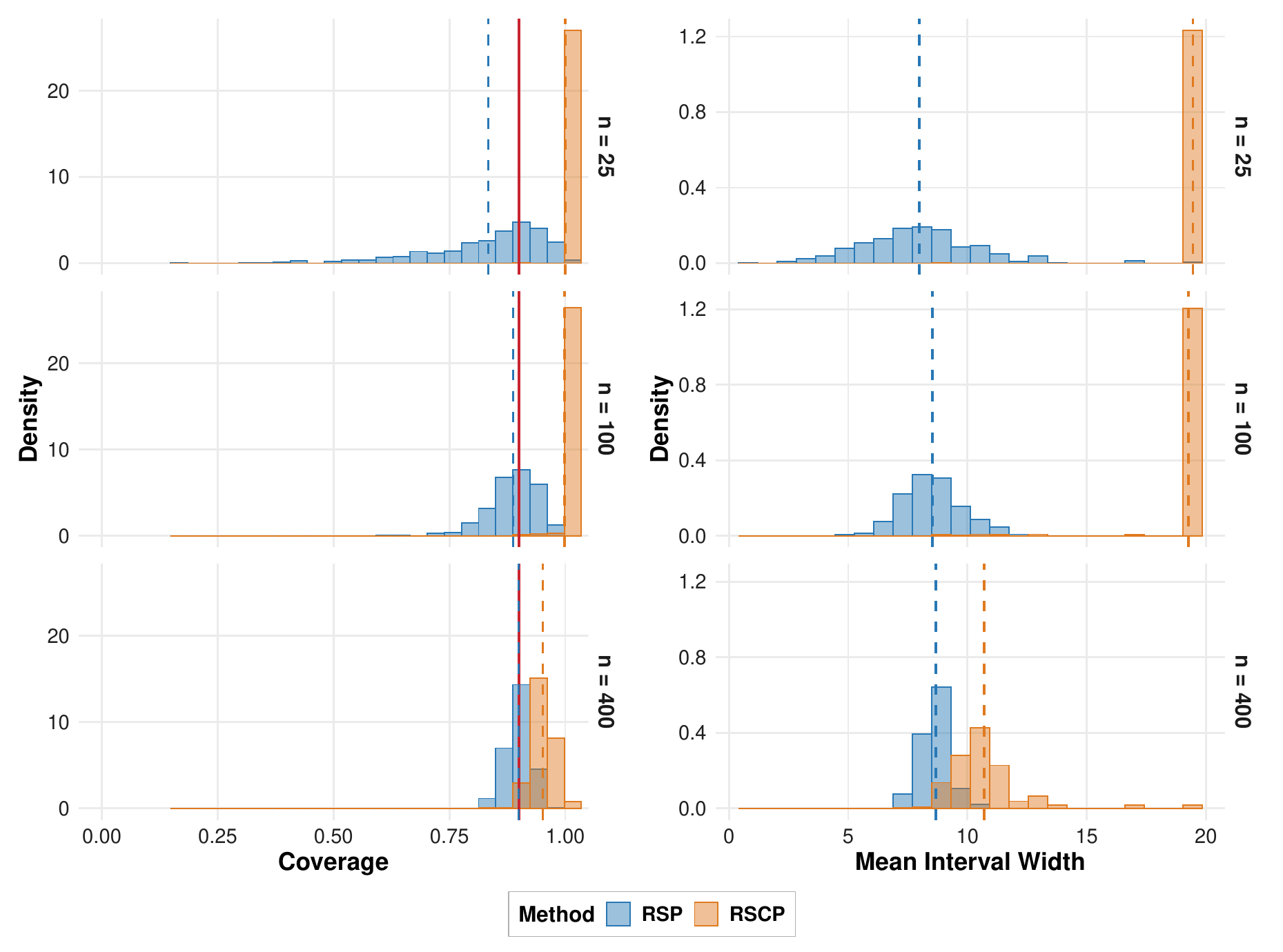}
  \caption{Histograms of coverage (left) and mean interval width (right) for $\rsp$ and $\rscp$ under Setting~$1$, using the absolute-residual non-conformity score. Rows correspond to sizes of the calibration data set $n=25,100,400$.}
  \label{fig:sim_setting1_absolute}
\end{figure}

\begin{figure}[ht]
  \centering
  \includegraphics[width=\textwidth]{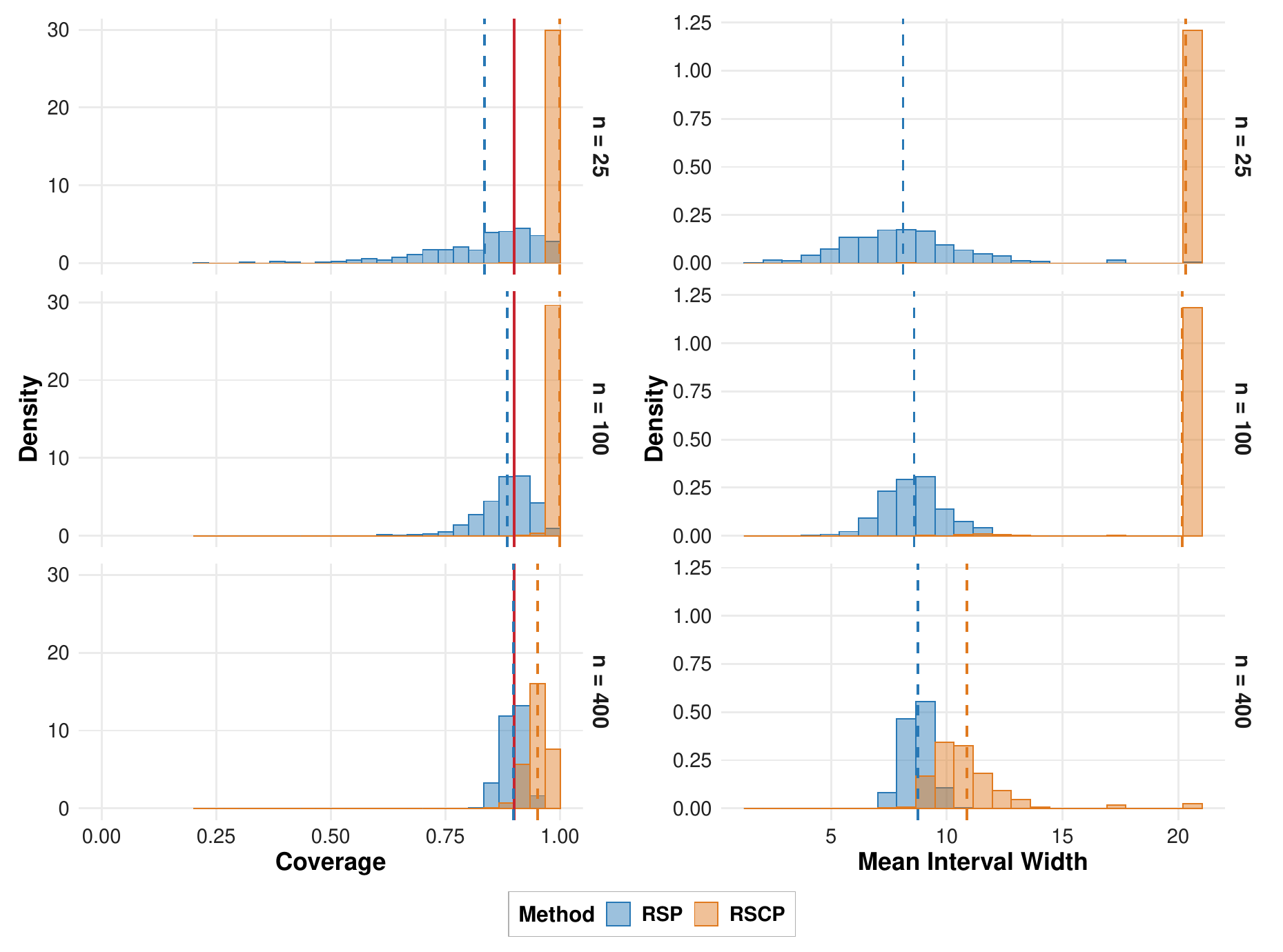}
  \caption{Histograms of coverage (left) and mean interval width (right) for $\rsp$ and $\rscp$ under Setting~$1$, using the CQR non-conformity score. Rows correspond to different sizes of calibration data set $n=25,100,400$.}
  \label{fig:sim_setting1_cqr}
\end{figure}

We set $1-\alpha=0.9$ and consider calibration sizes $n=25,100,400$, using an initialization sample of the same size. For each $n$, coverage conditional on $C=0$ and mean interval width are evaluated on $1500$ independently generated target observations and the experiment is repeated over $500$ Monte Carlo iterations. The colored dashed lines in the figures mark the method-specific Monte Carlo means, and the solid red line in the coverage panels marks the nominal level $0.9$. The results in \Cref{fig:sim_setting1_absolute,fig:sim_setting1_cqr} show the following.

\begin{enumerate}
    \item For $\rsp$, mean coverage approaches the nominal level as the calibration size increases. For the absolute-residual score, the mean coverages at $n=25,100,400$ are $0.834$, $0.887$, and $0.899$, respectively; the corresponding values for CQR are $0.835$, $0.885$, and $0.899$.
    \item The finite-sample correction makes $\rscp$ conservative at the smaller calibration sizes. Its mean coverage decreases toward $0.9$ as $n$ grows: at $n=400$ it is $0.951$ and $0.952$ for the absolute-residual and CQR scores, respectively. This conservativeness is accompanied by wider intervals, especially at $n=25$ and $n=100$. When the corrected equation does not cross zero on the numerical score grid, the upper endpoint of the grid is used; the resulting cases form the mass at the right edge of the width histograms for small sizes of calibration data.
    \item The absolute-residual and CQR scores have similar coverage and efficiency (width of the prediction set) in this setting. At $n=400$, the mean interval widths are $8.668$ and $10.701$ for $\rsp$ and $\rscp$, respectively, with the absolute-residual score, and $8.745$ and $10.880$ with CQR. The variability of both coverage and interval width decreases as $n$ increases.
\end{enumerate}
We provide additional simulations in \Cref{app:add_simul} to illustrate the behavior of both methods when the analyst observes misspecified transformations of the covariates.

\subsection{Real data}
\label{subsec:real_data}
In this subsection, we study the application of $\rsp$ and $\rscp$ to a real data set under both monotone and non-monotone missingness. We use a cohort study comprising $1240$ patients with rheumatoid arthritis from the Wichita Arthritis Center, an outpatient rheumatology facility. We refer readers to \cite{choi2002methotrexate} for a detailed description of the study. Patient characteristics (the baseline covariates) were collected during the initial visit. The intermediate covariates collected at subsequent visits include demographic variables (education level, smoking history, total income, and marital status), clinical variables (tender joint count, grip strength, morning stiffness, health assessment questionnaire disability index, arthritis impact measurement scales, depression scales, visual analogue scales for pain and global disease status, and body mass index), laboratory variables (erythrocyte sedimentation rate, white blood cell count, haemoglobin concentration, and rheumatoid factor), medication use, and self-reported data. The focus of the original study was to estimate the effect of methotrexate, a commonly used antirheumatic therapy, on mortality among the participants.

For our analysis, we use the data collected during the first $9$ months of follow-up. The response of interest, $Y$, is the total tender joint count measured at month $9$. We group the follow-up data into four-month intervals. The data collected at the initial visit form the baseline covariates $X$. In both analyses below, the indicator $\mathbf 1\{C=0\}$ is artificially generated from the baseline covariates according to
\[
\mathbb P(C=0\mid X=x)
=
\frac{\exp(0.01x^\top\mathbf 1)}{1+\exp(0.01x^\top\mathbf 1)}.
\]
Coverage is evaluated using the month-$9$ outcomes that are masked by this artificial coarsening and is therefore computed conditional on $C=0$.

We consider two non-conformity scores. The first is the generalized inverse-quantile score of \cite{romano2020classification}, a score which is well suited to the discrete outcome support $\{3,\ldots,9\}$. The second is the standard absolute-residual score of \cite{lei2018distribution}. The predictive model defining each score is trained on an independent split of the data. We set $1-\alpha=0.9$ and repeat each experiment $100$ times.

\paragraph{Three-stage monotone missingness.}
Let $Z^1$ denote the covariates recorded at month $1$ and let $Z^2$ denote those recorded at month $5$. Among individuals not assigned to $C=0$, we set $C=1$ for patients who drop out during months $1$--$4$, $C=2$ for patients who remain under follow-up beyond month $4$ but drop out during months $5$--$8$, and $C=\infty$ for patients observed through month $9$. Thus, $C=1$, $C=2$, and $C=\infty$ correspond respectively to observing $(X,Z^1)$, $(X,Z^1,Z^2)$, and $(X,Z^1,Z^2,Y)$, giving the three-stage monotone missingness structure described in \Cref{sec:monotone}.

For each repetition, the observed sample is split into a training sample for estimating nuisance functions, an independent sample used to obtain the preliminary $\rscp$ threshold, and a calibration sample. The stagewise propensity scores are estimated with probability forests. The conditional distribution of the non-conformity score at the final stage is estimated with a quantile regression forest, and the earlier conditional distribution functions are obtained through backward iterated regression. We then apply both the $\rsp$ estimating equation and the three-split $\rscp$ construction in \Cref{algg:split_robust_Mmdp}. In contrast to $\rsp$, $\rscp$ freezes the threshold-dependent nuisance regressions at the preliminary root and includes the finite-sample correction in the final calibration equation.

\begin{figure}[htbp]
    \centering
    \begin{subfigure}[t]{0.49\textwidth}
        \centering
        \includegraphics[width=\textwidth]{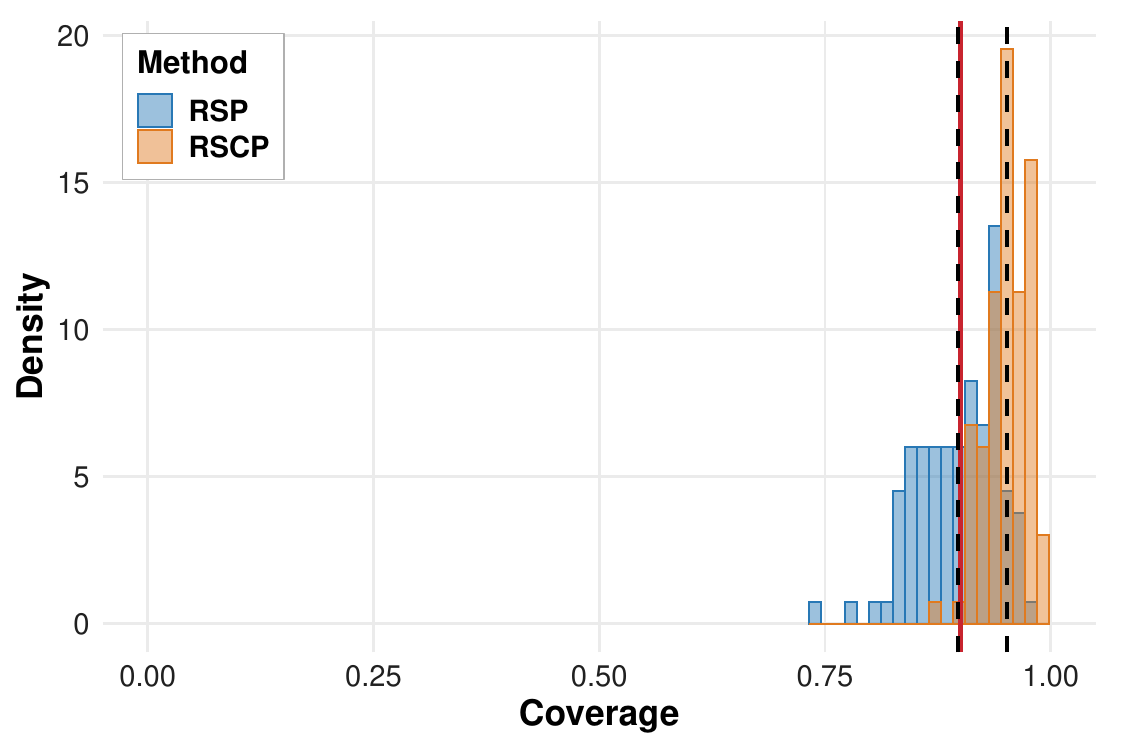}
        \caption{Absolute-residual score.}
        \label{fig:real_monotone_abs}
    \end{subfigure}
    \hfill
    \begin{subfigure}[t]{0.49\textwidth}
        \centering
        \includegraphics[width=\textwidth]{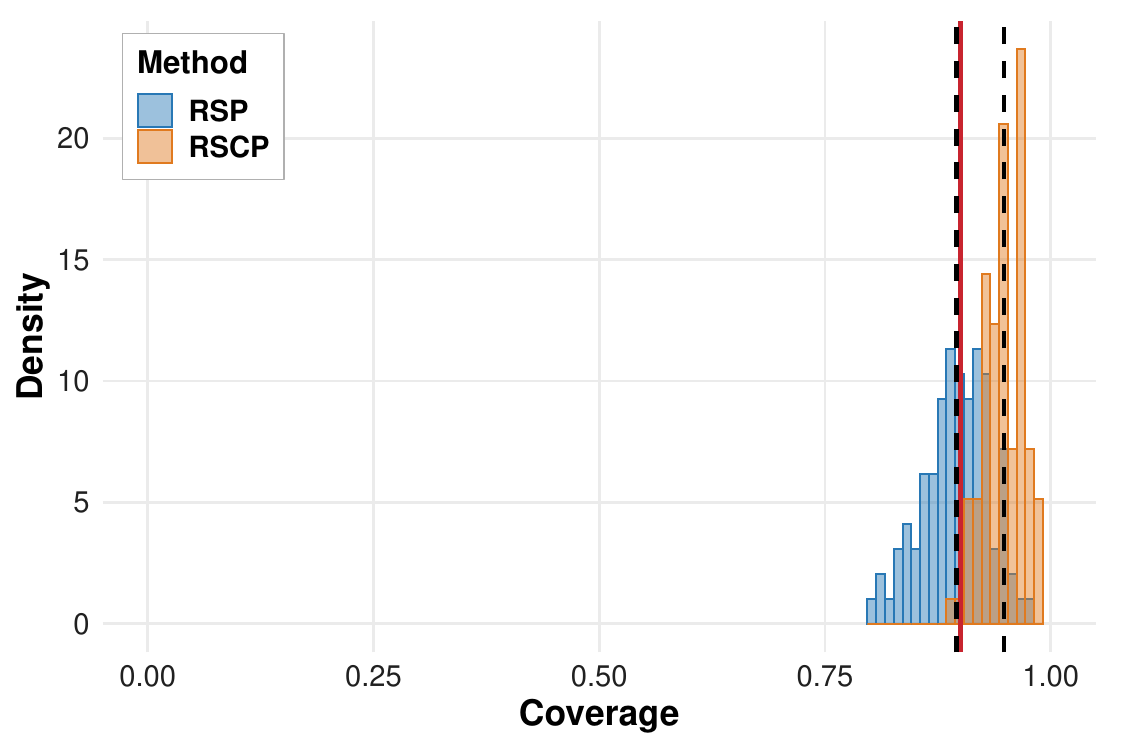}
        \caption{Generalized inverse-quantile score.}
        \label{fig:real_monotone_inv}
    \end{subfigure}
    \caption{Coverage under three-stage monotone missingness. In each panel, the overlapping histograms compare $\rsp$ and $\rscp$ over $100$ repetitions. The solid red line marks the nominal coverage $0.9$, and the dashed lines mark the method-specific mean coverage.}
    \label{fig:real_monotone_coverage}
\end{figure}

\paragraph{Three-stage non-monotone missingness.}
To obtain a non-monotone MAR experiment with known underlying outcomes, we begin with patients having complete measurements at months $0$, $1$, $5$, and $9$, and then mask variables according to an artificial model. Conditional on $C\geq1$, let $U_X$, $U_1$, and $U_2$ be the hyperbolic tangent of the average standardized components of $X$, $Z^1$, and $Z^2$, respectively, and let $U_Y=\tanh\{(Y-\bar Y)/s_Y\}$. We generate three incomplete-data pattern probabilities using the non-monotone MAR construction of \cite{sun2018inverse2}:
\begin{align*}
\rho_1(X,Z^1)
&=\operatorname{expit}(-2+0.25U_X+0.35U_1),\\
\rho_2(X,Z^1,Z^2)
&=\operatorname{expit}(-2-0.20U_X+0.20U_1+0.30U_2),\\
\rho_3(X,Z^1,Y)
&=\operatorname{expit}(-2+0.20U_X-0.25U_1+0.30U_Y),
\end{align*}
and set $\rho_\infty=1-\rho_1-\rho_2-\rho_3$. The bounded summaries and the specified coefficients ensure that $\rho_\infty$ is positive. Pattern $C=1$ reveals $(X,Z^1)$, pattern $C=2$ reveals $(X,Z^1,Z^2)$, pattern $C=3$ reveals $(X,Z^1,Y)$, and $C=\infty$ reveals the full data. In particular, pattern $C=3$ records the month-$9$ outcome despite the month-$5$ covariates being missing, so the resulting observation patterns are not nested. Each pattern probability depends only on variables observed under that pattern, and hence the mechanism is MAR.

The pattern-specific probabilities are re-estimated in the training sample by maximizing the observed-data likelihood proposed by \cite{sun2018inverse2}, with positivity truncation used to stabilize inverse-probability weights. To implement \Cref{algg:split_robust_nmdp}, we approximate the operator $\mathcal M$ in \Cref{sec:non_monotone} on the training data using nearest-neighbour conditional-expectation operators. For every candidate score threshold $\theta$, we solve the resulting finite-dimensional system for $\widehat d(\theta)=\widehat{\mathcal M}^{-1}\{\mathbf 1(R\leq\theta)-(1-\alpha)\}$. These estimates, the fitted pattern probabilities, and the baseline outcome regression are substituted into the efficient influence function in \eqref{eq:influence_nonmonotone_estimated}. The $\rsp$ threshold is the smallest $\theta$ for which the empirical mean of this estimated influence function is nonnegative, and the prediction set is $\{y:R(x,y)\leq\widehat r_\alpha\}$.

\begin{figure}[htbp]
    \centering
    \begin{subfigure}[t]{0.49\textwidth}
        \centering
        \includegraphics[width=\textwidth]{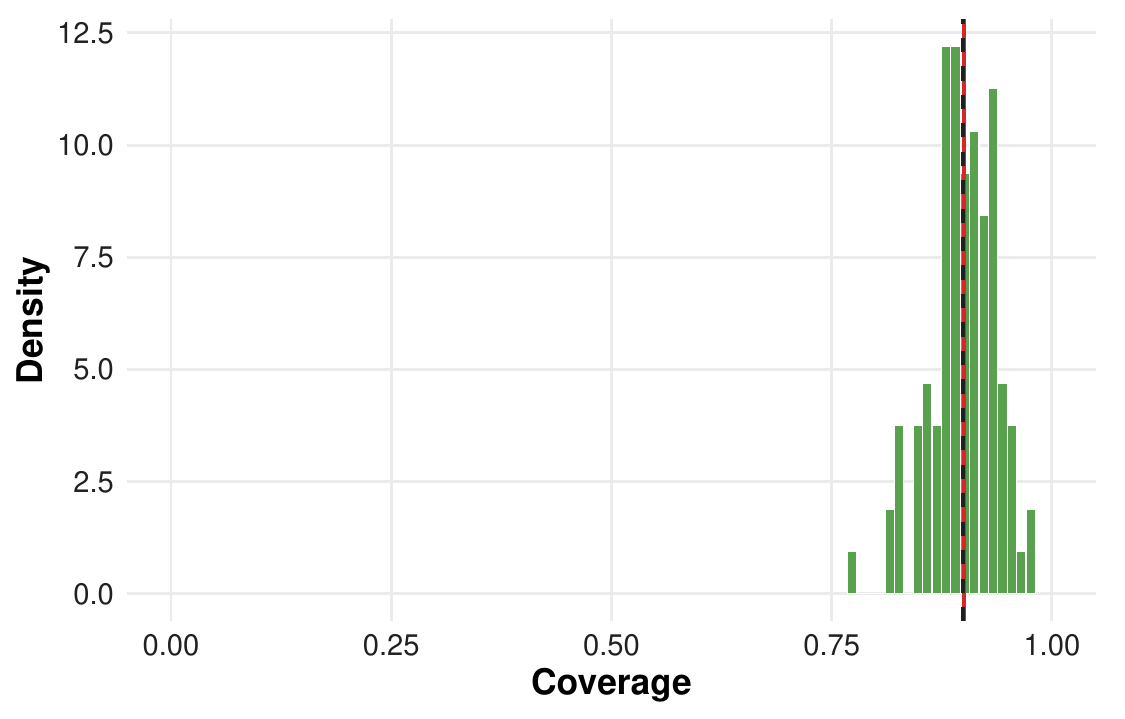}
        \caption{Absolute-residual score.}
        \label{fig:real_nonmonotone_abs}
    \end{subfigure}
    \hfill
    \begin{subfigure}[t]{0.49\textwidth}
        \centering
        \includegraphics[width=\textwidth]{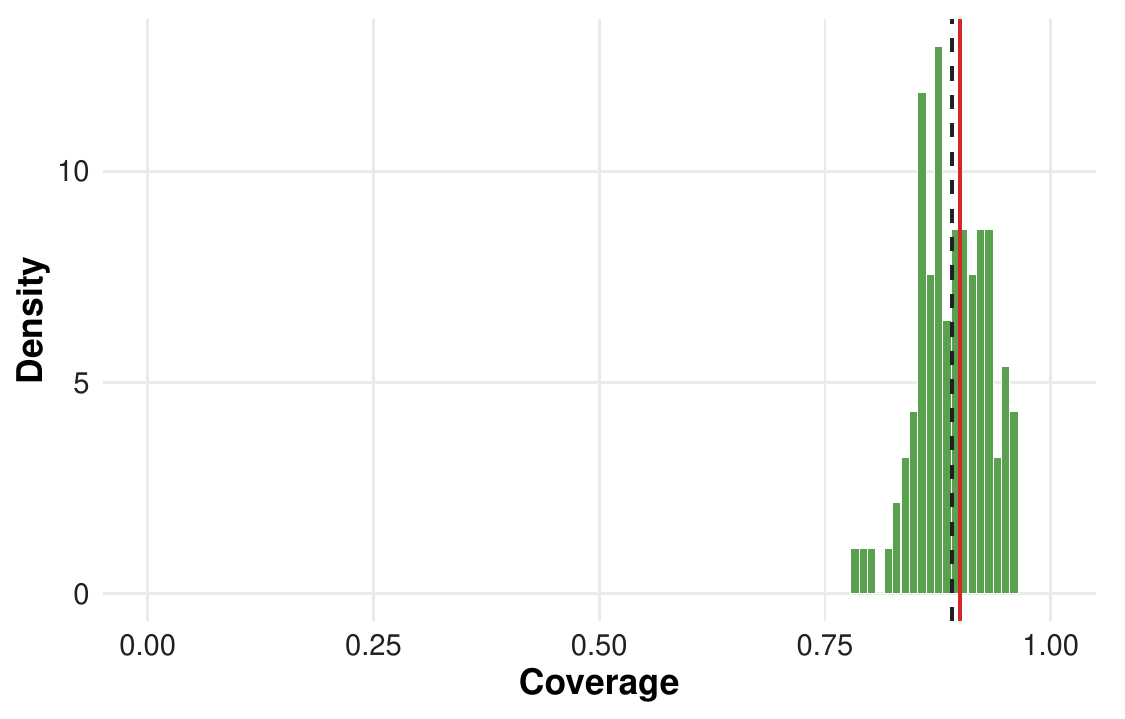}
        \caption{Generalized inverse-quantile score.}
        \label{fig:real_nonmonotone_inv}
    \end{subfigure}
    \caption{Coverage of $\rsp$ under three-stage non-monotone MAR over $100$ repetitions. The solid red line marks nominal coverage $0.9$, and the dashed black line marks mean coverage.}
    \label{fig:real_nonmonotone_coverage}
\end{figure}

\begin{table}[htbp]
\centering
\resizebox{\textwidth}{!}{
\begin{tabular}{lllccc}
\toprule
\textbf{Missingness} & \textbf{Method} & \textbf{Non-conformity score}
& \textbf{First quartile} & \textbf{Median} & \textbf{Third quartile}\\
\midrule
Monotone & $\rsp$  & Absolute residual & 0.864 & 0.904 & 0.933\\
Monotone & $\rsp$  & Inverse quantile  & 0.871 & 0.898 & 0.922\\
Monotone & $\rscp$ & Absolute residual & 0.939 & 0.953 & 0.972\\
Monotone & $\rscp$ & Inverse quantile  & 0.932 & 0.949 & 0.967\\
Non-monotone MAR & $\rsp$ & Absolute residual & 0.880 & 0.899 & 0.927\\
Non-monotone MAR & $\rsp$ & Inverse quantile  & 0.863 & 0.891 & 0.921\\
\bottomrule
\end{tabular}
}
\caption{Empirical coverage conditional on $C=0$ for the rheumatoid arthritis data. Entries are quartiles across $100$ repetitions at nominal coverage $0.9$.}
\label{tab:coverage_summary}
\end{table}

As shown in \Cref{fig:real_monotone_coverage,tab:coverage_summary}, under the multi-stage monotone missingness setting, the coverage histogram of $\rsp$ is centered close to the nominal level for both scores, whereas $\rscp$ is moderately conservative because of its finite-sample correction. Under non-monotone MAR setting, \Cref{fig:real_nonmonotone_coverage,tab:coverage_summary} show that the coverage of $\rsp$ remains centered near $0.9$. These findings support the theoretical robustness results for the monotone and non-monotone settings. The inverse-quantile score remains particularly natural for the discrete tender-joint-count outcome, while the absolute-residual score provides a useful continuous-score benchmark.

\paragraph*{Code availability.}
All code used to implement the proposed methods and reproduce the numerical experiments in this paper is publicly available at
\url{https://github.com/manitpaul/RSCP}.
The repository contains implementations of the $\rsp$ and $\rscp$ procedures together with the code used for the comparisons in the classical, monotone and non-monotone missing-data simulations, and the real-data analyses reported in this paper.
\section{Discussion}
\label{sec:extension}

In this work, we address the challenge of constructing valid prediction sets in the presence of general missing or coarsened data. Focusing on settings where a designated subset of covariates (referred to as baseline covariates) used to build a predictive algorithm is always observed, we introduce the Robust Split Prediction ($\rsp$) and the Robust Split Conformal Prediction ($\rscp$) algorithms.

The first method, Robust Split Prediction ($\rsp$), directly solves the
empirical estimating equation obtained from an efficient influence
function. In the classical missing-data problem this reduces to the
doubly robust construction of \cite{yang2024doubly}; our general
formulation extends the same principle to substantially more general coarsening patterns. The corresponding coverage error consists of
a population-bias term, quantifying the error in nuisance estimation, and an empirical-process remainder term. When the efficient influence function is multiply robust, the leading population-bias term inherits
the same higher-order mixed-bias structure. Moreover, we sharpen the coverage guarantees for the doubly robust conformal method of \cite{yang2024doubly} through a more refined empirical-process analysis, using peeling arguments to control the stochastic fluctuations of the estimating function.

The second method, Robust Split Conformal Prediction ($\rscp$), can be applied whenever the efficient influence function becomes monotone after freezing the threshold-dependent
nuisance functions at an initial root. An independent sample is used to obtain
a preliminary threshold (root), the threshold-dependent nuisance functions are evaluated at this preliminary value, and conformal risk control is
then applied to the resulting monotone estimating function. In the classical and monotone missing-data problems, $\rscp$ yields finite-sample coverage validity when the relevant nuisance
components are correctly specified, whereas the direct $\rsp$
construction retains a stochastic empirical-process remainder even when all the nuisance functions are correctly estimated.

We demonstrate that $\rsp$ or $\rscp$ can be applied to a variety of practically relevant settings, including classical (covariate-shift), monotone or non-monotone missing data problems under coarsening-at-random (CAR) assumption. We validate the empirical utility of $\rsp,\rscp$ through the application to synthetic data and a real-world longitudinal cohort study of rheumatoid arthritis patients, where the goal is to predict total tender joint counts at a future visit based on partially observed covariates. The positive results obtained from these studies highlight the flexibility and robustness of the proposed methods in addressing complex missing data scenarios in predictive inference. A key direction for future research is to formally define the notion of optimality for prediction sets in the presence of missing data, and to explore whether conformal prediction methods can be used to attain this target in the asymptotic regime.

\section{Acknowledgments}
Arun Kumar Kuchibhotla and Eric J. Tchetgen Tchetgen acknowledge support by National Science Foundation under Grant NSF DMS-2210662. 
\bibliography{references}
\bibliographystyle{plainnat}
\newpage
\appendix
\setcounter{section}{0}
\setcounter{equation}{0}
\setcounter{figure}{0}
\setcounter{remark}{0}
\renewcommand{\thesection}{S.\arabic{section}}
\renewcommand{\theequation}{E.\arabic{equation}}
\renewcommand{\thefigure}{A.\arabic{figure}}
\renewcommand{\theremark}{R.\arabic{remark}}
  \begin{center}
  \Large {\bf Supplement to ``Multiply Robust Conformal Risk Control with Coarsened Data''}
  \end{center}
\section{Background on conformal risk control}
\label{app:crc_background}

We briefly review the conformal risk control framework of
\cite{angelopoulos2022conformal} and explain the connection with the
methods developed in this paper. Let $Z_1,\ldots,Z_n,Z_{n+1}$ be exchangeable observations and let
$\{\mathcal C_\lambda:\lambda\in\Lambda\}$ be a nested family of
prediction rules indexed so that increasing $\lambda$ produces a more
conservative prediction set. Let $L_i(\lambda)
=
L(\mathcal C_\lambda,Z_i)$ denote the loss incurred on observation $i$. We assume that, for every $i$, $0\leq L_i(\lambda)\leq B$ and that the map $\lambda\longmapsto L_i(\lambda)$ is non-increasing and right-continuous. The population risk is $\mathcal R(\lambda)
=
\mathbb E[L_{n+1}(\lambda)]$. For a target risk level $\alpha$, we define the estimate,
\begin{equation}
\label{eq:crc_general_threshold}
\widehat\lambda
=
\inf
\left\{
\lambda\in\Lambda:
\frac{
\sum_{i=1}^n L_i(\lambda)+B
}{
n+1
}
\leq
\alpha
\right\}.
\end{equation}
Equivalently, we can define $\widehat \lambda$ as,
\[
\widehat\lambda
=
\inf
\left\{
\lambda:
\frac{n}{n+1}
\widehat{\mathcal R}_n(\lambda)
+
\frac{B}{n+1}
\leq
\alpha
\right\}, \qquad \widehat{\mathcal R}_n(\lambda)
=
\frac1n
\sum_{i=1}^n
L_i(\lambda).
\]
The conformal risk control argument gives,
\begin{equation}
\label{eq:crc_expected_risk}
\mathbb E
\left[
L_{n+1}(\widehat\lambda)
\right]
\leq
\alpha.
\end{equation}
The additional term $B/(n+1)$ is the finite-sample conformal
correction accounting for the unobserved loss of the $(n+1)$-st
exchangeable observation. For ordinary marginal coverage, we take the loss as,
\[
L_i(\theta)
=
\textbf{1}
\left\{
R(X_i,Y_i)>\theta
\right\},
\qquad
B=1.
\]
The corresponding risk is the usual miscoverage probability $\mathcal R(\theta)
=
\mathbb P
\left(
R(X,Y)>\theta
\right)$ and \eqref{eq:crc_general_threshold} becomes,
\[
\widehat\theta
=
\inf
\left\{
\theta:
\frac{
\sum_{i=1}^n
\textbf{1}\{R(X_i,Y_i)\leq\theta\}
}{
n+1
}
\geq
1-\alpha
\right\}.
\]
Thus split conformal prediction is a special case of conformal risk
control with the zero--one miscoverage loss.

\paragraph{Why CRC cannot be applied directly under coarsening.}
In the coarsened-data settings studied in this paper,
$R(X_i,Y_i)$ is not observed for every calibration unit, so the loss
$L_i(\theta)$ cannot be evaluated directly. Semiparametric theory
instead provides an observed-data function $f(O,\theta,\eta^*)$ satisfying,
\[
P f(O,\theta,\eta^*)
=
\mathbb P
\left(
R(X,Y)\leq\theta
\right).
\]
Consequently, $1-f(O,\theta,\eta^*)$ is an observed-data surrogate whose expectation equals the target
miscoverage risk. Importantly, the target remains the
actual miscoverage probability; the influence-function-based quantity
is only a device for identifying it from the coarsened observations.

There are two obstacles to applying
\eqref{eq:crc_general_threshold} directly to this surrogate. First,
$\eta^*$ is unknown and must be estimated. Second, the sample path $\theta\longmapsto f(O,\theta,\eta^*)$
need not be monotone. The standard CRC argument therefore does not apply
directly.

\paragraph{The $\rsp$ route.}
The first approach in this paper is to dispense with pointwise
monotonicity and solve the empirical moment equation directly:
\[
\widehat r_\alpha
=
\inf
\left\{
\theta:
P_n
f(O,\theta,\widehat\eta)
\geq
1-\alpha
\right\}.
\]
This is the Robust Split Prediction ($\rsp$) procedure. Its validity
uses tools from empirical process theory. This results in an additional stochastic
slack in the lower bound on the coverage probability.

\paragraph{The $\rscp$ route.}
Suppose instead that the nuisance functions can be partitioned as $\eta^*(\theta)
=
\left(
\eta_1^*(\theta),
\eta_2^*
\right)$, where $\eta_1^*(\theta)$ contains the threshold-dependent nuisance components. The Robust Split Conformal Prediction ($\rscp$)
procedure first uses an independent sample to compute a preliminary
threshold $\widehat r_{\alpha,0}$ and then replaces
$\widehat\eta_1(\theta)$ by
$\widehat\eta_1(\widehat r_{\alpha,0})$ in the final estimating
equation. The resulting map $\theta
\longmapsto
f
(
O,\theta,
\widehat\eta_1(\widehat r_{\alpha,0}),
\widehat\eta_2
)$ is monotone in the classical and monotone missing-data examples in
this paper. The usual conformal correction term in standard CRC literature can therefore be
applied to the updated estimating function.

This explains both the advantage and the limitation of $\rscp$. When
the required decomposition exists, it converts an
influence-function based estimating equation into a conformal
risk control problem with monotone loss function and can recover finite-sample validity when
the corresponding nuisance functions are estimated correctly. The
$O((n+1)^{-1})$ correction term is the main ingredient of standard CRC that guarantees finite-sample
validity; the substantive contribution here is the construction of the new estimating function (by freezing the threshold-dependent nuisance funtions at an initial root) that makes the multiply robust efficient influence function compatible
with the CRC argument. The same freezing principle is not restricted to monotone coarsening patterns. In the non-monotone MAR problem considered in this paper, freezing both the ordinary threshold-dependent outcome regression and the threshold-dependent operator-derived regression functions leaves only a positively weighted indicator term depending on the final threshold. The resulting estimating function is therefore monotone and is also compatible with the $\rscp$ construction. We nevertheless use $\rsp$ for the detailed non-monotone development and numerical experiments in order to keep the exposition and implementation focused.
\section{Classical missing data problem and its connection to covariate shift}
\label{app:cov_shift}
Covariate shift refers to the situation when the distribution of the test data differs from that of the training data, and this distributional shift may or may not be known a priori. We suppose that the calibration set $\mathcal{D}^{cal} = \{O_i\}_{i = 1}^n$ is composed of two parts $\mathcal{D}^{cal}_P$ and $\mathcal{D}^{cal}_Q$ where,
\[
\mathcal{D}^{cal}_P \coloneqq \{O_i = (X_i, Y_i): 1 \leq i \leq n_1 \}, \quad \mathcal{D}^{cal}_Q \coloneqq \{O_i = X_i: n_1 + 1 \leq i \leq n   \}.
\]
The random variables $(X_i,Y_i) \in \mathcal{D}^{cal}_P$ are i.i.d.\ from $P_X \otimes P_{Y|X}$ and the random variables $X_i \in \mathcal{D}^{cal}_Q $ are i.i.d.\ from $Q_X$. The outcome is missing in all the data points in $D_Q^{cal}$. The different marginal distributions of $X$ in $\mathcal{D}^{cal}_P$ and $\mathcal{D}^{cal}_Q$ imply that there is a covariate-shift between the distribution of the points in these two data sets. \cite{tibshirani2019conformal} studies the covariate shift problem with the objective of constructing a prediction set $\widehat{\mathcal C}_{n,\alpha}(x)$ using the calibration set $D^{cal}$ such that, 
\[
\mathbb{P}\left( Y_{n+1} \in \widehat{\mathcal C}_{n,\alpha}(X_{n+1}) \right) \geq 1 - \alpha, \quad \mbox{where} \quad (X_{n+1},Y_{n+1}) \sim Q_X \otimes P_{Y|X}. 
\]
Let $R(x, y)$ be a non-conformity score function (see \Cref{rem:nc_score}). The prediction set $\widehat{\mathcal C}_{n,\alpha}(x)$ is generally of the form $\widehat{\mathcal C}_{n,\alpha}(x) = \{y \in \mathbb{R}: R(x,y) \leq \widehat r_{\alpha}\}$ where $\widehat r_{\alpha}$ is estimated from the calibration set $\mathcal{D}^{cal}$. Therefore, the validity of the prediction set $\widehat{\mathcal C}_{n,\alpha}(x)$ implies,
\[
\mathbb{P}\left( R(X_{n+1}, Y_{n+1}) \leq \widehat r_{\alpha} \right) \geq 1 - \alpha, \quad \mbox{where} \quad (X_{n+1},Y_{n+1}) \sim Q_X \otimes P_{Y|X}. 
\]
In other words, our goal is to estimate the $(1-\alpha)$-th quantile of $R(X_{n+1}, Y_{n+1})$ from the calibration set where $ (X_{n+1},Y_{n+1}) \sim Q_X \otimes P_{Y|X}$. There are broadly two ways of approaching this problem in the literature. The first one is to model the problem in the following way, 
\[
\mathbb{P}_{(X_{n+1}, Y_{n+1}) \sim Q_X \otimes P_{Y|X}}(R(X_{n+1},Y_{n+1}) \leq \widehat r_{\alpha}) = \mathbb{E}_{(X_{n+1},Y_{n+1}) \sim P_X \otimes P_{Y|X}} \left[\textbf{1}\{R(X_{n+1},Y_{n+1}) \leq \widehat r_{\alpha} \} \frac{dQ_X}{dP_X}(X_{n+1}) \right].
\] 
This requires the assumption that $Q_X \ll P_X$ for the radon-nikodym derivative to exist. This approach has been studied extensively in \cite{tibshirani2019conformal} where the authors propose the weighted conformal prediction (WCP) algorithm to construct valid outcome prediction set when the likelihood ratio $w^*(x) \coloneqq (dQ_X/dP_X)(x)$ is known. The quantile of $R(X_{n+1}, Y_{n+1})$ ($ (X_{n+1},Y_{n+1}) \sim Q_X \otimes P_{Y|X}$) is estimated by the function $\widehat r_{\alpha}(x)$,
\begin{equation}
\label{eq:tibshirani_ralpha}
  \widehat r_{\alpha}(x, w^*) = \mbox{Quantile}(1 - \alpha; \sum_{i = 1}^{n_1} p_i(x, w^*) \delta_{R_i} + p_{n+1}(x, w^*)\delta_{\infty}),   
\end{equation}
where $\{p_i(x, w)\}_{i = 1}^{n_1}$ and $p_{n+1}(x, w)$ are defined as, 
\[
p_i(x, w) = \frac{ w(X_i)}{\sum_{i = 1}^{n_1}  w(X_i) +  w(x)} \mbox{  for  } i \in [n_1] \quad \mbox{and} \quad  p_{n+1}(x, w) = \frac{ w(x)}{\sum_{i = 1}^{n_1}  w(X_i) +  w(x)}. 
\]
Corollary-$1$ in \cite{tibshirani2019conformal} shows that the $\widehat r_{\alpha}(x, w^*)$ (defined in \eqref{eq:tibshirani_ralpha}) gives finite sample valid coverage,
\[
\mathbb{P}_{(X_{n+1}, Y_{n+1}) \sim Q_X \otimes P_{Y|X}}(R(X_{n+1},Y_{n+1}) \leq \widehat r_{\alpha}(X_{n+1}, w^*)) \geq 1 - \alpha.  
\]
However, it is often the case that the true likelihood ratio $w^*(x)$ is not known to us. \cite{lei2021conformal} shows that it is sufficient to have a consistent estimate $\widehat w(x)$ of $w^*(x)$. Specifically, the coverage slack resulting from the estimation of $w^*$ can be characterized explicitly,
\[
\mathbb{P}_{(X_{n+1}, Y_{n+1}) \sim Q_X \otimes P_{Y|X}}( R(X_{n+1},Y_{n+1}) \leq \widehat r_{\alpha}(X_{n+1}, \widehat w) ) \geq  1- \alpha - \frac{1}{2} \mathbb{E}_{X \sim P_X} |\widehat w (X) - w^*(X)|. 
\]
The second way to go about this problem is to write the coverage probability as the expectation of the conditional distribution of the outcome,
\[
\mathbb{P}_{(X_{n+1}, Y_{n+1}) \sim Q_X \otimes P_{Y|X}}(R(X_{n+1},Y_{n+1}) \leq \widehat r_{\alpha}) = \mathbb{E}_{Q_X}[m_R^*(\widehat r_{\alpha}, X_{n+1})], 
\]
where $m_R^*(\theta, x) = \mathbb{P}(R(X,Y) \leq \theta | X = x)$. Being a conditional distribution function, $m_R^*(\theta, x)$ is a non-decreasing function of $\theta \in [0,1]$ or in other words, $1 - m_R^*(\theta, x) = \mathbb{P}(R(X,Y) > \theta | X = x)$ is a non-increasing function of $\theta \in [0,1]$. The process of obtaining suitable $\widehat r_{\alpha}$ in this scenario has been studied in \cite{angelopoulos2022conformal} under the assumption that $m_R^*(\theta, x)$ is known. The quantile of $R(X_{n+1}, Y_{n+1})$ ($ (X_{n+1},Y_{n+1}) \sim Q_X \otimes P_{Y|X}$) is estimated by $\widehat r_{\alpha}(m_R^*)$,
\begin{equation}
\label{eq:crc_eq}
\widehat r_{\alpha}(m_R^*) = \inf\left\{\theta: \frac{n-m}{n-m+1}\widehat S_n(\theta, m_R^*) + \frac{1}{n - m + 1} \leq \alpha \right\},
\end{equation}
where $\widehat S_n(\theta, m_R)$ is defined as,
\[
\widehat S_n(\theta) = 1 - \frac{\sum_{i= m+1}^n m_R(\theta, X_i)}{n - m}. 
\]
Theorem $1$ of \cite{angelopoulos2022conformal} states that $\widehat r_{\alpha}(m_R^*)$ (defined in \eqref{eq:crc_eq}) satisfies,
\[
\mathbb{P}_{(X_{n+1}, Y_{n+1}) \sim Q_X \otimes P_{Y|X}}(R(X_{n+1},Y_{n+1}) \leq \widehat r_{\alpha}(m_R^*)) = \mathbb{E}_{Q_X}[m_R^*(\widehat r_{\alpha}(m_R^*), X_{n+1})] \geq 1- \alpha. 
\]
It can be easily shown that knowing the exact conditional distribution $m_R^*(\theta, x)$ is not required for asymptotically valid coverage. In fact, having a consistent estimate $\widehat m_R(\theta, x)$ of $m_R^*(\theta, x)$ suffices,
\begin{equation*}
    \begin{split}
        &\mathbb{P}_{(X_{n+1}, Y_{n+1}) \sim Q_X \otimes P_{Y|X}}(R(X_{n+1},Y_{n+1}) \leq \widehat r_{\alpha}(\widehat m_R)) \\
        = &\mathbb{E}_{Q_X}[m_R^*(\widehat r_{\alpha}(\widehat m_R), X_{n+1})] \\
        = & \mathbb{E}_{Q_X}[\widehat m_R(\widehat r_{\alpha}(\widehat m_R), X_{n+1})] - (\mathbb{E}_{Q_X}[\widehat m_R(\widehat r_{\alpha}(\widehat m_R), X_{n+1})] - \mathbb{E}_{Q_X}[m_R^*(\widehat r_{\alpha}(\widehat m_R), X_{n+1})]) \\
        \geq & 1- \alpha - \sup_{\theta \in \mathbb{R}}\left|\mathbb{E}_{Q_X}[\widehat m_R(\theta, X_{n+1})] - \mathbb{E}_{Q_X}[m_R^*(\theta, X_{n+1})] \right| .
    \end{split}
\end{equation*}
Therefore we see that consistent estimation of the $(1-\alpha)$-th quantile of $R(X_{n+1}, Y_{n+1})$ ($ (X_{n+1},Y_{n+1}) \sim Q_X \otimes P_{Y|X}$) is possible if either $w^*(x)$ or $m_R^*(\theta, x)$ is consistently estimated. 

We can reformulate the covariate shift problem as the classical missing data problem. To do this, we define variables $R_i$ for all $i \in [n_1]$. For $1 \leq i \leq n_1$ (i.e.\ for all points in $\mathcal{D}^{cal}_P$) we set $R_i = R(X_i, Y_i)$. For $n_1+1 \leq i \leq n$ (i.e.\ for all points in $\mathcal{D}^{cal}_Q$) we do not observe the response $Y_i$ and hence do not observe $R_i$. The complete information on the $i$-th data point is $\tilde O_i = (X_i , R_i)$. Following the coarsening notation introduced in \Cref{tab:missing_data_type}, the observed calibration data is $\mathcal{D}^{cal} = \{O_i\}_{i =1}^n$ where $O_i = (C_i, \mathfrak G_{C_i}(\tilde O_i))$ for all $i \in [n]$. Here $C_i$ are the random coarsening levels taking values in $\{0, \infty\}$ and $\mathfrak G_0(\tilde O_i) = X_i, \mathfrak G_{\infty}(\tilde O_i) = (X_i, R_i)$. Moreover the distributional assumptions on the points in $\mathcal{D}^{cal}_P$ and $\mathcal{D}^{cal}_Q$ can be re-stated as,
\begin{align*}
    X_i|C_i = \infty \stackrel{d}{=} P_X \quad \mbox{and} \quad X_i|C_i = 0 \stackrel{d}{=} Q_X \quad &\mbox{for all } i \in [n], \\
    R_i|C_i = \infty,X_i = x
\stackrel{d}{=}
R_i|C_i = 0,X_i = x
\stackrel{d}{=}
P_{R|X = x} \quad \mbox{a.e. } x \quad &\mbox{ for all } i \in [n] .
\end{align*}
The above essentially means that $R_i \perp C_i|X_i$. This is the classical missing at random (MAR) assumption.

\section{Multi-stage monotone missingness problem}
\label{app:multi_stage_monotone}

Having understood the application of $\rscp$ in the two-stage monotone missingness problem and analyzed the resulting coverage guarantee, we are now ready to apply $\rscp$ to the general multi-stage monotone missingness problem. Recall that the calibration set $\mathcal{D}^{cal}$ comprises the $n$ i.i.d.\ data points, $\mathcal{D}^{cal} = \{O_1, \cdots, O_n\}$ where $O_i = (C_i, \mathfrak G_{C_i}(\tilde O_i))$ for $i \in [n]$. For the $i$-th data point, the complete data vector is $\tilde O_i = (X_i, Z^1_i,\cdots, Z^D_i, Y_i)$ and the coarsening random variable $C_i$ takes value in the set $\{0, 1,\cdots, D, \infty\}$ where $C_i = 0, k, \infty$ corresponds to the coarsened data $X, (X, Z_1, \cdots, Z_k ), (X, Z, Y)$ respectively for any $k \in [D]$. We assume that the MAR assumption \eqref{eq:mar_monotone} holds. We consider the following definitions of the multi-stage propensity scores,
\begin{equation}
\label{eq:propensity_multi_stage}
    \begin{split}
      \pi_1^*(x) &= \frac{\mathbb{P}(C = 0|X = x)}{\mathbb{P}(C \geq 1|X = x)}, \\
      \pi_{j+1}^*(x, z_1, \cdots, z_j)
      &= \mathbb{P}(C > j| X = x, Z^1 = z_1, \cdots, Z^j = z_j,C \geq j),
    \end{split}
\end{equation}
for $j = 1, \cdots, D$. As in the two-stage monotone missingness problem, all the multi-stage propensity scores defined above can be written in terms of $\omega(\cdot)$ in \eqref{eq:mar_monotone}. It can be easily verified that for $j \in [D]$,
\begin{equation*}
    \begin{split}
      \pi_1^*(x)
      &= \frac{\omega(0, x)}{1 - \omega(0, x)}, \\
      \pi_{j+1}^*(x, z_1, \cdots, z_j)
      &=1 - \frac{\omega(j, (x, z_1, \cdots, z_j))}
      {1 - \omega(0,x) - \sum_{k = 1}^{j-1} \omega(k, (x, z_1, \cdots, z_k))}.
    \end{split}
\end{equation*}
We also consider the following definitions of the outcome regression functions,
\begin{equation}
\label{eq:outcome_multi_stage}
    \begin{split}
        m_{D+1}^*(\theta, x, z_1, \cdots, z_D)
        &=
        \mathbb{P}[Y \leq \theta | X = x, Z ^1 = z_1, \cdots, Z^D = z_D, C = \infty], \\
        m_j^*(\theta, x, z_1, \cdots, z_{j -1})
        &=
        \int m_{j + 1}^*(\theta, x, z_1, \cdots, z_j)
        dF(z_j | x, z_1, \cdots, z_{j - 1}, C \geq j),
    \end{split}
\end{equation}
for $j = 1, \cdots, D$. 
Proceeding analogously to the two-stage case, one can show that the efficient influence function for the $(1-\alpha)$-th quantile of $Y|C=0$ under the multi-stage MAR monotone missingness setup is
\[
\mathrm{IF}
\left(
r_\alpha,
c,
\mathfrak G_c((x,z^1,\ldots,z^D,y)),
\{m_j^*\}_{j=1}^{D+1},
\{\pi_j^*\}_{j=1}^{D+1}
\right),
\]
where the function $\mathrm{IF}(\cdots)$ is defined at
$o = (c, \mathfrak G_c((x, z^1, \cdots, z^D, y)))$ as,
\begin{equation}
\label{eq:influ_Mmdp}
    \begin{split}
        &\mathrm{IF}(\theta, o, \{m_j\}_{j = 1}^{D+1},
        \{\pi_j\}_{j = 1}^{D + 1}) \\
      ={}&
        \frac{\textbf{1}\{c= \infty\} \pi_1(x)}
        {\pi_2(x, z^1)\cdots \pi_{D+1}(x, z^1, \cdots, z^D)}
        (\textbf{1}\{y \leq \theta \} - ( 1- \alpha)) \\
        &-
        \frac{\textbf{1}\{c \geq D\} \pi_1(x)}
        {\pi_2(x, z^1)\cdots \pi_{D}(x, z^1, \cdots, z^{D-1})}
        \left(
        \frac{\textbf{1}\{c = \infty\}}
        {\pi_{D+1}(x, z^1, \cdots, z^D)}
        - 1
        \right)
        [m_{D+1}(\theta, x, z^1, \cdots, z^D)- ( 1- \alpha)] \\
        &- \cdots \\
        &-
        \textbf{1}\{c \geq 1\}\pi_1(x)
        \left(
        \frac{\textbf{1}\{c \geq 2\}}
        {\pi_{2}(x, z^1)}
        - 1
        \right)
        [m_{2}(\theta, x, z^1)- ( 1- \alpha)] \\
        &-
        (\textbf{1}\{c \geq 1\}\pi_1(x) - \textbf{1}\{c = 0\})
        [m_{1}(\theta, x)- ( 1- \alpha)].
    \end{split}
\end{equation}
Prediction sets of the form
\[
    \widehat{\mathcal C}_{n,\alpha}
    =
    \{y\in\mathbb R:y\leq\widehat r_\alpha\}
\]
are constructed using the multi-stage $\rscp$ procedure described in
\Cref{algg:split_robust_Mmdp}. As illustrated in
\Cref{rem:nc_score}, the same construction can be readily generalized
to prediction sets of the form
$\{y\in\mathbb R:R(x,y)\leq\widehat r_\alpha\}$ for a general
non-conformity score $R(x,y)$. We introduce compact notation for the nuisance-estimation errors. Let us define,
\[
    H_0=X,
    \qquad
    H_j=(X,Z^1,\ldots,Z^j),
    \qquad j=1,\ldots,D.
\]
We define the following,
\begin{equation}
\label{eq:multi_stage_Aj}
\begin{split}
A_1
&:=
\|\widehat\pi_1-\pi_1^*\|_2,\\
A_j
&:=
\left\|
\frac{
\widehat\pi_1(X)
\prod_{\ell=2}^{j-1}\pi_\ell^*(H_{\ell-1})
}{
\prod_{\ell=2}^{j-1}\widehat\pi_\ell(H_{\ell-1})
}
\right\|_4
\left\|
\frac{\pi_j^*(H_{j-1})}
{\widehat\pi_j(H_{j-1})}
-1
\right\|_4,
\qquad j=2,\ldots,D+1,
\end{split}
\end{equation}
where an empty product is defined to be one. Additionally for $j=1,\ldots,D+1$, we define,
\begin{equation}
\label{eq:multi_stage_EjDj}
\begin{split}
E_j
&:=
\sup_{\theta\in\mathbb R}
\left\|
m_j^*(\theta,H_{j-1})
-
\widehat m_j(\theta,H_{j-1})
\right\|_2,\quad D_j
:=
\left\|
m_j^*(\widehat r_\alpha,H_{j-1})
-
m_j^*(\widehat r_{\alpha,0},H_{j-1})
\right\|_2.
\end{split}
\end{equation}

We also define the following for each $j=1,\ldots,D$,
\begin{equation}
\label{eq:multi_stage_mdagger}
\widehat m_j^\dagger(\theta,H_{j-1})
:=
\mathbb E
\left[
\widehat m_{j+1}(\theta,H_j)
\,\middle|\,
H_{j-1},C\geq j
\right],
\end{equation}
and
\begin{equation}
\label{eq:multi_stage_regression_remainder}
\mathcal R_m
:=
\sum_{j=1}^{D}
\sup_{\theta\in\mathbb R}
\left\|
\widehat m_j^\dagger(\theta,H_{j-1})
-
\widehat m_j(\theta,H_{j-1})
\right\|_2.
\end{equation}

\begin{theorem}[Coverage guarantee for $\rscp$ under multi-stage monotone missingness]
\label{thm:Mmdp_validity}
Suppose that the MAR condition \eqref{eq:mar_monotone} holds and
$\mathbb P(C=0)>0$. Assume that $ \widehat m_j(\theta,H_{j-1})\in[0,1]$ for $j=1,\ldots,D+1$ and that the estimated propensity scores are bounded away from zero. Suppose that, for some
$\overline\pi\geq1$,
\begin{equation}
\label{eq:multi_stage_weight_bound}
\widehat\pi_1(X)\leq\overline\pi,
\qquad
\frac{\widehat\pi_1(X)}
{\widehat\pi_2(H_1)\cdots\widehat\pi_j(H_{j-1})}
\leq\overline\pi,
\qquad j=2,\ldots,D+1,
\end{equation}
almost surely. Suppose further that $Y$ is continuously distributed
and that, for every $j=1,\ldots,D+1$, the map $ \theta\mapsto
    \widehat m_j(\theta,H_{j-1})$ is right continuous. Then the threshold $\widehat r_\alpha$ returned by
\Cref{algg:split_robust_Mmdp} satisfies the following conditional on both the training splits,
\begin{equation}
\label{eq:multi_stage_simple_coverage}
\begin{split}
&
\left|
\mathbb P
\left(
Y_{n+1}\leq\widehat r_\alpha
\,\middle|\,
C_{n+1}=0
\right)
-(1-\alpha)
\right| \leq
\frac{
\displaystyle
\sum_{j=1}^{D+1}A_j(E_j+D_j)
}{
\mathbb P(C=0)
}
+
\frac{
(D+3)\overline\pi
}{
\mathbb P(C=0)(n+1)
}.
\end{split}
\end{equation}
Suppose, in addition, that there exist $\eta>0$ and positive constants
$\mathfrak C_0,\mathfrak C_1,\mathfrak C_2,\ldots,
\mathfrak C_{D+2}$ such that the following conditions hold on $\mathcal N_\alpha
    =
    [r_\alpha-\eta,r_\alpha+\eta]$: first, for all $\theta,\theta'\in\mathcal N_\alpha$,
\begin{equation}
\label{eq:multi_stage_lipschitz_gprime}
\left\|
g'
\left(
O,\theta,
\{\widehat\pi_j\}_{j=1}^{D},
\{\widehat m_j\}_{j=1}^{D+1}
\right)
-
g'
\left(
O,\theta',
\{\widehat\pi_j\}_{j=1}^{D},
\{\widehat m_j\}_{j=1}^{D+1}
\right)
\right\|_2
\leq
\mathfrak C_0|\theta-\theta'|,
\end{equation}
where $g'$ is the function used in
\Cref{algg:split_robust_Mmdp}. Second, for every
$j=1,\ldots,D+1$,
\begin{equation}
\label{eq:multi_stage_lipschitz_m}
\left\|
m_j^*(\theta,H_{j-1})
-
m_j^*(\theta',H_{j-1})
\right\|_2
\leq
\mathfrak C_{j+1}|\theta-\theta'|,
\qquad
\theta,\theta'\in\mathcal N_\alpha.
\end{equation}
Third, the weighted distribution function is locally well separated:
\begin{equation}
\label{eq:multi_stage_margin}
\begin{split}
&
\left|
P
\left[
\textbf{1}\{C=\infty\}
\frac{
\widehat\pi_1(X)
}{
\widehat\pi_2(H_1)\cdots
\widehat\pi_{D+1}(H_D)
}
\left\{
\textbf{1}\{Y\leq\theta\}
-
\textbf{1}\{Y\leq\theta'\}
\right\}
\right]
\right| \geq
\mathfrak C_1|\theta-\theta'|,
\end{split}
\end{equation}
for $\theta,\theta'\in\mathcal N_\alpha$. We also assume that, for every fixed $\theta$ and
$j=1,\ldots,D+1$,
\begin{equation}
\label{eq:multi_stage_no_fixed_jump}
\mathbb P
\left(
J(\widehat m_j(\cdot,H_{j-1}),\theta)>0
\right)
=0.
\end{equation}
Then, for every $\delta\in(0,1)$, there exists a finite constant
$C_\delta$, depending only on $\delta$, $D$, $\overline\pi$, $\eta$,
and the constants in
\eqref{eq:multi_stage_lipschitz_gprime}--\eqref{eq:multi_stage_margin},
such that, with probability at least $1-\delta$ over the training split used to construct $\widehat r_{\alpha, 0}$,
\begin{equation}
\label{eq:multi_stage_high_prob_coverage}
\begin{split}
&
\left|
\mathbb P
\left(
Y_{n+1}\leq\widehat r_\alpha
\,\middle|\,
C_{n+1}=0
\right)
-(1-\alpha)
\right|\leq
\frac{
(D+3)\overline\pi
}{
\mathbb P(C=0)(n+1)
} +
\frac{C_\delta}{\mathbb P(C=0)}
\sum_{j=1}^{D+1}
A_j
\left\{
E_j
+
\mathcal R_m
+
\frac1{\sqrt n}
+
\frac1{\sqrt{n_2}}
\right\}.
\end{split}
\end{equation}
\end{theorem}
\Cref{thm:Mmdp_validity} implies that if the stagewise regression estimators satisfy the following rate condition,
\[
\sup_{\theta\in\mathbb R}
\left\|
\widehat m_j^\dagger(\theta,H_{j-1})
-
\widehat m_j(\theta,H_{j-1})
\right\|_2
=
O(n_2^{-1/2}),
\qquad
j=1,\ldots,D,
\]
then we have the simplified coverage guarantee,
\begin{equation}
\label{eq:multi_stage_high_prob_coverage_simplified}
\begin{split}
&
\left|
\mathbb P
\left(
Y_{n+1}\leq\widehat r_\alpha
\,\middle|\,
C_{n+1}=0
\right)
-(1-\alpha)
\right|\leq
\frac{
(D+3)\overline\pi
}{
\mathbb P(C=0)(n+1)
}
+
\frac{C_\delta}{\mathbb P(C=0)}
\sum_{j=1}^{D+1}
A_j
\left(
E_j+n^{-1/2}+n_2^{-1/2}
\right).
\end{split}
\end{equation}
The proof of \Cref{thm:Mmdp_validity} follows the same argument as
the proof of \Cref{thm:2mmdp} and is therefore skipped for brevity.  

\begin{remark}[Multiple robustness and finite-sample validity]
\label{rem:multi_stage_multiple_robustness}
As in the two stage monotone missingness case, it is possible to obtain the following sharper lower bound on the coverage guarantee,
\begin{equation}
\label{eq:multi_stage_finite_lower}
\mathbb P
\left(
Y_{n+1}\leq\widehat r_\alpha
\,\middle|\,
C_{n+1}=0
\right)
\geq
1-\alpha
-
\frac{
\displaystyle
\sum_{j=1}^{D+1}
A_j(E_j+D_j)
}{
\mathbb P(C=0)
}.
\end{equation}
Consequently, $\rscp$ has exact finite-sample validity whenever $A_1=A_2=\cdots=A_{D+1}=0$. In particular, this occurs when all the stagewise propensity-score
functions are correctly specified. More generally, the leading nuisance-estimation bias in
\eqref{eq:multi_stage_high_prob_coverage_simplified}, $\sum_{j=1}^{D+1}A_jE_j$, vanishes whenever, for every $j=1,\ldots,D+1$, we have either $A_j = 0$ or $E_j = 0$. There are therefore
$2^{D+1}$ combinations of nuisance components for which the leading term in the coverage slack vanishes. This is the $2^{D+1}$-robustness
property of the multi-stage monotone missingness construction.

It is useful to distinguish this multiple robustness property from
exact finite-sample validity.  If an outcome-regression component is
correctly specified at a stage while its corresponding propensity
score is not, the product $A_jE_j$ vanishes, but the
coverage slack can still contain the contribution $ A_j
    (
    n^{-1/2}+n_2^{-1/2}
    )$ from the estimation of the initial root. Hence, in any of these
$2^{D+1}$ configurations, the remaining slack is, up to
constants, $ (
    \sum_{j=1}^{D+1}A_j
    )
    \left(
    n^{-1/2}+n_2^{-1/2}
    \right)
    +
    O(n^{-1})$ under additional regularity conditions. When all the propensity-score components are correctly specified,
all the $A_j$'s vanish, recovering exact finite-sample validity.
\end{remark}

\section{Non-monotone missing data problem}
\label{sec:non_monotone}

In practice, missingness of data may not be monotone, i.e.\ there may be no nested missingness pattern. In other words unlike the monotone missingness case, the missingness of one of the intermediate covariates $Z^j$ in a longitudinal study does not imply that the future covariates $\{Z_{j'}\}_{j' = j+1}^{D}$ and the response of interest $Y$ are unobserved. A typical example of non-monotone missingness is a longitudinal cohort study in which a participant is absent for one or more interim visits but re-enters the study at later time points. Semiparametric theory for non-monotone missingness problems has extensively been studied under the missing at random (MAR) assumption in \cite{robins1994estimation, robins1997non, tsiatis2006semiparametric, sun2018inverse2}, and is known to present several new challenges for estimation and inference. In the following, we provide a general framework for predictive inference under nonmonotone MAR, however, for the most part we do not directly address some of these inferential and implementation challenges, which we leave to future work.

As in \Cref{sec:monotone}, the complete data vector $\tilde O = (X, Z, Y)$ comprises the baseline covariates $X$, the intermediate covariates $Z = (Z^1, \cdots, Z^D)$, and the response of interest $Y$. The observed data vector is $O = (C, \mathfrak G_C(\tilde O))$ where $C$ is a random coarsening level takes values in $\mathcal K = \{0,1, \cdots, K, \infty\}$ for some $K \leq 2^{D+1} - 2$. There exists a collection of subsets $\{\mathscr S_k\}_{k \in [K]}$ satisfying $X \subset \mathscr S_k \subset (X, Z, Y)$ for all $k \in [K]$ such that the coarsening levels $C = 0, k, \infty$ correspond to the coarsened data $X, \mathscr S_k, (X, Z, Y)$ respectively for all $k \in [K]$. It is easy to see that the monotone missing data structure is a special case of this set-up. For instance if $K = D$ and $\mathscr S_k = (X, Z^1, \cdots, Z^k)$ for $k \in [D]$, we are back in the monotone missing data problem. We make the following coarsening at random assumption (\cite{robins1994estimation}, \cite{tsiatis2006semiparametric}),
\begin{equation}
\label{eq:car}
    \mathbb{P}(C = k|X, Z, Y)
    =
    \omega(k,  \mathfrak{G}_k((X,Z, Y)))
    \quad
    \mbox{for $k = 0, 1, \cdots, K, \infty$}.
\end{equation}
We also consider the following definition for all $k \in [K]$,
\begin{equation*}
    \begin{split}
        \omega^*(k, \mathfrak{G}_k((X,Z, Y)))
        \coloneqq&
        \mathbb{P}(C = k|X, Z, Y, C \geq 1) \\
        =&
        \frac{\mathbb{P}(C = k|X, Z, Y)}
        {\mathbb{P}(C \geq 1 |X, Z, Y)}
        =
        \frac{\omega(k,  \mathfrak G_k((X,Z, Y)))}
        {1 - \omega(0, X)}.
    \end{split}
\end{equation*}
The efficient influence function for $r_{\alpha}$, the $(1 - \alpha)$-th quantile of $Y|C = 0$, has been derived in Chapter 10 of \cite{tsiatis2006semiparametric}. We introduce some new definitions before stating the influence function. Suppose $h(\cdots)$ denotes any real valued function on the full data $\tilde O = (X, Z, Y)$. Consider the following real valued linear operator,
\[
\mathcal{L}\{ h(X, Z, Y)\}
=
\sum_{k = 1}^{\infty}
\textbf{1}\{C = k\}
\mathbb{E}[h(X, Z, Y) |\mathfrak G_k((X,Z, Y))].
\]
We consider another linear operator $\mathcal M$ defined on $h(\cdots)$,,
\[
\mathcal{M}\{ h(X, Z, Y)\}
=
\mathbb{E}[\mathcal{L}\{h\}|X, Z, Y, C \geq 1]
=
\sum_{k = 1}^{\infty}
\omega^*(k, \mathfrak G_k((X,Z, Y)))
\mathbb{E}[h(X, Z, Y) | \mathfrak G_k((X,Z, Y))].
\]
Let $d(r_{\alpha}, X, Z, Y)$ be the result of applying the inverse operator $\mathcal{M}^{-1}$ on the function of observed data
\[
h(X, Z, Y)
=
\textbf{1}\{ Y \leq r_{\alpha}\} - (1 - \alpha),
\]
i.e.\
\[
d(r_{\alpha}, X, Z, Y)
=
\mathcal{M}^{-1}
\{
\textbf{1}\{ Y \leq r_{\alpha}\} - (1 - \alpha)
\}.
\]
The efficient influence function [\cite{robins1994estimation, tsiatis2006semiparametric}] of the $(1 - \alpha)$-th quantile of $Y|C = 0$ under \eqref{eq:car} is given below,
\begin{equation}
\label{eq:influence_nonmonotone}
\begin{split}
    \mathrm{IF}(r_{\alpha}, O)
    ={}&
    \textbf{1}\{C \geq 1\}\pi_1^*(X)
    \Bigg[
    \frac{\textbf{1}\{C = \infty\}}
    {\omega^*(\infty, X, Z, Y)}
    (\textbf{1}\{Y \leq r_{\alpha}\} - (1 - \alpha))
    \Bigg.
    \\
    &-
    \frac{\textbf{1}\{C = \infty\}}
    {\omega^*(\infty, X, Z, Y)}
    \Bigg[
    \sum_{k = 1}^K
    \omega^*(k, \mathfrak G_k((X,Z, Y)))
    \mathbb{E}[d(r_{\alpha}, X, Z, Y)
    | \mathfrak G_k((X,Z, Y))]
    \Bigg]
    \\
    &\Bigg.
    +
    \Bigg[
    \sum_{k =1}^K
    \textbf{1}\{C = k\}
    \mathbb{E}[d(r_{\alpha}, X, Z, Y)
    | \mathfrak G_k((X, Z, Y))]
    \Bigg]
    -
    [m_1^*(r_{\alpha}, X) - (1 - \alpha)]
    \Bigg]
    \\
    &+
    \textbf{1}\{C = 0\}
    [m_1^*(r_{\alpha}, X) - (1 - \alpha)] ,
\end{split}
\end{equation}
where the nuisance functions $\pi_1^*(x), m_1^*(\theta, x)$ have the same definitions as in \Cref{sec:cov_shift}. We have suppressed the dependence of the influence function $\mathrm{IF}(\cdots)$ on the nuisance functions $\{\omega^*(k,\mathfrak G_k((X, Z, Y)))\}_{k = 1}^{\infty}$ and $\{\mathbb{E}[d(\theta, X, Z, Y)| \mathfrak G_k((X, Z, Y))]\}_{k= 1}^K$ to keep the notations compact. The above influence function has a double robustness property in the sense that if either all the outcome regression functions $m_1^*(r_{\alpha}, X)$, $\{\mathbb{E}[d(r_{\alpha}, X, Z, Y)| \mathfrak G_k((X, Z, Y))]\}_{k = 1}^{K}$ or all the propensity scores $\pi_1^*(X), \{\omega^*(k, \mathfrak G_k((X, Z, Y)))\}_{k = 1}^{\infty}$ are correctly specified then $\mathbb{E}[\mathrm{IF}(r_{\alpha}, O)] = 0$. The following lemma explicitly states this property of the efficient influence function.

\begin{lemma}
\label{lem:non_monotone_dr}
The following equation holds true,
\begin{equation*}
\begin{split}
& \mathbb{E} \Bigg[
\textbf{1}\{C \geq 1\}\widehat \pi_1(X)
\Bigg[
\frac{\textbf{1}\{C = \infty\}}
{\widehat \omega(\infty, X, Z, Y)}
(\textbf{1}\{Y \leq r_{\alpha}\} - (1 - \alpha))
\Bigg.
\\
&-
\frac{\textbf{1}\{C = \infty\}}
{\widehat \omega(\infty, X, Z, Y)}
\Bigg[
\sum_{k = 1}^K
\widehat \omega(k,  \mathfrak G_k((X, Z, Y)))
\widehat{\mathbb{E}}
[d(r_{\alpha}, X, Z, Y)|\mathfrak G_k((X, Z, Y))]
\Bigg]
\\
&\Bigg.
+
\Bigg[
\sum_{k =1}^K
\textbf{1}\{C = k\}
\widehat{\mathbb{E}}
[d(r_{\alpha}, X, Z, Y)|\mathfrak G_k((X, Z, Y))]
\Bigg]
-
[\widehat m_1(r_{\alpha}, X) - (1 - \alpha)]
\Bigg]
\\
&+
\Bigg.
\textbf{1}\{C = 0\}
[\widehat m_1(r_{\alpha}, X) - (1 - \alpha)]
\Bigg]
=0,
\end{split}
\end{equation*}
if either all the propensity score functions are known, i.e.\
\[
\widehat \pi_1(X)=\pi_1^*(X),
\qquad
\{\widehat \omega(k,\mathfrak G_k((X,Z,Y)))\}_{k=1}^{\infty}
=
\{\omega^*(k,\mathfrak G_k((X,Z,Y)))\}_{k=1}^{\infty},
\]
or all the outcome regression functions are known, i.e.\
\[
\widehat m_1(\theta,x)=m_1^*(\theta,X),
\]
and
\[
\{
\widehat{\mathbb E}
[d(\theta,X,Z,Y)|\mathfrak G_k((X,Z,Y))]
\}_{k=1}^{K}
=
\{
\mathbb E
[d(\theta,X,Z,Y)|\mathfrak G_k((X,Z,Y))]
\}_{k=1}^{K}.
\]
\end{lemma}

The proof of \Cref{lem:non_monotone_dr} is discussed in
\Cref{app:proof_dr_non_monotone}. Given estimates
$\widehat \pi_1(X)$,
$\{\widehat \omega(k,\mathfrak G_k((X,Z,Y)))\}_{k=1}^{\infty}$
of
$\pi_1^*(X)$,
$\{\omega^*(k,\mathfrak G_k((X,Z,Y)))\}_{k=1}^{\infty}$
and estimates
$\widehat m_1(\theta,x)$,
$\{\widehat{\mathbb E}[d(\theta,X,Z,Y)|
\mathfrak G_k((X,Z,Y))]\}_{k=1}^{K}$
of
$m_1^*(\theta,X)$,
$\{\mathbb E[d(\theta,X,Z,Y)|
\mathfrak G_k((X,Z,Y))]\}_{k=1}^{K}$ respectively, we define
$\widehat{\mathrm{IF}}(\theta,O_i)$ for all data points $i\in[n]$ in the
calibration data set as,
\begin{equation}
\label{eq:influence_nonmonotone_estimated}
\begin{split}
\widehat{\mathrm{IF}}(\theta,O_i)
={}&
\textbf{1}\{C_i\geq1\}\widehat\pi_1(X_i)
\Bigg[
\frac{\textbf{1}\{C_i=\infty\}}
{\widehat\omega(\infty,X_i,Z_i,Y_i)}
(\textbf{1}\{Y_i\leq\theta\}-(1-\alpha))
\Bigg.
\\
&-
\frac{\textbf{1}\{C_i=\infty\}}
{\widehat\omega(\infty,X_i,Z_i,Y_i)}
\Bigg[
\sum_{k=1}^{K}
\widehat\omega(k,\mathfrak G_k((X_i,Z_i,Y_i)))
\widehat{\mathbb E}
[d(\theta,X_i,Z_i,Y_i)|
\mathfrak G_k((X_i,Z_i,Y_i))]
\Bigg]
\\
&\Bigg.
+
\Bigg[
\sum_{k=1}^{K}
\textbf{1}\{C_i=k\}
\widehat{\mathbb E}
[d(\theta,X_i,Z_i,Y_i)|
\mathfrak G_k((X_i,Z_i,Y_i))]
\Bigg]
-
[\widehat m_1(\theta,X_i)-(1-\alpha)]
\Bigg]
\\
&+
\textbf{1}\{C_i=0\}
[\widehat m_1(\theta,X_i)-(1-\alpha)].
\end{split}
\end{equation}

Prediction sets of the form
\[
\widehat{\mathcal C}_{n,\alpha}
=
\{y\in\mathbb R:y\leq\widehat r_\alpha\}
\]
are generated using $\rsp$, as detailed in
\Cref{algg:split_robust_nmdp} in \Cref{app:algo}. As in the earlier
sections, this construction can be modified to incorporate a general
non-conformity score $R(x,y)$.

The coverage guarantee for the non-monotone missingness problem follows directly from the general $\rsp$ result in \Cref{thm:gen_influ_func}, together with the
conditional on a coarsening-pattern argument in
\Cref{rem:general_pattern_conditional}.

To see this, let $p_0:=\mathbb P(C=0)$ and define the normalized observed-data function,
\begin{equation}
\label{eq:nmdp_f_normalized}
f_{\mathrm{nm}}^*(O,\theta)
:=
(1-\alpha)
+
\frac{
\mathrm{IF}(\theta,O)
}{
p_0
}.
\end{equation}
It is easy to check that, 
\[
P f_{\mathrm{nm}}^*(O,\theta)
=
\mathbb P(Y\leq\theta\mid C=0).
\]
We similarly define the following using the estimated nuisance functions,
\[
\widehat f_{\mathrm{nm}}(O,\theta)
:=
(1-\alpha)
+
\frac{
\widehat{\mathrm{IF}}(\theta,O)
}{
p_0
}.
\]
The $\rsp$ estimating equation based on
$\widehat f_{\mathrm{nm}}$ is $(1/n) \sum_{i=1}^{n}
\widehat f_{\mathrm{nm}}(O_i,\theta)
\geq
1-\alpha$ which is equivalent to $(1/n)\sum_{i=1}^{n}
\widehat{\mathrm{IF}}(\theta,O_i)
\geq0 $. Thus the threshold in \Cref{algg:split_robust_nmdp} is exactly the
$\rsp$ threshold associated with the normalized influence function
representation \eqref{eq:nmdp_f_normalized}. Therefore if the threshold-indexed class $\{\widehat f_{\mathrm{nm}}(O,\theta):\theta\in\mathbb R\}$
satisfies the bounded-variation condition in
Assumption~\ref{assump:gen_bv}, and the corresponding
centered influence function satisfies the second-moment condition in
Assumption~\ref{assump:gen_local_variance}, with variance-inflation
term $d'_{\mathrm{nm}}$, then \Cref{thm:gen_influ_func} gives the following coverage guarantee conditional on the training data,
\begin{equation}
\label{eq:nmdp_general_coverage}
\begin{split}
&
\left|
\mathbb P
\left(
Y_{n+1}\leq\widehat r_\alpha
\,\middle|\,
C_{n+1}=0
\right)
-(1-\alpha)
\right| \leq
d_{\mathrm{nm}}
+
\mathfrak C_n
\left(
\sqrt{d_{\mathrm{nm}}}+n^{-1/4}
\right)
\left(
n^{-1/2}
+
n^{-1/4}\sqrt{d'_{\mathrm{nm}}}
\right),
\end{split}
\end{equation}
where $\mathfrak C_n=\mathfrak C\log(en)$ and $d_{\mathrm{nm}}
:=
(1/{p_0})
\sup_{\theta\in\mathbb R}
|
P
[
\widehat{\mathrm{IF}}(\theta,O)
-
\mathrm{IF}(\theta,O)
]
|$. The bias term $d_{\mathrm{nm}}$ inherits the double
robustness of the efficient influence function. In particular, the
same calculation underlying \Cref{lem:non_monotone_dr}, applied
uniformly in $\theta$, shows that $d_{\mathrm{nm}}=0$ if either all
the propensity-score components $\pi_1^*,
\{\omega^*(k,\mathfrak G_k)\}_{k=1}^{\infty}$ are correctly specified, or all the outcome-regression components $m_1^*,
\{
\mathbb E[d(\theta,X,Z,Y)|\mathfrak G_k]
\}_{k=1}^{K}$ are correctly specified.  

The procedure developed above is the $\rsp$ version of the non-monotone construction. Consequently, even when the nuisance-estimation errors vanish
($d_{\mathrm{nm}}=d'_{\mathrm{nm}}=0$), its coverage slack retains the empirical-process remainder
$O(\log n/n^{3/4})$ appearing in \Cref{thm:gen_influ_func}.

The $\rscp$ construction can also be applied under non-monotone MAR. In particular, in addition to freezing
$\widehat m_1(\theta,X)$ at an independent preliminary root
$\widehat r_{\alpha,0}$, we freeze each of the threshold-dependent
operator-derived regression functions $
\widehat{\mathbb E}
[
d(\theta,X,Z,Y)
\,|\,
\mathfrak G_k((X,Z,Y))
]$ (for $k=1,\ldots,K$) at $\theta=\widehat r_{\alpha,0}$. After these nuisance functions are frozen, the only remaining dependence of the estimated influence function on the final threshold $\theta$ is through,
\[
\frac{
\mathbf 1\{C=\infty\}\widehat\pi_1(X)
}{
\widehat\omega(\infty,X,Z,Y)
}
\mathbf 1\{Y\leq\theta\},
\]
which is non-decreasing in $\theta$ whenever the estimated propensity weights are nonnegative. Hence the resulting estimating function satisfies the monotonicity requirement underlying the general $\rscp$ construction in \Cref{algg:split_robust_three_split}. Under the corresponding nuisance-bias, range, and regularity conditions, the general coverage result in \Cref{thm:gen_three_split} can therefore be applied. For ease of exposition, we restrict the detailed development and the numerical experiments for non-monotone missingness to the $\rsp$ construction. 
\section{Additional simulations}
\label{app:add_simul}

We next compare $\rsp$ and $\rscp$ when the response and missingness indicators depend on latent Gaussian covariates, but the analyst observes only nonlinear transformations of those covariates. For each unit $i\in\mathcal D$, independently generate the following variables.

\textbf{Setting $2$:}
\begin{equation}
\label{eq:sim_dgp_mis}
    \begin{split}
        X_i',Z_{1i}',\epsilon_{1i}',\epsilon_{2i}'
        \stackrel{\mathrm{iid}}{\sim}{}&\mathcal N(0,1),\\
        Z_{2i}'={}&2X_i'+Z_{1i}'+0.2\epsilon_{1i}',\\
        X_i={}&\exp(X_i'/2),\\
        Z_{1i}={}&(Z_{1i}'Z_{2i}'+0.6)^3,\\
        Z_{2i}={}&(Z_{1i}'+Z_{2i}')^2,\\
        Y_i={}&2X_i'+2Z_{1i}'+0.6Z_{2i}'+0.2\epsilon_{2i}',\\
        T_i\mid X_i' \sim{}&\operatorname{Bernoulli}
        \!\left\{\operatorname{expit}(0.2X_i')\right\},\\
        S_i\mid X_i',Z_{1i}',Z_{2i}',T_i \sim{}&
        \begin{cases}
            \operatorname{Bernoulli}\!\left\{
            \operatorname{expit}(0.1X_i'+0.1Z_{1i}'+0.05Z_{2i}')
            \right\}, & T_i=0,\\
            0, & T_i=1.
        \end{cases}
    \end{split}
\end{equation}
As in Setting~$1$, $T_i=1$ corresponds to $C_i=0$, $T_i=0,S_i=0$ corresponds to $C_i=1$, and $T_i=0,S_i=1$ corresponds to $C_i=\infty$. Thus the observed data retain the two-stage monotone pattern, but the analyst observes $(X_i,Z_{1i},Z_{2i})$ rather than the latent covariates $(X_i',Z_{1i}',Z_{2i}')$ that generate $Y_i,T_i,$ and $S_i$. For example, logistic regression of $T_i$ on $X_i'$ and linear regression of $Y_i$ on $(X_i',Z_{1i}',Z_{2i}')$ are correctly specified, whereas the corresponding working models based on $(X_i,Z_{1i},Z_{2i})$ are misspecified.

We use the same score-learning, nuisance-estimation, calibration, and testing protocol as in \Cref{subsec:simulation}: the target level is $0.9$, the calibration and independent initialization sample sizes are $n=25,100,400$, coverage and mean interval width are evaluated using $1500$ independently generated observations with $C=0$, and each experiment is repeated $500$ times. The $\rsp$ method evaluates the full efficient-influence-function equation at every candidate threshold. The $\rscp$ method instead freezes its threshold-dependent nuisance functions at an independent preliminary threshold and applies the finite-sample correction in \Cref{alg:split_robust_mmdp}.

\begin{figure}[ht]
  \centering
  \includegraphics[width=\textwidth]{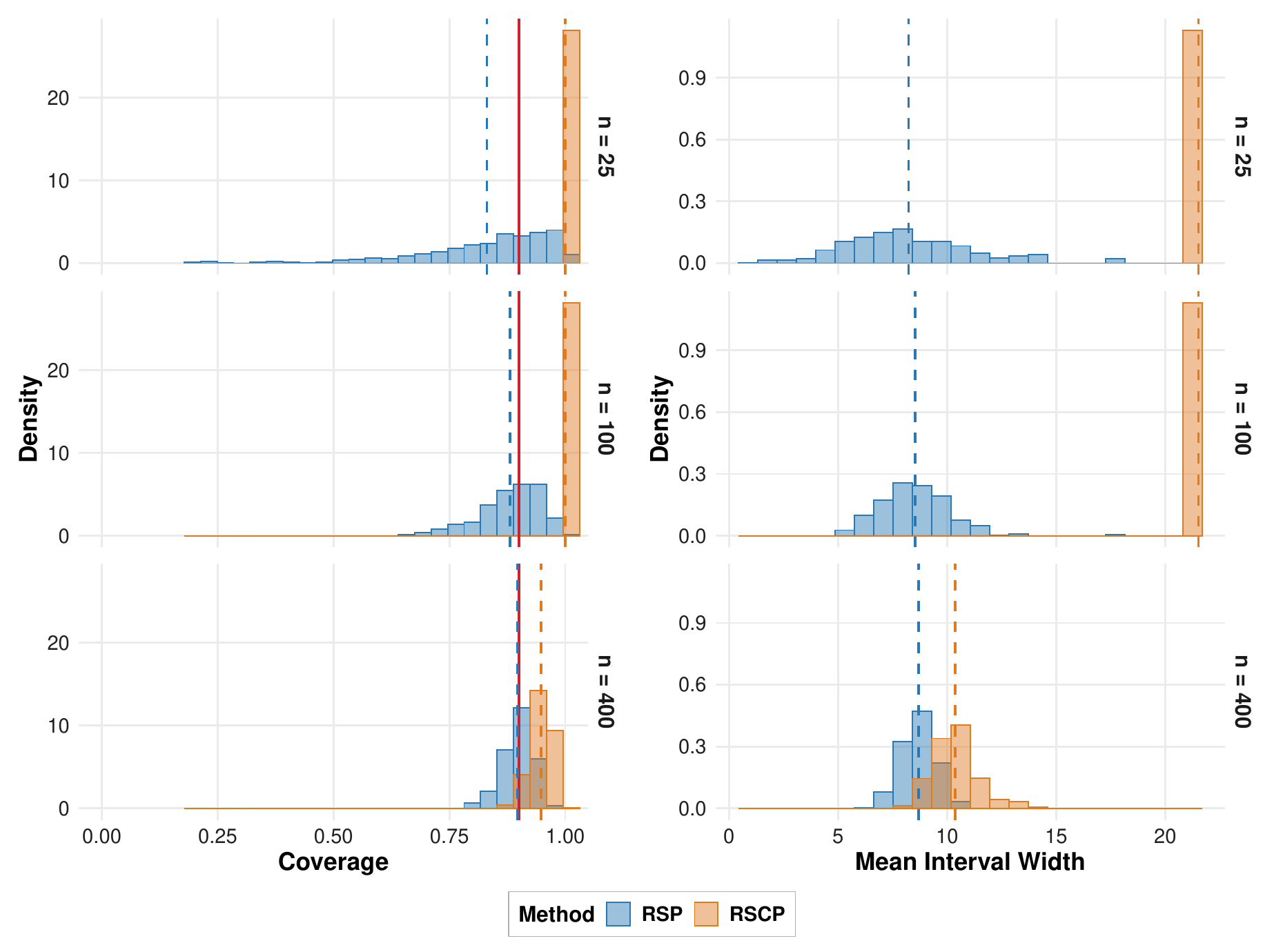}
  \caption{Histograms of coverage (left) and mean interval width (right) for $\rsp$ and $\rscp$ under Setting~$2$, using the absolute-residual non-conformity score. Rows correspond to sizes of calibration data set $n=25,100,400$; colored dashed lines mark method-specific means and the solid red line marks nominal coverage $0.9$.}
  \label{fig:sim_setting2_absolute}
\end{figure}

\begin{figure}[ht]
  \centering
  \includegraphics[width=\textwidth]{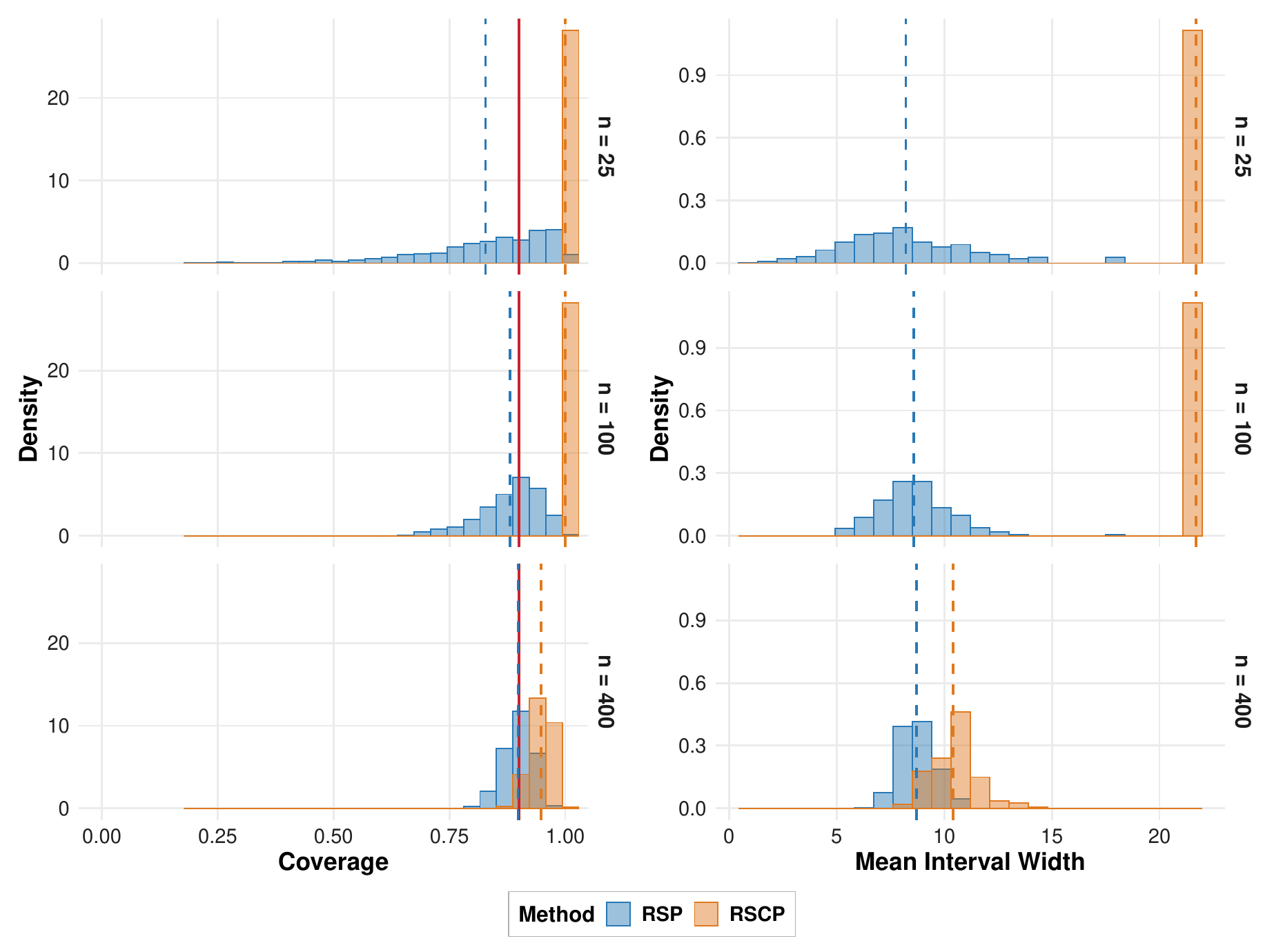}
  \caption{Histograms of coverage (left) and mean interval width (right) for $\rsp$ and $\rscp$ under Setting~$2$, using the CQR non-conformity score. Rows correspond to different sizes of calibration data set $n=25,100,400$; colored dashed lines mark method-specific means and the solid red line marks nominal coverage $0.9$.}
  \label{fig:sim_setting2_cqr}
\end{figure}

The results in \Cref{fig:sim_setting2_absolute,fig:sim_setting2_cqr} are qualitatively similar to those under Setting~$1$.
\begin{enumerate}
    \item The mean coverage of $\rsp$ improves with calibration size despite the misspecified observed covariates. For the absolute-residual score, the values at $n=25,100,400$ are $0.831$, $0.881$, and $0.898$; for CQR they are $0.828$, $0.881$, and $0.898$.
    \item The $\rscp$ correction again yields conservative intervals at small $n$, with mean coverage $1.000$ at $n=25$ and $n=100$. At $n=400$, the mean coverage decreases to $0.947$ and $0.948$ for the absolute-residual and CQR scores, respectively.
    \item The two scores have comparable efficiency. At $n=400$, the mean widths are $8.685$ and $10.382$ for $\rsp$ and $\rscp$, respectively, with the absolute-residual score, and $8.718$ and $10.422$ with CQR. As in Setting~$1$, corrected equations that do not cross zero on the numerical grid are represented by its upper endpoint, producing the right-edge mass in the $\rscp$ width histograms for small sizes of the calibration data.
\end{enumerate}


We now retain the latent Gaussian outcome model and nonlinear observed covariates of Setting~$2$, but replace monotone dropout by a non-monotone MAR mechanism. Specifically, independently for each unit $i$, 

\textbf{Setting $3$:}
\begin{equation}
\label{eq:sim_dgp_nonmonotone_mis}
\begin{split}
X_i',Z_{1i}',\epsilon_{1i}',\epsilon_{2i}'
\stackrel{\mathrm{iid}}{\sim}{}&\mathcal N(0,1),\\
Z_{2i}'={}&2X_i'+Z_{1i}'+0.2\epsilon_{1i}',\\
X_i={}&\exp(X_i'/2),\\
Z_{1i}={}&(Z_{1i}'Z_{2i}'+0.6)^3,\\
Z_{2i}={}&(Z_{1i}'+Z_{2i}')^2,\\
Y_i={}&2X_i'+2Z_{1i}'+0.6Z_{2i}'+0.2\epsilon_{2i}'.
\end{split}
\end{equation}
Thus, as in Setting~$2$, the analyst observes $(X_i,Z_{1i},Z_{2i})$ rather than the latent covariates $(X_i',Z_{1i}',Z_{2i}')$ entering the response model. We first generate,
\[
T_i\mid X_i'
\sim\operatorname{Bernoulli}\{\operatorname{expit}(0.2X_i')\}
=\operatorname{Bernoulli}\{\operatorname{expit}(0.4\log X_i)\},
\]
and set $C_i=0$ when $T_i=1$, so only $X_i$ is revealed. Conditional on $T_i=0$, define the bounded summaries,
\[
U_{Xi}=\tanh(\log X_i),\qquad
U_{1i}=\tanh\!\left\{\frac{\operatorname{sgn}(Z_{1i})|Z_{1i}|^{1/3}}{2}\right\},
\qquad
U_{2i}=\tanh\!\left\{\frac{\log(1+Z_{2i})}{2}\right\},
\qquad
U_{Yi}=\tanh(Y_i/4),
\]
and the pattern-specific probabilities,
\begin{align*}
\rho_1(X_i,Z_{1i})
&=\operatorname{expit}(-2.40+0.25U_{Xi}+0.35U_{1i}),\\
\rho_2(X_i,Z_{2i})
&=\operatorname{expit}(-2.40-0.25U_{Xi}+0.35U_{2i}),\\
\rho_3(X_i,Z_{1i},Y_i)
&=\operatorname{expit}(-2.50+0.20U_{Xi}-0.20U_{1i}+0.30U_{Yi}),\\
\rho_\infty(X_i,Z_{1i},Z_{2i},Y_i)
&=1-\rho_1-\rho_2-\rho_3.
\end{align*}
The bounded summaries ensure $\rho_\infty>0.57$. We draw $C_i\in\{1,2,3,\infty\}$ from these probabilities conditional on $T_i=0$, where,
\[
\mathfrak G_1=(X_i,Z_{1i}),\qquad
\mathfrak G_2=(X_i,Z_{2i}),\qquad
\mathfrak G_3=(X_i,Z_{1i},Y_i),\qquad
\mathfrak G_\infty=(X_i,Z_{1i},Z_{2i},Y_i).
\]
These patterns are not nested: for example, pattern $2$ records $Z_{2i}$ without $Z_{1i}$, while pattern $3$ records $Y_i$ without $Z_{2i}$. Moreover, each $\rho_k$ depends only on variables contained in $\mathfrak G_k$, so the construction satisfies the pattern-specific non-monotone MAR condition used by \citet{sun2018inverse2}.

We estimate the nuisance functions using the transformed covariates. In particular, the $C=0$ propensity is fitted as a polynomial logistic regression in raw $X$, while the pattern probabilities are fitted by the observed-data likelihood using standardized raw components of the variables available under each pattern. These models do not contain the logarithmic, cube-root, or hyperbolic-tangent transformations in the generating mechanism and are therefore misspecified. Conditional-expectation operators are approximated by inverse-probability-weighted nearest-neighbour smoothers on complete nuisance-training observations. We form a finite-dimensional approximation to $\mathcal M$, solve for
\[
\widehat d(\theta)
=\widehat{\mathcal M}^{-1}
\{\mathbf 1(R\leq\theta)-(1-\alpha)\},
\]
substitute the resulting nuisance estimates into the non-monotone efficient influence function, and take the $\rsp$ threshold to be the smallest $\theta$ for which its calibration-sample mean is nonnegative. We use the absolute-residual score of \citet{lei2018distribution} and the CQR score of \citet{romano2019conformalized}.

\begin{figure}[ht]
  \centering
  \includegraphics[width=\textwidth]{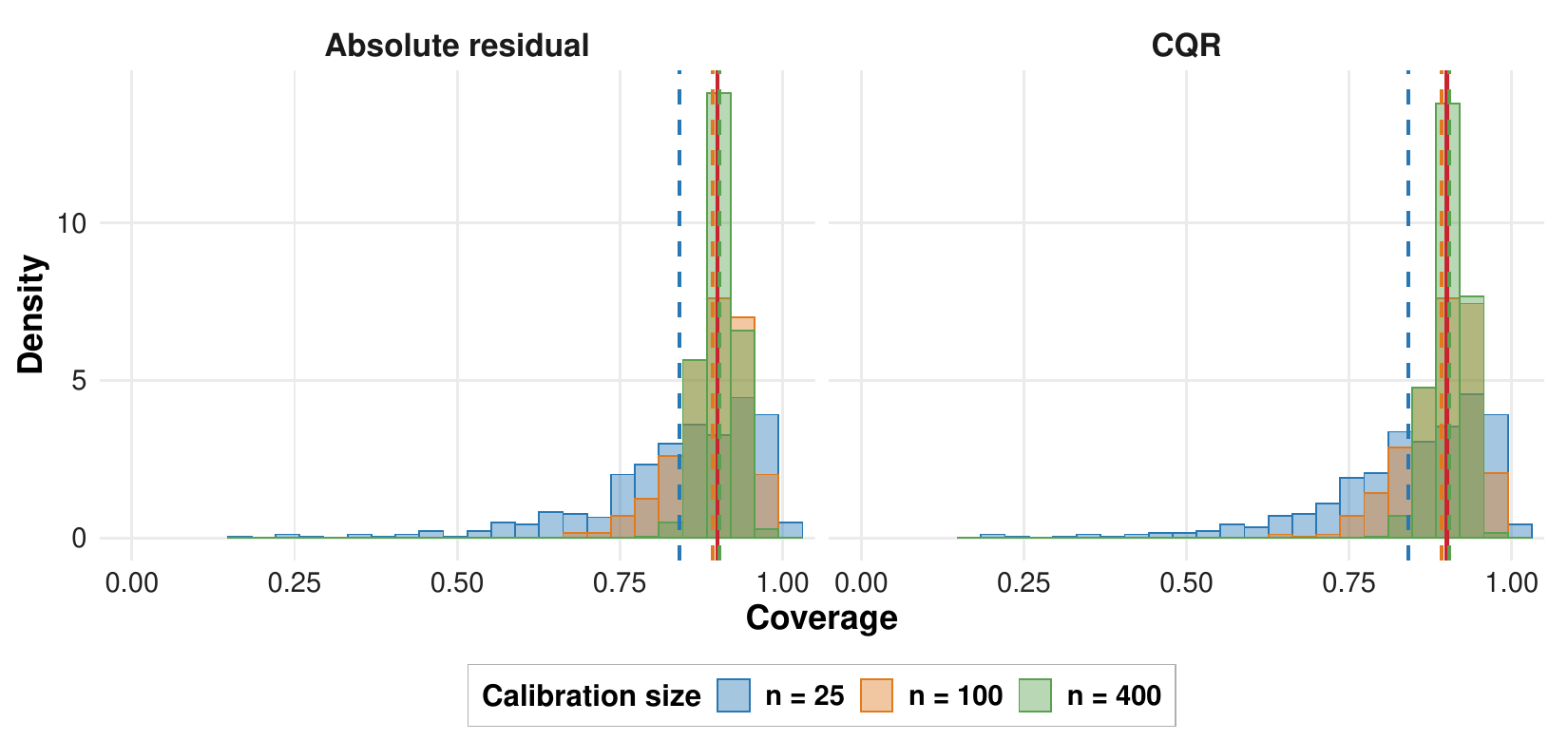}
  \caption{Coverage histograms for $\rsp$ under non-monotone MAR with misspecified covariates (Setting~$3$). The left and right panels use the absolute-residual and CQR non-conformity scores, respectively; colors identify the sizes of the calibration data set $n=25,100,400$.}
  \label{fig:sim_setting3_nonmonotone_mis}
\end{figure}

We set $1-\alpha=0.9$ and consider calibration sizes $n=25,100,400$. Coverage is evaluated conditional on $C=0$ using $1500$ independently generated target observations in each repetition, and the experiment is repeated $500$ times. The Monte Carlo repetitions are run in parallel. In \Cref{fig:sim_setting3_nonmonotone_mis}, the solid red line marks nominal coverage $0.9$, while the colored dashed lines mark the mean coverage for each calibration size.

The results show that coverage improves and becomes more concentrated as the size of the calibration data set increases. For the absolute-residual score, the mean coverages at $n=25,100,400$ are $0.841$, $0.893$, and $0.903$, with corresponding standard deviations $0.133$, $0.055$, and $0.027$. For CQR, the mean coverages are $0.842$, $0.891$, and $0.903$, with standard deviations $0.132$, $0.056$, and $0.027$. Thus, although both scores exhibit appreciable lower-tail variability at $n=25$, their coverage distributions are centered close to the nominal level by $n=100$ and concentrate tightly around $0.9$ at $n=400$. The similarity of the two panels also indicates that the observed robustness is not sensitive to the choice between the absolute-residual and CQR scores in this experiment.

\section{Proof of \Cref{thm:influence_old_cmdp}}
\label{app:thm_old_influ_cmdp}
Recall the definition of $\widehat r_{\alpha}$ in \eqref{eq:dr_cmdp}, 
\[
\widehat r_{\alpha} = \inf\left\{ \theta: \frac{1}{n} \sum_{i = 1}^n f(O_i, \theta, \widehat \pi, \widehat m_R) \geq 0 \right\},
\]
where, 
\[
f(O, \theta, \pi, m_R) = \textbf{1}\{C = \infty\} \pi(X) [\textbf{1}\{R \leq \theta \} - m_R(\theta, X)] + \textbf{1}\{C = 0\} [m_R(\theta, X) - (1 - \alpha)].
\]
From here-on we use $Pf(O, \theta, \pi, m_R) = (1/n) \sum_{i = 1}^n f(O_i, \theta, \pi, m_R)$. Following \cite{yang2024doubly} we have the following regarding the coverage guarantee, 
\begin{equation*}
    \begin{split}
        \mathbb P(R_{n+1} \leq \widehat r_{\alpha} |C_{n+1} = 0) = & 1- \alpha + \frac{\mathbb E[Pf (O,\widehat r_{\alpha}, \pi^*, m_R^* )]}{\mathbb P(C = 0)} \\
        \geq & 1 - \alpha + \frac{\mathbb E[Pf (O,\widehat r_{\alpha}, \widehat \pi,\widehat m_R )]}{\mathbb P(C = 0)} - \frac{\| \widehat \pi - \pi^*\|_2 \sup_{\theta}\| \widehat m_R(\theta, \cdot) - m_R^*(\theta, \cdot)\|_2}{\mathbb P(C = 0)} .
    \end{split}
\end{equation*}
The last inequality follows from Theorem $2$ of \cite{yang2024doubly},
\begin{equation}
    \label{eq:yang_dr}
    \sup_{\theta \in \mathbb R}|P[f(O, \theta, \widehat \pi, \widehat m_R) - f(O, \theta, \pi^*, m_R^*)]| \leq \| \widehat \pi - \pi^*\|_2 \sup_{\theta \in \mathbb R}\| \widehat m_R(\theta, \cdot) - m_R^*(\theta, \cdot)\|_2
\end{equation}
In the above equation, $R_{n+1}$ is used to define the random variable $R(X_{n+1}, Y_{n+1})$. From the definition of $\widehat r_{\alpha}$ we have the following, 
\begin{equation*}
    \begin{split}
        & \mathbb E[Pf (O,\widehat r_{\alpha}, \widehat \pi,\widehat m_R )]  \\
        \geq & \mathbb E[(P - P_n)f (O,\widehat r_{\alpha}, \widehat \pi,\widehat m_R )] \\
        = & \mathbb E[(P - P_n)[f (O,\widehat r_{\alpha}, \widehat \pi,\widehat m_R ) - f (O,r_{\alpha}, \widehat \pi,\widehat m_R )]] \\
        = & \mathbb E\left[\frac{(P - P_n)[f (O,\widehat r_{\alpha}, \widehat \pi,\widehat m_R ) - f (O,r_{\alpha}, \widehat \pi,\widehat m_R )] }{\sqrt{|P[f (O,\widehat r_{\alpha}, \widehat \pi,\widehat m_R ) - f (O,r_{\alpha}, \widehat \pi,\widehat m_R )]| + \varepsilon_n}} \sqrt{|P[f (O,\widehat r_{\alpha}, \widehat \pi,\widehat m_R ) - f (O,r_{\alpha}, \widehat \pi,\widehat m_R )]| + \varepsilon_n} \right] \\
        \stackrel{(i)}{\geq} & - \sqrt{\mathbb E |P[f (O,\widehat r_{\alpha}, \widehat \pi,\widehat m_R ) - f (O,r_{\alpha}, \widehat \pi,\widehat m_R )]| + \varepsilon_n} \sqrt{\mathbb E \sup_{\theta}\left[\frac{|(P - P_n)[f (O,\theta, \widehat \pi,\widehat m_R ) - f (O,r_{\alpha}, \widehat \pi,\widehat m_R )]|^2 }{|P[f (O,\theta, \widehat \pi,\widehat m_R ) - f (O,r_{\alpha}, \widehat \pi,\widehat m_R )]| + \varepsilon_n}  \right]} \\
        =& - \sqrt{\mathbb E |P[g (O,\widehat r_{\alpha}, \widehat \pi,\widehat m_R ) ]| + \varepsilon_n} \sqrt{\mathbb E \sup_{\theta}\left[\frac{|(P - P_n)[g (O,\theta, \widehat \pi,\widehat m_R ) ]|^2 }{|P[g (O,\theta, \widehat \pi,\widehat m_R )]| + \varepsilon_n}  \right]},
    \end{split}
\end{equation*}
where $g (O,\theta, \widehat \pi,\widehat m_R ) = f (O,\theta, \widehat \pi,\widehat m_R ) - f (O,r_{\alpha}, \widehat \pi,\widehat m_R )$ and $\varepsilon_n = 1/\sqrt{n}$. The step-$(i)$ follows from the Cauchy-schwartz inequality. By Theorem 3 of \cite{yang2024doubly} we have the following,
\begin{equation*}
    \begin{split}
        \mathbb E[\sup_{\theta} |(P_n - P) f(O, \theta, \widehat \pi, \widehat m_R) |] \leq \frac{\mathcal C}{\sqrt{n}},
    \end{split}
\end{equation*}
where $\mathcal C $ is a constant. With slight abuse of notation we will use $\mathcal C$ to denote constants in this proof. Under the assumption that $R$ is continuous and $\mathbb P(J(\widehat m_R(\theta, X) > 0) = 0$ for all $\theta \in \mathbb R$ we have the following using the Jump Lemma of \cite{angelopoulos2022conformal}, 
\[
\sup_{\theta} J(P_nf(O, \theta, \widehat \pi, \widehat m_R) , \theta) \stackrel{a.s.}{\leq} \frac{\overline \pi}{n}. 
\]
The above expression combined with the definition of $\widehat r_{\alpha}$ implies the following, 
\[
0 \leq P_nf(O, \widehat r_{\alpha}, \widehat \pi, \widehat m_R)  \stackrel{a.s.}{\leq} \frac{\overline \pi}{n}. 
\]
This combined with \eqref{eq:yang_dr} implies the following upper bound, 
\begin{equation*}
    \begin{split}
       & \mathbb E |P[g (O,\widehat r_{\alpha}, \widehat \pi,\widehat m_R ) ]| \\
       \leq & \E |P[f (O,\widehat r_{\alpha}, \widehat \pi,\widehat m_R ) ]| + \E |P[f (O, r_{\alpha}, \widehat \pi,\widehat m_R ) ]| \\
       \leq & \E |P[f (O,\widehat r_{\alpha}, \widehat \pi,\widehat m_R ) ]|  + \| \widehat \pi - \pi^*\|_2 \sup_{\theta}\| \widehat m_R(\theta, \cdot) - m_R^*(\theta, \cdot)\|_2 \\
       \leq & \E |P_n[f (O,\widehat r_{\alpha}, \widehat \pi,\widehat m_R ) ]| + \mathbb E[\sup_{\theta} |(P_n - P) f(O, \theta, \widehat \pi, \widehat m_R) |] + \| \widehat \pi - \pi^*\|_2 \sup_{\theta}\| \widehat m_R(\theta, \cdot) - m_R^*(\theta, \cdot)\|_2  \\
       \leq & \frac{\overline \pi}{n} + \frac{\mathcal C}{\sqrt{n}} + \| \widehat \pi - \pi^*\|_2 \sup_{\theta}\| \widehat m_R(\theta, \cdot) - m_R^*(\theta, \cdot)\|_2  .
    \end{split}
\end{equation*}
Therefore we have the following bound with slight abuse of notation, 
\begin{equation}
    \label{eq:bound_first_emp}
    \mathbb E |P[g (O,\widehat r_{\alpha}, \widehat \pi,\widehat m_R ) ]| \leq \frac{\mathcal C}{\sqrt{n}} + \| \widehat \pi - \pi^*\|_2 \sup_{\theta}\| \widehat m_R(\theta, \cdot) - m_R^*(\theta, \cdot)\|_2. 
\end{equation}
For ease of computation we define $\mu(\theta) = |Pg (O, \theta, \widehat \pi, \widehat m_R)|$. Let $\Delta m^*_R(\theta, X) = m^*_R(\theta, X) - m_R^*(r_{\alpha}, X)$ and $\Delta \widehat m_R(\theta, X) = \widehat m_R(\theta, X) - \widehat m_R(r_{\alpha}, X)$. We also use the following notations for the differences $\delta_{\pi}(X) = \widehat \pi(X) - \pi^*(X)$ and $\delta_m(\theta, X) = \Delta m_R^*(\theta, X) - \Delta \widehat m_R(\theta, X) $. The mean of $g (O, \theta, \widehat \pi, \widehat m_R)$ can be computed as follows, 
\begin{equation*}
    \begin{split}
       & Pg (O, \theta, \widehat \pi, \widehat m_R) \\
       = & \mathbb E [\mathbb P(C = \infty|X) (\widehat \pi(X) - \pi^*(X))(\Delta m_R^*(\theta, X) - \Delta \widehat m_R(\theta, X) ) + \mathbb P(C = 0|X)\Delta m_R^*(\theta, X)  ].
    \end{split}
\end{equation*}
Note that because of the monotonicity of $\theta \mapsto \Delta m_R^*(\theta, \cdot) $ we have the following using triangle inequality and cauchy schwartz inequality, 
\[
\mathbb E [\mathbb P(C = 0|X) |\Delta m_R^*(\theta, X)|] \leq | Pg (O, \theta, \widehat \pi, \widehat m_R)| + 2\| \widehat \pi - \pi^*\|_2 \sup_{\theta}\| \widehat m_R(\theta, \cdot) - m_R^*(\theta, \cdot)\|_2.
\]
Using the boundedness of $\pi^*(X), \widehat \pi(X) \leq \overline \pi $, we have $\mathbb P(C = \infty |X) \widehat \pi(X)^2 \leq \mathcal C\mathbb P(C = 0|X)$. Therefore we have the following inequality,
\begin{equation*}
    \begin{split}
         & P g^2(O, \theta, \widehat \pi, \widehat m_R) \\
      = & \mathbb E [\textbf{1}\{C = \infty\} \widehat \pi(X)^2 [\{\textbf{1}\{Y \leq \theta \} - \textbf{1}\{Y \leq r_{\alpha}\} - \{\widehat m_R(\theta, X) - \widehat m_R(r_{\alpha}, X) \}]^2] \\
      & + \mathbb E[\textbf{1}\{C = 0\} \{\widehat m_R(\theta, X) - \widehat m_R(r_{\alpha}, X) \}^2] \\
      = & \mathbb E [\mathbb P\{C = \infty|X\} \widehat \pi(X)^2\mathbb E[ [\{\textbf{1}\{Y \leq \theta \} - \textbf{1}\{Y \leq r_{\alpha}\} - \Delta \widehat m_R(\theta, X)]^2|X]] \\
      & + \mathbb E[\mathbb{P}\{C = 0|X\} \Delta \widehat m_R(\theta, X)^2]  \\
      \leq & \mathbb E [\mathbb P\{C = \infty|X\} \widehat \pi(X)^2\mathbb E[ [\{\textbf{1}\{Y \leq \theta \} - \textbf{1}\{Y \leq r_{\alpha}\} - \Delta \widehat m_R(\theta, X)]^2|X]] \\
      & + \mathbb E[\mathbb{P}\{C = 0|X\} \{2|\Delta m_R^*(\theta, X)| + 2 \delta_m(\theta, X)^2\}] \\
      \leq & \mathcal C \mathbb E[\mathbb{P}\{C = 0|X\}|\Delta m_R^*(\theta, X)| ]  + \mathcal C \mathbb E[\delta_m(\theta, X)^2].
    \end{split}
\end{equation*}
Combining this with the prior bounds we have,
\begin{equation*}
   P g^2(O, \theta, \widehat \pi, \widehat m_R) \leq K|Pg (O, \theta, \widehat \pi, \widehat m_R) | + K_1 \{ \| \widehat \pi - \pi^*\|_2 \sup_{\theta}\| \widehat m_R(\theta, \cdot) - m_R^*(\theta, \cdot)\|_2 + \sup_{\theta}\| \widehat m_R(\theta, \cdot) - m_R^*(\theta, \cdot)\|_2^2 \},
\end{equation*}
where $K, K_1$ are suitable constants. For ease of explanation we will denote the term $ K_1 \{ \| \widehat \pi - \pi^*\|_2 \sup_{\theta}\| \widehat m_R(\theta, \cdot) - m_R^*(\theta, \cdot)\|_2 + \sup_{\theta}\| \widehat m_R(\theta, \cdot) - m_R^*(\theta, \cdot)\|_2^2 \}$ by $v$ in the remaining portion of the proof. Consider the family of functions $\mathscr F = \{g(O, \theta, \widehat \pi, \widehat m_R)|\theta \in \mathbb R \}$ whose envelope is $\overline \pi$. Let $H(\delta) = \sup_{Q} \log \mathcal N(\delta \overline \pi, \mathscr F, L_2(Q))$ and $J(\delta) = \int_{0}^{\delta} \sqrt{1 + H(u)}du $. We do the peeling using $\mu(\theta)$. Let us define, 
\begin{equation*}
    \begin{split}
        q_0 = &\sup_{\theta} \mu(\theta) \leq \overline \pi,\quad  q_j = 2^{-j}q_0 \quad \mbox{for} \quad j \geq 1 , \\
        A_j = &\{\theta \in \mathbb R| q_{j+1} < \mu(\theta) \leq q_j\} \quad \mbox{for} \quad j \geq 0. 
    \end{split}
\end{equation*}
From our analysis we know that the following holds on the shell $A_j$ ($j \geq 0$), 
\[
Pg^2(O, \theta, \widehat \pi, \widehat m_R) \leq K q_j + v \coloneqq q_j^+ . 
\]
We define $\delta_j = \sqrt{q_j^+}/\overline \pi$ for $j \geq 0$. The following two results hold for any local family of functions $\mathscr F_{loc}$ whose envelope is $\overline \pi$, 
\begin{equation}
    \label{eq:useful_results}
    \begin{split}
        \mathbb E \left[ \sup_{f \in \mathscr F_{loc}} |(P_n - P)f \right] \leq & \mathcal C \left[\frac{\overline \pi}{\sqrt{n}} J(\delta) + \frac{J(\delta)^2}{\delta^2 n} \right], \\
        \mathbb E \left[ \sup_{f \in \mathscr F_{loc}} |(P_n - P)f \right]^2
        \leq & 2 \left(  \mathbb E \left[ \sup_{f \in \mathscr F_{loc}} |(P_n - P)f \right]\right)^2 + \frac{2\overline q}{n} + \frac{8\overline \pi^2}{9 n^2}, 
    \end{split}
\end{equation}
where $\overline q = \sup_{f \in \mathscr F_{loc}} Pf^2$. The first result is Theorem $2.14.2$ from \cite{van1996weak} and the second result is obtained from Theorem $12.5$ (Bousquet's inequality) \cite{boucheron2003concentration} by integrating the tail bounds. Note that if $\mathscr F_{loc} \in A_j$ then $\overline q = \sup_{\theta \in A_j} Pg^2(O, \theta, \widehat \pi, \widehat m_R) \leq K q_j + v = q_j^+ $. Using \eqref{eq:useful_results} we have the following on the shell $A_j$ ($j \geq 0$), 
\begin{equation*}
    \begin{split}
        \mathbb E \sup_{\theta \in A_j} |(P_n - P)g(O, \theta, \widehat \pi, \widehat m_R)|^2 \leq \frac{4\mathcal C^2 \overline \pi^2}{n}J(\delta_j)^2 + \frac{2q_j^+}{n} + \frac{8\overline \pi^2}{9 n^2} + \frac{4\mathcal C^2J(\delta_j)^4}{\delta_j^4 n^2}. 
    \end{split}
\end{equation*}
Let $\overline H(\delta) = (1/\delta)\int_{0}^{\delta} H(u)du$. Using Jensen's inequality we have $J(\delta) \leq \delta \sqrt{1 + \overline H(\delta)}$. Therefore we have the following upper bounds for any $j \geq 0$,
\begin{equation*}
    \begin{split}
        J(\delta_j)^2 \leq & \delta_j^2 (1 + \overline H(\delta_j)) = (2Kq_{j+1} + v)(1 + \overline H(\delta_j))/\overline \pi^2, \\
        q_j^+ =& Kq_j + v = 2Kq_{j+1} + v, \\
        J(\delta_j)^4/\delta_j^4 \leq & (1 + \overline H(\delta_j))^2. 
    \end{split}
\end{equation*}
Let $J^* = \min\{j \in \mathbb N: q_{j+1} \leq \varepsilon_n\}$ and $\delta_{\varepsilon_n} = \sqrt{K\varepsilon_n +v}/\overline \pi$. By peeling we can bound the ratio type empirical process term, 
\begin{equation*}
    \begin{split}
        & \mathbb E \sup_{\theta}\left[\frac{|(P - P_n)[g (O,\theta, \widehat \pi,\widehat m_R ) ]|^2 }{|P[g (O,\theta, \widehat \pi,\widehat m_R )]| + \varepsilon_n}  \right] \\
        \leq & \sum_{j \leq J^*} \frac{\mathbb E \sup_{\theta \in A_j} |(P_n - P)g(O, \theta, \widehat \pi, \widehat m_R)|^2}{q_{j+1} } + \frac{\mathbb E \sup_{\mu(\theta) \leq \varepsilon_n} |(P_n - P)g(O, \theta, \widehat \pi, \widehat m_R)|^2}{\varepsilon_n} \\
        \coloneqq& \mathbf I + \mathbf{II}.
    \end{split}
\end{equation*}
Using the bounds that we have, we obtain the following,
\begin{equation*}
    \begin{split}
        \mathbf I = & \sum_{j \leq J^*} \frac{\mathbb E \sup_{\theta \in A_j} |(P_n - P)g(O, \theta, \widehat \pi, \widehat m_R)|^2}{q_{j+1} } \\
        \leq & \frac{1}{n}\sum_{j \leq J^*} \left[(1 + \overline H(\delta_j)) \left\{8\mathcal C^2K +  \frac{4 \mathcal C^2v}{\varepsilon_n} \right\} + 4K + \frac{2v}{\varepsilon_n} +  \frac{8\overline \pi^2}{9 n} + \frac{4\mathcal C^2}{n \varepsilon_n} (1 + \overline H(\delta_j))^2  \right] \\
        \leq & \frac{1}{n} \left[\int_{\delta_{\varepsilon_n}}^{\delta_0}(1 + \overline H(\delta)) d\log(1/\delta) \left\{8\mathcal C^2K +  \frac{4\mathcal C^2v}{\varepsilon_n} \right\} + 4KJ^* + \frac{2vJ^*}{\varepsilon_n} +  \frac{8\overline \pi^2 J^*}{9 n} + \frac{4\mathcal C^2}{n \varepsilon_n} \int_{\delta_{\varepsilon_n}}^{\delta_0}(1 + \overline H(\delta))^2d\log(1/\delta)   \right].
    \end{split}
\end{equation*}
In the above bound $J^* = O(\log(n))$. Moreover if $\theta \mapsto \widehat m(\theta, X)$ is a monotonic function, it is staightforward to show that $H(\delta) \leq \mathcal C + a \log(b/\delta)$ where $\mathcal C, a, b$ are constants. Therefore $\overline H(\delta) = (1/\delta)\int_{0}^{\delta} H(u)du \leq \mathcal C + a(1 + \log(b/\delta))$. This in turn implies that, 
\[
\int_{\delta_{\varepsilon_n}}^{\delta_0}(1 + \overline H(\delta))d\log(1/\delta) \leq \int_{\delta_{\varepsilon_n}}^{\delta_0} (1 + \mathcal C + a(1 + \log(b/\delta)))d \log(1/\delta) \leq \mathcal C(\log n)^2. 
\]
We can similarly bound the second term, 
\begin{equation*}
    \begin{split}
        \mathbf {II} = & \frac{\mathbb E \sup_{\mu(\theta) \leq \varepsilon_n} |(P_n - P)g(O, \theta, \widehat \pi, \widehat m_R)|^2}{\varepsilon_n} \\
        \leq & \frac{1}{n\varepsilon_n} \left[4\mathcal C^2(2K \varepsilon_n + v)(1 + \overline H(\delta_{\varepsilon_n})) + 2(K\varepsilon_n +v) + \frac{8 \overline \pi^2}{9 n} + \frac{4\mathcal C^2}{n}(1 + \overline H(\delta_{\varepsilon_n}))^2\right].
    \end{split}
\end{equation*}
Combining these bounds we have the following, 
\begin{equation*}
    \begin{split}
        \mathbb E \sup_{\theta}\left[\frac{|(P - P_n)[g (O,\theta, \widehat \pi,\widehat m_R ) ]|^2 }{|P[g (O,\theta, \widehat \pi,\widehat m_R )]| + \varepsilon_n}  \right]  \leq \mathfrak C_n'\left(\frac{1}{n} + \frac{v}{n \varepsilon_n} \right) = \mathfrak C_n'\left(\frac{1}{n} + \frac{v}{\sqrt{n}} \right),
    \end{split}
\end{equation*}
where $\mathfrak C_n' = \mathcal C(\log n)^2$ for some constant $\mathcal C$. Let $a_\pi = \|\widehat \pi - \pi^*\|_2$ and $a_m = \sup_{\theta}\| \widehat m_R(\theta, \cdot) - m_R^*(\theta, \cdot)\|_2$ respectively denote the errors in estimating the nuisance functions. Combining all our previous results we have the following, 
\begin{equation*}
    \begin{split}
        & \mathbb E[Pf (O,\widehat r_{\alpha}, \widehat \pi,\widehat m_R )]  \\
        \geq & - \sqrt{\mathbb E |P[g (O,\widehat r_{\alpha}, \widehat \pi,\widehat m_R ) ]| + \varepsilon_n} \sqrt{\mathbb E \sup_{\theta}\left[\frac{|(P - P_n)[g (O,\theta, \widehat \pi,\widehat m_R ) ]|^2 }{|P[g (O,\theta, \widehat \pi,\widehat m_R )]| + \varepsilon_n}  \right]} \\
        \geq & -\left\{\frac{\mathcal C}{n^{1/4}} + \sqrt{a_\pi a_m} \right\}\left\{\frac{\mathfrak C_n'^{1/2}}{n^{1/2}} + \frac{\mathfrak C_n'^{1/2}\{\sqrt{a_\pi a_m} + a_m\} }{n^{1/4}} \right\} \\
        \geq & - \mathfrak C_n\left\{\frac{1}{n^{3/4}} + \frac{\sqrt{a_\pi a_m} + a_m}{n^{1/2}} + \frac{\sqrt{a_\pi a_m} \{ \sqrt{a_\pi a_m} + a_m \}}{n^{1/4}} \right\} \\
        = & - \mathfrak C_n\left\{\frac{1}{n^{3/4}} +  \frac{\{ \sqrt{a_\pi a_m} + a_m \}}{n^{1/4}} \left[\frac{1}{n^{1/4}}+ \sqrt{a_\pi a_m} \right]\right\},
    \end{split}
\end{equation*}
where $\mathfrak C_n = \mathcal C(\log n)$ for some suitable constant $C$. Therefore we get the following lower bound on the coverage using, 
\begin{equation*}
    \begin{split}
        &\mathbb P(R_{n+1} \leq \widehat r_{\alpha} |C_{n+1} = 0) \\
        \geq & 1 - \alpha - \frac{\| \widehat \pi - \pi^*\|_2 \sup_{\theta}\| \widehat m_R(\theta, \cdot) - m_R^*(\theta, \cdot)\|_2}{\mathbb P(C = 0)} + \frac{\mathbb E[Pf (O,\widehat r_{\alpha}, \widehat \pi,\widehat m_R )]}{\mathbb P(C = 0)}  \\
        \geq & 1 - \alpha - \frac{a_\pi a_m}{\mathbb P(C = 0)} \\
        & - \frac{\mathfrak C_n}{\mathbb P(C = 0)}\left\{\frac{1}{n^{3/4}} +  \frac{\{ \sqrt{a_\pi a_m} + a_m \}}{n^{1/4}} \left[\frac{1}{n^{1/4}}+ \sqrt{a_\pi a_m} \right]\right\}. 
    \end{split}
\end{equation*}
Following the same proof technique we also have the following upper bound on coverage,
\begin{equation*}
    \begin{split}
        &\mathbb P(R_{n+1} \leq \widehat r_{\alpha} |C_{n+1} = 0) \\
        \leq & 1 - \alpha + \frac{\| \widehat \pi - \pi^*\|_2 \sup_{\theta}\| \widehat m_R(\theta, \cdot) - m_R^*(\theta, \cdot)\|_2}{\mathbb P(C = 0)} + \frac{\mathbb E[Pf (O,\widehat r_{\alpha}, \widehat \pi,\widehat m_R )]}{\mathbb P(C = 0)}  \\
        \leq & 1 - \alpha + \frac{a_\pi a_m}{\mathbb P(C = 0)} \\
        & + \frac{\mathfrak C_n}{\mathbb P(C = 0)}\left\{\frac{1}{n^{3/4}} +  \frac{\{ \sqrt{a_\pi a_m} + a_m \}}{n^{1/4}} \left[\frac{1}{n^{1/4}}+ \sqrt{a_\pi a_m} \right]\right\}. 
    \end{split}
\end{equation*}
This completes the proof of the theorem. 

\section{Proof of \Cref{thm:validity_proof}}
\label{app:proof_valid_cmdp}
The final estimate $\widehat r_{\alpha}$ of $r_{\alpha}$ is defined using the calibration set $\mathcal{D}^{cal}$,  
  \begin{equation*}
      \begin{split}
      \widehat  r_{\alpha} = &\inf\Bigg\{\theta \in \mathbb{R}: \frac{\sum_{i = 1}^n f(O_i, \theta, \widehat r_{\alpha, 0}, \widehat \pi, \widehat m_R(\theta))}{n+1} -  \frac{\overline \pi}{n+1}  \geq 0\Bigg\},  
      \end{split}
  \end{equation*}
where,
  \[
  f(O, \theta, \pi, m_R(\theta')) = \textbf{1}\{C = \infty\}\pi(X) [\textbf{1}\{R(X, Y) \leq \theta\} - m_R(\theta', X)] + \textbf{1}\{C = 0\} [m_R(\theta', X) - (1 - \alpha)].
  \]
We note that if $\widehat \pi$ and $\widehat m_R$ are estimates of the true functions $\pi^*$ and $m_R^*$ respectively, then we can bound the difference between $\mathbb{E}[f(O,\theta, \pi^*, m_R^*(\theta))]$ and $\mathbb{E}[f(O,\theta, \widehat \pi, \widehat m_R(\theta'))]$ for any $\theta, \theta' \in \mathbb R$.
\begin{equation*}
    \begin{split}
        &\mathbb{E}[f(O,\theta, \widehat \pi, \widehat m_R(\theta'))] \\
     = &     \mathbb E[\mathbb{P}(C = \infty |X) \widehat \pi(X) [\mathbb P(R \leq \theta| X) - \widehat m_R(\theta', X)] + \mathbb{P}(C = 0 |X) [\widehat m_R(\theta', X) - (1 - \alpha)]] \\
     = &  \mathbb E[\mathbb{P}(C = \infty |X)\{ \widehat \pi(X) - \pi^*(X)\} [m_R^*(\theta, X) - \widehat m_R(\theta', X)]] \\
     & + \mathbb E [ \mathbb{P}(C = 0 |X) [m_R^*(\theta, X) - \widehat m_R(\theta', X)] ] \\
     & + \mathbb E[\mathbb{P}(C = 0 |X) [\widehat m_R(\theta', X) - (1 - \alpha)]] \\
     = & \mathbb E[\mathbb{P}(C = \infty |X)\{ \widehat \pi(X) - \pi^*(X)\} [m_R^*(\theta, X) - \widehat m_R(\theta', X)]] \\
     & + \mathbb E [ \mathbb{P}(C = 0 |X) [m_R^*(\theta, X) - (1 - \alpha)] .
    \end{split}
\end{equation*}
Therefore we have the following upper bound for any $\theta, \theta' \in \mathbb R$,
\begin{equation}
    \label{eq:dr_one_step}
    \begin{split}
       & | \mathbb{E}[f(O,\theta, \widehat \pi, \widehat m_R(\theta'))] -  \mathbb{E}[f(O,\theta, \pi^*, m_R^*(\theta))]| \\
      = & | \mathbb E[\mathbb{P}(C = \infty |X)\{ \widehat \pi(X) - \pi^*(X)\} [m_R^*(\theta, X) - \widehat m_R(\theta', X)]] |    \\
     \leq & \| \widehat \pi(X) - \pi^*(X) \|_2 \| m_R^*(\theta, X) - \widehat m_R(\theta', X) \|_2 \\
     \leq &  \| \widehat \pi(X) - \pi^*(X) \|_2 \| m_R^*(\theta, X) - \widehat m_R(\theta, X) \|_2 + \| \widehat \pi(X) - \pi^*(X) \|_2 \| m_R^*(\theta, X) -  m_R^*(\theta', X) \|_2 .
    \end{split}
\end{equation}
Note that if we set $\theta = \theta'$ then we get back the result of Theorem 2 of \cite{yang2024doubly},
\begin{equation}
\label{eq:double_robust}
    \sup_{\theta\in \mathbb{R}} \left|\mathbb{E}[f(O,\theta, \pi^*, m_R^*(\theta))]- \mathbb{E}[f(O,\theta, \widehat \pi, \widehat m_R(\theta))] \right| \leq \| \widehat \pi  - \pi^*\|_2 \sup_{\theta \in \mathbb{R}} \| \widehat m_R (\theta, \cdot) - m_R^*(\theta, \cdot)\|_2.
\end{equation}
Let us define the following, 
\[
\widehat r_{\alpha}' = \inf \left\{\theta: \frac{1}{n+1} \sum_{i = 1}^{n+1} f(O_i, \theta, \widehat \pi, \widehat m_R(\widehat r_{\alpha,0})) \geq 0 \right\}. 
\]
We note that each of the functions $\theta \mapsto f(O, \theta, \widehat \pi, \widehat m_R(\widehat r_{\alpha,0}))$ is a right-continuous function in $\theta$ and the same property also holds for the mean function $\theta \mapsto (1/(n+1)) \sum_{i = 1}^{n+1} f(O_i, \theta, \widehat \pi, \widehat m_R(\widehat r_{\alpha,0}))$. Thus by taking limits we get the following, 
\[
\frac{1}{n+1} \sum_{i = 1}^{n+1}f(O_i, \widehat r_{\alpha}', \widehat \pi, \widehat m_R(\widehat r_{\alpha,0})) \geq 0
\]
We observe that for the set $E = \{O_1, \cdots, O_{n+1}\}$ we have the following because of the exchangeability of the random variables $O_i$'s, 
\begin{equation*}
    \begin{split}
        & O_i | E \sim \frac{1}{n+1} \sum_{i = 1}^{n+1} \delta_{O_i} \\
        \implies & \mathbb{E}[f(O_{n+1}, \theta, \widehat \pi, \widehat m_R(\widehat r_{\alpha,0})) | E] = \frac{1}{n+1} \sum_{i = 1}^{n +1} f(O_i, \theta, \widehat \pi, \widehat m_R(\widehat r_{\alpha,0})). 
 \end{split}
\end{equation*}
 Combining with the previous deductions we obtain that, 
 \begin{equation*}
     \begin{split}
         & \mathbb{E}[f(O_{n+1}, \widehat r_{\alpha}', \widehat \pi, \widehat m_R(\widehat r_{\alpha,0}))|E] = \frac{1}{n+1}\sum_{i = 1}^{n +1} f(O_{i}, \widehat r_{\alpha}', \widehat \pi, \widehat m_R(\widehat r_{\alpha,0}))\geq 0 \\
         \implies & \mathbb{E}[f(O_{n+1}, \widehat r_{\alpha}', \widehat \pi, \widehat m_R(\widehat r_{\alpha,0}))] \geq 0.
     \end{split}
 \end{equation*}
The last line follows because of tower property. We now observe that $\widehat r_{\alpha}$ has the following property, 
\begin{equation*}
    \begin{split}
        &\frac{1}{n+1} \sum_{i = 1}^{n+1} f(O_{i}, \widehat r_{\alpha}, \widehat \pi, \widehat m_R(\widehat r_{\alpha,0})) \\
        =& \frac{1}{n+1} \sum_{i = 1}^{n} f(O_{i}, \widehat r_{\alpha}, \widehat \pi, \widehat m_R(\widehat r_{\alpha,0})) + \frac{f(O_{n+1}, \widehat r_{\alpha}, \widehat \pi, \widehat m_R(\widehat r_{\alpha,0}))}{n+1} \\
        \geq &  \frac{1}{n+1} \sum_{i = 1}^{n} f(O_{i}, \widehat r_{\alpha}, \widehat \pi, \widehat m_R(\widehat r_{\alpha,0})) - \frac{\overline \pi}{n+1} \\
        \geq & 0. 
    \end{split}
\end{equation*}
The last line follows from the definition of $\widehat r_{\alpha}$. Thus we see that $\widehat r_{\alpha}$ belongs to the set whose infimum is $\widehat r_{\alpha}'$ and hence $\widehat r_{\alpha}' \leq \widehat r_{\alpha}$. Using the monotonicity of $\theta \mapsto f(O_{n+1}, \theta, \widehat \pi, \widehat m_R(\widehat r_{\alpha,0}))$ we have the following coverage guarantee, 
\begin{equation}
\label{eq:cmdp_lower_first}
    \begin{split}
       & \mathbb{P}( R(X_{n+1},Y_{n+1}) \leq \widehat r_{\alpha}| C_{n+1} = 0) \\
       =& 1 - \alpha +   \frac{\mathbb{E}[f(O_{n+1}, \widehat r_{\alpha},  \pi^*,  m_R^*(\widehat r_{\alpha}))]}{\mathbb{P}(C = 0)} \quad \mbox{(using \eqref{coverage_link_cmdp})}\\
       \stackrel{\eqref{eq:dr_one_step}}{\geq} & 1 -\alpha \\
       & + \frac{\mathbb{E}[f(O_{n+1}, \widehat r_{\alpha},  \widehat \pi, \widehat m_R(\widehat r_{\alpha, 0}))] - \|\widehat \pi  - \pi^*\|_2 \sup_{\theta} \| \widehat m_R (\theta, \cdot) - m_R^*(\theta, \cdot)\|_2 -  \| \widehat \pi - \pi^*\|_2 \|m_R^*(\widehat r_{\alpha}, \cdot) - m_R^*(\widehat r_{\alpha, 0}, \cdot) \|_2}{\mathbb{P}(C = 0)}    \\  
       \stackrel{(i)}{\geq} &  1 -\alpha \\
       &+ \frac{\mathbb{E}[f(O_{n+1}, \widehat r_{\alpha}',  \widehat \pi, \widehat m_R(\widehat r_{\alpha, 0}))] - \|\widehat \pi  - \pi^*\|_2 \sup_{\theta} \| \widehat m_R (\theta, \cdot) - m_R^*(\theta, \cdot)\|_2 -  \| \widehat \pi - \pi^*\|_2 \|m_R^*(\widehat r_{\alpha}, \cdot) - m_R^*(\widehat r_{\alpha, 0}, \cdot) \|_2}{\mathbb{P}(C = 0)}    \\
       \stackrel{(ii)}{\geq} & 1 - \alpha - \frac{  \|\widehat \pi  - \pi^*\|_2 \sup_{\theta} \| \widehat m_R (\theta, \cdot) - m_R^*(\theta, \cdot)\|_2}{\mathbb{P}(C = 0)} - \frac{\| \widehat \pi - \pi^*\|_2 \|m_R^*(\widehat r_{\alpha}, \cdot) - m_R^*(\widehat r_{\alpha, 0}, \cdot) \|_2 }{\mathbb{P}(C = 0)}.
    \end{split}
\end{equation}

In step-$(i)$ we use the monotonicity of the function $\theta \mapsto f(O_{n+1}, \theta, \widehat \pi, \widehat m_R(\widehat r_{\alpha,0}))$ and $\widehat r_{\alpha} \geq \widehat r_{\alpha}'$. In step-$(ii)$ we have used $\mathbb{E}[f(O_{n+1}, \widehat r_{\alpha}',  \widehat \pi, \widehat m_R(\widehat r_{\alpha, 0}))] \geq 0$. We will obtain an upper bound on $\|m_R^*(\widehat r_{\alpha}, X) - m_R^*(\widehat r_{\alpha, 0}, X) \|_2$ later in the proof. Let us denote the mean function $(1/n)\sum_{i = 1}^n f(O_i,\theta, \pi, m_R(\theta'))$ by $M_n(\theta; \pi, m_R(\theta'))$. For any given function $f$ we define the jump function $J(f, x) = \lim_{\epsilon \downarrow 0} | f(x - \epsilon) - f(x)|$. For the function $f(O_i,\theta, \widehat \pi, \widehat m_R(\widehat r_{\alpha, 0}))$ the only jump is at $\theta = R(X_i, Y_i)$ and the magnitude of the jump is $|\widehat \pi(X_i)| \leq \overline \pi $. Thus $J(G_i, \theta) \leq \overline \pi$ for $i = 1 , \cdots, n$ where $G_i = f(O_i,\theta, \widehat \pi, \widehat m_R(\widehat r_{\alpha, 0}))$. Moreover because of the continuity of $R(X_i,Y_i)$ we have $\mathbb{P}(J(G_i, \theta) > 0) = \mathbb{P}(R(X_i, Y_i) = \theta) = 0$ for all $\theta \in \mathbb{R}$ and for $i = 1 , \cdots, n$. Using the jump lemma in \cite{angelopoulos2022conformal} we get that, 
\begin{equation}
    \label{eq:jump_eq}
    \sup_{\theta}J(M_n(\cdot;\widehat \pi, \widehat m_R(\widehat r_{\alpha, 0})) , \theta) \stackrel{a.s.}{\leq} \frac{\overline{ \pi}}{n}. 
\end{equation}
We now define the following, 
\[
\widehat r_{\alpha}'' = \inf \left\{\theta: M_{n+1}(\theta;\widehat \pi, \widehat m_R(\widehat r_{\alpha, 0})) - \frac{2\overline \pi}{n+1} \geq 0\right\}.
\]
We observe that we can write $M_{n+1}(\widehat r_{\alpha}'';\widehat \pi, \widehat m_R)$ in the following manner, 
\[
M_{n+1}(\widehat r_{\alpha}''; \widehat \pi, \widehat m_R(\widehat r_{\alpha, 0})) = \frac{2 \overline \pi}{n+1} + \left( M_{n+1}(\widehat r_{\alpha}''; \widehat \pi, \widehat m_R(\widehat r_{\alpha, 0})) - \frac{ 2 \overline \pi}{n+1}\right).
\]
We see that the first term is $2\overline \pi/(n+1)$ and the second term can be at most $\overline \pi/(n+1)$ because of \eqref{eq:jump_eq}. Thus we see that $M_{n+1}(\widehat r_{\alpha}''; \widehat \pi, \widehat m_R(\widehat r_{\alpha, 0}))$ is bounded above by $3 \overline \pi/(n+1)$. From our previous discussions we observe the following, 
\begin{equation*}
    \begin{split}
        \mathbb{E}[f(O_{n+1}, \widehat r_{\alpha}'',  \widehat \pi, \widehat m_R(\widehat r_{\alpha, 0}))] = \mathbb E[M_{n+1}(\widehat r_{\alpha}'';  \widehat \pi, \widehat m_R(\widehat r_{\alpha, 0})) |E]  \leq \frac{3 \overline \pi}{n+1} .
    \end{split}
\end{equation*} 
We now observe the following regarding $\widehat r_{\alpha}''$, 
\begin{equation*}
    \begin{split}
       &\frac{\sum_{i = 1}^n f(O_{i}, \widehat r_{\alpha}'',  \widehat \pi, \widehat m_R(\widehat r_{\alpha, 0}))}{n+1} - \frac{\overline \pi}{n+1} \\ 
       = &  \frac{\sum_{i = 1}^{n+1} f(O_{i}, \widehat r_{\alpha}'',  \widehat \pi, \widehat m_R(\widehat r_{\alpha, 0}))}{n+1} - \frac{f(O_{n+1}, \widehat r_{\alpha}'',  \widehat \pi, \widehat m_R(\widehat r_{\alpha, 0}))}{n+1} - \frac{\overline \pi}{n+1}         \\
        \geq & M_{n+1}(\widehat r_{\alpha}''; \widehat \pi, \widehat m_R (\widehat r_{\alpha, 0})) - \frac{2 \overline \pi}{n+1} \\
         \geq & 0. 
    \end{split}
\end{equation*}

In the last line we used the right continuity of the mean function $\theta \mapsto M_{n+1}(\theta; \widehat \pi, \widehat m_R(\widehat r_{\alpha, 0}))$. This shows that $\widehat r_{\alpha}''$ belongs to the set whose infimum is $\widehat r_{\alpha}$. Therefore we get that $\widehat r_{\alpha} \leq \widehat r_{\alpha}''$. By following the same steps as in the first part of the proof we can get, 
\begin{equation}
\label{eq:cmdp_upper_first}
    \begin{split}
       & \mathbb{P}( R(X_{n+1},Y_{n+1}) \leq \widehat r_{\alpha}| C_{n+1} = 0) \\
       =& 1 - \alpha +   \frac{\mathbb{E}[f(O_{n+1}, \widehat r_{\alpha},  \pi^*,  m_R^*(\widehat r_{\alpha}))]}{\mathbb{P}(C = 0)} \quad \mbox{(using \eqref{coverage_link_cmdp})}\\
       \stackrel{\eqref{eq:dr_one_step}}{\leq} & 1 -\alpha \\
       &+ \frac{\mathbb{E}[f(O_{n+1}, \widehat r_{\alpha},  \widehat \pi, \widehat m_R(\widehat r_{\alpha, 0}))] + \|\widehat \pi  - \pi^*\|_2 \sup_{\theta} \| \widehat m_R (\theta, \cdot) - m_R^*(\theta, \cdot)\|_2 +  \| \widehat \pi - \pi^*\|_2 \|m_R^*(\widehat r_{\alpha}, \cdot) - m_R^*(\widehat r_{\alpha, 0}, \cdot) \|_2}{\mathbb{P}(C = 0)}    \\  
       \stackrel{(i)}{\leq} &  1 -\alpha \\
       &+ \frac{\mathbb{E}[f(O_{n+1}, \widehat r_{\alpha}'',  \widehat \pi, \widehat m_R(\widehat r_{\alpha, 0}))] + \|\widehat \pi  - \pi^*\|_2 \sup_{\theta} \| \widehat m_R (\theta, \cdot) - m_R^*(\theta, \cdot)\|_2 +  \| \widehat \pi - \pi^*\|_2 \|m_R^*(\widehat r_{\alpha}, \cdot) - m_R^*(\widehat r_{\alpha, 0}, \cdot) \|_2}{\mathbb{P}(C = 0)}    \\
       \stackrel{(ii)}{\leq} & 1 - \alpha + \frac{3 \overline \pi}{\mathbb P(C = 0)(n+1)} + \frac{  \|\widehat \pi  - \pi^*\|_2 \sup_{\theta} \| \widehat m_R (\theta, \cdot) - m_R^*(\theta, \cdot)\|_2}{\mathbb{P}(C = 0)} + \frac{\| \widehat \pi - \pi^*\|_2 \|m_R^*(\widehat r_{\alpha}, \cdot) - m_R^*(\widehat r_{\alpha, 0}, \cdot) \|_2}{\mathbb{P}(C = 0)}.
    \end{split}
\end{equation}
In step-$(i)$ we use the monotonicity of the function $\theta \mapsto f(O_{n+1}, \theta, \widehat \pi, \widehat m_R(\widehat r_{\alpha,0}))$ and $\widehat r_{\alpha} \leq \widehat r_{\alpha}''$. In step-$(ii)$ we have used $\mathbb{E}[f(O_{n+1}, \widehat r_{\alpha}'',  \widehat \pi, \widehat m_R(\widehat r_{\alpha, 0}))] \leq 3 \overline \pi/(n+1)$. This completes the proof of the theorem. 
\section{Proof of Proposition~\ref{prop:rscp_init_error}}
\label{app:proof_bound_init_cmdp}
We bound $\|m_R^*(\widehat r_{\alpha}, X) - m_R^*(\widehat r_{\alpha, 0}, X) \|_2 $ in \Cref{lem:cmdp_init_final}. 
\begin{lemma}
\label{lem:cmdp_init_final}
Under assumptions that $(1)$ $\theta \mapsto \widehat m_R(\theta, X)$ is a right continuous function, $(2)$ $|P[\textbf 1\{C = \infty\} \widehat \pi(X) [\textbf{1}\{R\leq \theta\} - \textbf{1}\{R\leq \theta'\}] \}]| \geq \mathfrak C_1| \theta - \theta'|$ for $\theta, \theta' \in [r_{\alpha }- \eta, r_{\alpha}+ \eta]$, $(3)$ $\|m_R^*(\theta, X) - m_R^*(\theta', X)\|_2 \leq \mathfrak C_2 | \theta - \theta'| $ for $\theta, \theta' \in [r_{\alpha }- \eta, r_{\alpha}+ \eta]$ we have the following upper bound with probability at least $ 1- 2 \delta$ ($c'$ is a constant depending only on $\overline \pi$), 
\begin{equation*}
    \begin{split}
          & \|m_R^*(\widehat r_{\alpha} , X) - m_R^*(\widehat r_{\alpha, 0}, X) \|_2 \\
      \leq & \mathfrak C_2  \left[\mathcal{C} \frac{\log n + \log(1/\delta)}{n} + \frac{2 \overline \pi}{n} + \frac{3 - 2\alpha}{n_2 + 1} \right. \\
        & \left. + (1 + \|\widehat \pi  - \pi^*\|_2) \sup_{\theta} \| \widehat m_R (\theta, \cdot) - m_R^*(\theta, \cdot)\|_2 +  \frac{\sup_\theta \sigma_P[f(O, \theta, \widehat \pi, \widehat m_R(\theta))]}{\sqrt{n \delta}} \right]. 
    \end{split}
\end{equation*}
\end{lemma}

\begin{proof}[Proof of \Cref{lem:cmdp_init_final}]
 From here on we use $P_nf(O, \theta, \pi, m_R(\theta'))$ to denote the average $(1/n) \sum_{i = 1}^n f(O_i, \theta, \pi, m_R(\theta'))$ computed over all the calibration points. We also use $P f(O, \theta, \pi, m_R(\theta'))$ to denote the population mean $\mathbb E[f(O, \theta, \pi, m_R(\theta'))]$. We observe the following decomposition, 
 \begin{equation*}
     \begin{split}
       P_nf(O, \widehat r_{\alpha}, \widehat \pi, \widehat m_R(\widehat r_{\alpha, 0})) = &  P_nf(O, \widehat r_{\alpha, 0}, \widehat \pi, \widehat m_R(\widehat r_{\alpha, 0})) \\
       & + P[f(O, \widehat r_{\alpha}, \widehat \pi, \widehat m_R(\widehat r_{\alpha, 0})) - f(O, \widehat r_{\alpha, 0}, \widehat \pi, \widehat m_R(\widehat r_{\alpha, 0}))] \\
       & + (P_n - P) [ f(O, \widehat r_{\alpha}, \widehat \pi, \widehat m_R(\widehat r_{\alpha, 0}))  -   f(O, \widehat r_{\alpha, 0}, \widehat \pi, \widehat m_R(\widehat r_{\alpha, 0}))].
     \end{split}
 \end{equation*}
Rearranging the above and using triangle inequality we obtain the following, 
\begin{equation}
\label{eq:cmdp_main_ineq}
    \begin{split}
   &| P[f(O, \widehat r_{\alpha}, \widehat \pi, \widehat m_R(\widehat r_{\alpha, 0})) - f(O, \widehat r_{\alpha, 0}, \widehat \pi, \widehat m_R(\widehat r_{\alpha, 0}))]| \\
   = &    | P_nf(O, \widehat r_{\alpha}, \widehat \pi, \widehat m_R(\widehat r_{\alpha, 0})) |  + | P_nf(O, \widehat r_{\alpha, 0}, \widehat \pi, \widehat m_R(\widehat r_{\alpha, 0}))| \\
   &+ |(P_n - P) [ f(O, \widehat r_{\alpha}, \widehat \pi, \widehat m_R(\widehat r_{\alpha, 0}))  -   f(O, \widehat r_{\alpha, 0}, \widehat \pi, \widehat m_R(\widehat r_{\alpha, 0}))] |  \\
   \leq & | P_nf(O, \widehat r_{\alpha}, \widehat \pi, \widehat m_R(\widehat r_{\alpha, 0})) |  + | P_nf(O, \widehat r_{\alpha, 0}, \widehat \pi, \widehat m_R(\widehat r_{\alpha, 0}))| \\
   &+\mathcal C \sup_{I\in \mathcal{I} }|(P_n - P) [ \textbf{1}\{C = \infty\} \widehat \pi(X) (\textbf{1}\{R \in I\})/\sigma_I] |  \sqrt{|\widehat r_{\alpha}- \widehat r_{\alpha, 0}|}, 
    \end{split}
\end{equation}
where $\mathcal I$ is the set of all intervals in $\mathbb R$, $\sigma_I = \max\{\sqrt{\mbox{Var}_P(\textbf{1}\{C = \infty\} \widehat \pi(X) (\textbf{1}\{R \in I\}))}, \sqrt{\log n/n} \}$ for $I \in \mathcal{I}$ where $F_R(\cdot)$ is the distribution function of $R(X, Y)$, and $\mathcal C$ is a constant. The last inequality follows from the assumption that $\|m_R^*(\theta, X) - m_R^*(\theta', X) \|_2 \leq \mathfrak C_2 |\theta - \theta'|$ for $\theta, \theta' \in [r_\alpha - \eta, r_\alpha + \eta]$. From the definition of $\widehat r_{\alpha}$ we have,
\[
P_nf(O, \widehat r_{\alpha}, \widehat \pi, \widehat m_R(\widehat r_{\alpha, 0}))  = \frac{n+1}{n}\left[ \frac{1}{n+1} \sum_{i = 1}^n f(O_i, \widehat r_{\alpha}, \widehat \pi, \widehat m_R(\widehat r_{\alpha, 0})) \right] \geq \frac{\overline \pi}{n}. 
\]
Using the fact that $\sup_{\theta} J(M_n(\cdot, \widehat \pi, \widehat m_R(\widehat r_{\alpha, 0})), \theta) \leq \overline \pi/n$ we have the following upper bound,
\begin{equation*}
  P_n f(O, \widehat r_{\alpha}, \widehat \pi, \widehat m_R(\widehat r_{\alpha, 0}))  = \frac{(n+1)}{n}\left\{ \frac{\overline \pi}{n+1} + \left\{\frac{1}{n+1} \sum_{i = 1}^n f(O_i, \widehat r_{\alpha}, \widehat \pi, \widehat m_R(\widehat r_{\alpha, 0}))  - \frac{\overline \pi}{n+1}\right\}\right\} \leq \frac{2 \overline \pi}{n}. 
\end{equation*}
Therefore we have the following bounds on $P_n f(O, \widehat r_{\alpha}, \widehat \pi, \widehat m_R(\widehat r_{\alpha, 0}))$,
\[
\frac{\overline \pi}{n} \leq P_n f(O, \widehat r_{\alpha}, \widehat \pi, \widehat m_R(\widehat r_{\alpha, 0})) \leq \frac{2 \overline \pi}{n}. 
\]
Similarly we can bound $  |P_nf(O, \widehat r_{\alpha, 0}, \widehat \pi, \widehat m_R(\widehat r_{\alpha, 0})) |$ in the following manner, 
\begin{equation*}
    |P_nf(O, \widehat r_{\alpha, 0}, \widehat \pi, \widehat m_R(\widehat r_{\alpha, 0})) | \leq  |P f(O, \widehat r_{\alpha, 0}, \widehat \pi, \widehat m_R(\widehat r_{\alpha, 0})) |  + |(P_n - P) f(O, \widehat r_{\alpha, 0}, \widehat \pi, \widehat m_R(\widehat r_{\alpha, 0}))  |.
\end{equation*}
Let $\sigma_P[f(O, \theta, \widehat \pi, \widehat m_R(\theta))] = \sqrt{\mbox{Var}[f(O, \theta, \widehat \pi, \widehat m_R(\theta))]}$. Using Chebyshev's inequality we have the following with probability at least $1 - \delta$,
\[
|(P_n - P) f(O, \widehat r_{\alpha, 0}, \widehat \pi, \widehat m_R(\widehat r_{\alpha, 0}))  | \leq \frac{\sigma_P[f(O, \widehat r_{\alpha, 0}, \widehat \pi, \widehat m_R(\widehat r_{\alpha, 0}))]}{\sqrt{n \delta}} \leq \frac{\sup_{\theta} \sigma_P[f(O, \theta, \widehat \pi, \widehat m_R(\theta))]}{\sqrt{n \delta}}. 
\]
We define the following,
\[
\widehat r_{\alpha, 0}' = \inf\left\{   \theta:   \sum_{i = 1}^{n_2}\frac{ \textbf{1}\{C_i^{t_2} = 0\} [\widehat m_R(\theta, X_i^{t_2}) - (1 - \alpha)] }{n_2+1} + \frac{\textbf{1}\{C_{n+1} = 0\} [\widehat m_R(\theta, X_{n+1}) - (1 - \alpha)]}{n_2 + 1} \geq 0\right\}. 
\]
Using the right continuity of $\theta \mapsto \widehat m_R(\theta, \cdot) $ we have the following lower bound,
\begin{equation*}
    \begin{split}
       & P[\textbf{1}\{C_{n+1} = 0\} [\widehat m_R(\widehat r_{\alpha, 0}', X_{n+1}) - (1 - \alpha)]] \\
      = & \mathbb E \left[  \sum_{i = 1}^{n_2}\frac{ \textbf{1}\{C_i^{t_2} = 0\} [\widehat m_R(\widehat r_{\alpha, 0}', X_i^{t_2}) - (1 - \alpha)] }{n_2+1} + \frac{\textbf{1}\{C_{n+1} = 0\} [\widehat m_R(\widehat r_{\alpha, 0}', X_{n+1}) - (1 - \alpha)]}{n_2 + 1}\right]   \\
      \geq & 0. 
    \end{split}
\end{equation*}
We can see from the definition of $\widehat r_{\alpha, 0}$ that $\widehat r_{\alpha, 0}$ satisfies the condition of the set whose infimum is $\widehat r_{\alpha, 0}'$. Therefore we have the inequality $\widehat r_{\alpha, 0} \geq \widehat r_{\alpha, 0}'$. Using the monotonicity of the function $\textbf{1}\{C_{n+1} = 0\} [\widehat m_R(\widehat r_{\alpha, 0}', X_{n+1}) - (1 - \alpha)]$ we have,
\[
P[\textbf{1}\{C_{n+1} = 0\} [\widehat m_R(\widehat r_{\alpha, 0}, X_{n+1}) - (1 - \alpha)]] \geq P[\textbf{1}\{C_{n+1} = 0\} [\widehat m_R(\widehat r_{\alpha, 0}', X_{n+1}) - (1 - \alpha)]] \geq 0. 
\]
Similarly we can define,
\begin{equation*}
    \begin{split}
 \widehat r_{\alpha, 0}'' = & \inf\left\{   \theta: \sum_{i = 1}^{n_2}\frac{ \textbf{1}\{C_i^{t_2} = 0\} [\widehat m_R(\theta, X_i^{t_2}) - (1 - \alpha)] }{n_2+1} + \frac{\textbf{1}\{C_{n+1} = 0\} [\widehat m_R(\theta, X_{n+1}) - (1 - \alpha)]}{n_2 + 1} \right. \\
 & \quad \quad\left. - \frac{2(1-\alpha)}{n_2 +1}\geq 0\right\}.        
    \end{split}
\end{equation*}
Since $\mathbb P(J(\widehat m_R(\theta, X), \theta) > 0) = 0$ for all $\theta$ we have the following using jump lemma in \cite{angelopoulos2022conformal}, 
\[
J\left( \sum_{i = 1}^{n_2}\frac{ \textbf{1}\{C_i^{t_2} = 0\} [\widehat m_R(\theta, X_i^{t_2}) - (1 - \alpha)] }{n_2+1} + \frac{\textbf{1}\{C_{n+1} = 0\} [\widehat m_R(\theta, X_{n+1}) - (1 - \alpha)]}{n_2 + 1}, \theta\right) \stackrel{a.s.}{\leq} \frac{1}{n_2 + 1}. 
\]
Therefore we have the following upper bound,
\begin{equation*}
    \begin{split}
       &  \sum_{i = 1}^{n_2}\frac{ \textbf{1}\{C_i^{t_2} = 0\} [\widehat m_R(\widehat r_{\alpha, 0}'', X_i^{t_2}) - (1 - \alpha)] }{n_2+1} + \frac{\textbf{1}\{C_{n+1} = 0\} [\widehat m_R(\widehat r_{\alpha, 0}'', X_{n+1}) - (1 - \alpha)]}{n_2 + 1} \\
     = &\frac{2(1-\alpha)}{n_2 +1} \\
     &+ \left\{   \sum_{i = 1}^{n_2}\frac{ \textbf{1}\{C_i^{t_2} = 0\} [\widehat m_R(\widehat r_{\alpha, 0}'', X_i^{t_2}) - (1 - \alpha)] }{n_2+1} + \frac{\textbf{1}\{C_{n+1} = 0\} [\widehat m_R(\widehat r_{\alpha, 0}'', X_{n+1}) - (1 - \alpha)]}{n_2 + 1} - \frac{2(1-\alpha)}{n_2 +1} \right\} \\
     \leq & \frac{2(1-\alpha)}{n_2 +1} + \frac{1}{n_2 + 1} \\
     = & \frac{3 - 2\alpha}{n_2 + 1}. 
    \end{split}
\end{equation*}
Using the right continuity of $\theta \mapsto \widehat m_R(\theta, \cdot) $ we can see that $\widehat r_{\alpha, 0}''$ satisfies the condition of the set whose infimum is $\widehat r_{\alpha, 0}$,
\begin{equation*}
    \begin{split}
        &  \sum_{i = 1}^{n_2}\frac{ \textbf{1}\{C_i^{t_2} = 0\} [\widehat m_R(\widehat r_{\alpha, 0}'', X_i^{t_2}) - (1 - \alpha)] }{n_2+1} - \frac{(1 - \alpha)}{n_2 + 1} \\
     = &  \sum_{i = 1}^{n_2}\frac{ \textbf{1}\{C_i^{t_2} = 0\} [\widehat m_R(\widehat r_{\alpha, 0}'', X_i^{t_2}) - (1 - \alpha)] }{n_2+1} + \frac{\textbf{1}\{C_{n+1} = 0\} [\widehat m_R(\widehat r_{\alpha, 0}'', X_{n+1}) - (1 - \alpha)]}{n_2 + 1} \\
     & -  \frac{\textbf{1}\{C_{n+1} = 0\} [\widehat m_R(\widehat r_{\alpha, 0}'', X_{n+1}) - (1 - \alpha)]}{n_2 + 1}    - \frac{(1 - \alpha)}{n_2 + 1}  \\
     \geq &  \sum_{i = 1}^{n_2}\frac{ \textbf{1}\{C_i^{t_2} = 0\} [\widehat m_R(\widehat r_{\alpha, 0}'', X_i^{t_2}) - (1 - \alpha)] }{n_2+1} + \frac{\textbf{1}\{C_{n+1} = 0\} [\widehat m_R(\widehat r_{\alpha, 0}'', X_{n+1}) - (1 - \alpha)]}{n_2 + 1} \\
     & - \frac{2(1 - \alpha)}{n_2 + 1} \\
     \geq & 0. 
    \end{split}
\end{equation*}
Therefore $\widehat r_{\alpha, 0} \leq \widehat r_{\alpha, 0}''$ and hence monotonicity of the function $\textbf{1}\{C_{n+1} = 0\} [\widehat m_R(\widehat r_{\alpha, 0}', X_{n+1}) - (1 - \alpha)]$ implies the following upper bound,
\begin{equation*}
    \begin{split}
     & P[\textbf{1}\{C_{n+1} = 0\} [\widehat m_R(\widehat r_{\alpha, 0}, X_{n+1}) - (1 - \alpha)]]  \\
     \leq & P[\textbf{1}\{C_{n+1} = 0\} [\widehat m_R(\widehat r_{\alpha, 0}'', X_{n+1}) - (1 - \alpha)]]  \\
     = & \mathbb E \left[ \sum_{i = 1}^{n_2}\frac{ \textbf{1}\{C_i^{t_2} = 0\} [\widehat m_R(\widehat r_{\alpha, 0}'', X_i^{t_2}) - (1 - \alpha)] }{n_2+1} + \frac{\textbf{1}\{C_{n+1} = 0\} [\widehat m_R(\widehat r_{\alpha, 0}'', X_{n+1}) - (1 - \alpha)]}{n_2 + 1} \right] \\
     \leq & \frac{3 - 2\alpha}{n_2 + 1}. 
    \end{split}
\end{equation*}
Using this bound we obtain the following,
\begin{equation*}
    \begin{split}
        & |P f(O, \widehat r_{\alpha, 0}, \widehat \pi, \widehat m_R(\widehat r_{\alpha, 0})) | \\
        \leq &  |P f(O, \widehat r_{\alpha, 0},  \pi^*,  m_R^*(\widehat r_{\alpha, 0})) | +  \|\widehat \pi  - \pi^*\|_2 \sup_{\theta} \| \widehat m_R (\theta, \cdot) - m_R^*(\theta, \cdot)\|_2 \\
        = & |P [\textbf{1}\{C_{n+1} = 0\} [ m_R^*(\widehat r_{\alpha, 0}, X_{n+1}) - (1 - \alpha)]]  | +  \|\widehat \pi  - \pi^*\|_2 \sup_{\theta} \| \widehat m_R (\theta, \cdot) - m_R^*(\theta, \cdot)\|_2 \\ 
        \leq & | P[\textbf{1}\{C_{n+1} = 0\} [\widehat m_R(\widehat r_{\alpha, 0}, X_{n+1}) - (1 - \alpha)]] | + (1 + \|\widehat \pi  - \pi^*\|_2) \sup_{\theta} \| \widehat m_R (\theta, \cdot) - m_R^*(\theta, \cdot)\|_2 \\
        \leq & \frac{3 - 2\alpha}{n_2 + 1} + (1 + \|\widehat \pi  - \pi^*\|_2) \sup_{\theta} \| \widehat m_R (\theta, \cdot) - m_R^*(\theta, \cdot)\|_2. 
    \end{split}
\end{equation*}
The last step is to bound the empirical process term,
\begin{equation*}
    \begin{split}
  \sup_{I\in \mathcal{I} }|(P_n - P) [ \textbf{1}\{C = \infty\} \widehat \pi(X) (\textbf{1}\{R \in I\})/\sigma_I] |. 
    \end{split}
\end{equation*}
Here $\mathcal{I}$ is the set of all intervals in $\mathbb R$. For every interval $I\subseteq\mathbb R$, we define,
\[
    g_I(O)
    = \1\{C=\infty\}\widehat\pi(X)\1\{R\in I\}.
\]
We note that,
\[
    v_I = \Var_P(g_I(O)),
    \qquad
    \sigma_I = \max\left\{\sqrt{v_I},\sqrt{\frac{\log n}{n}}\right\}.
\]
The ratio-type empirical process of interest is
\[
    Z_n
    = \sup_{I\subseteq\mathbb R\text{ interval}}
    \frac{\sqrt n\,\left|(P_n-P)g_I\right|}{\sigma_I}.
\]
Here and below $\calC$ denotes a positive constant whose value may change from line to line.  The constant $\calC$ does not depend on $n$, $\delta$, $P$, $I$, or $\widehat\pi$, except through the displayed occurrence of $\overline\pi$.

\begin{lemma}[High-probability bound on ratio type empirical process term]
There exists a universal constant $\calC>0$ such that, for all $n\geq 3$ and all $\delta\in(0,1)$, with probability at least $1-\delta$,
\[
    Z_n
    \leq
    \calC\left\{
        \sqrt{\log n + \log\log(en) + \log\frac{1}{\delta}}
        +
        \overline\pi\,
        \frac{\log n + \log\log(en) + \log\frac{1}{\delta}}{\sqrt{\log n}}
    \right\}.
\]
Consequently, since $\log\log(en)\leq \log n$ for all $n\geq 3$, the simpler bound
\[
    Z_n
    \leq
    \calC(1+\overline\pi)
    \left\{
        \sqrt{\log n + \log\frac{1}{\delta}}
        +
        \frac{\log(1/\delta)}{\sqrt{\log n}}
    \right\}
\]
holds with probability at least $1-\delta$.  
\end{lemma}

\begin{proof}
For an interval $I$, we write,
\[
    S_I = \sum_{i=1}^n \{g_I(O_i)-Pg_I\}
    = n(P_n-P)g_I.
\]
Since $0\leq g_I\leq \overline\pi$, we have,
\[
    |g_I(O)-Pg_I|\leq \overline\pi
    \qquad P\text{-almost surely},
\]
and,
\[
    \sum_{i=1}^n \Var\{g_I(O_i)\} = n v_I.
\]
Bernstein's inequality therefore gives, for every fixed interval $I$ and every $y>0$,
\[
    \Pbb\left(
        |S_I|>
        \sqrt{2n v_I y} + \frac{2\overline\pi y}{3}
    \right)
    \leq 2e^{-y}.
\]
Since $\sigma_I^2\geq v_I$ and $\sigma_I\geq \sqrt{\log n/n}$, division by $\sqrt n\sigma_I$ yields,
\[
    \Pbb\left(
        \frac{\sqrt n |(P_n-P)g_I|}{\sigma_I}
        >
        \sqrt{2y}
        +
        \frac{2\overline\pi y}{3\sqrt{\log n}}
    \right)
    \leq 2e^{-y}.
\]
The rest of the proof uses a peeling argument over the possible variance sizes of $g_I$ and a union bound over the interval traces. Let us set,
\[
    a_n=\frac{\log n}{n}.
\]
For integers $k\geq 0$, we define the variance slices,
\[
    \mathcal I_0
    =
    \{I: v_I\leq a_n\},\quad \mathcal I_k
    =
    \{I:2^{k-1}a_n < v_I\leq 2^k a_n\} \mbox{  for  } k \geq 1.
\]
Because $0\leq g_I\leq \overline\pi$, we have $v_I\leq \overline\pi^2$.  Hence only those $k$ satisfying $2^{k-1}a_n\leq \overline\pi^2$ are relevant.  The number of relevant slices is at most
\[
    K_n
    \leq
    1+\frac{\log(en\overline\pi^2/\log n)}{\log 2}.
\]
If $\overline\pi<1$, the right hand side can be replaced by the same expression with $\overline\pi$ set equal to one.  Thus, after enlarging the universal constant, it is enough to use the bound $ K_n\leq \calC\log(en).$ For a fixed realization $O_1,\ldots,O_n$, the sets $\{i:R_i\in I\}$ generated by intervals $I\subseteq\mathbb R$ are contiguous blocks after ordering the observed values $R_1,\ldots,R_n$.  Therefore the number of distinct subsets of $\{1,\ldots,n\}$ obtained by intersecting the sample with intervals is no larger than,
\[
    1+\frac{n(n+1)}{2}\leq n^2,
    \qquad n\geq 2.
\]
Equivalently, on the sample points, the class $\{\1\{R\in I\}:I\text{ interval}\}$ has growth at most $n^2$.  Multiplication by the fixed weights $\1\{C=\infty\}\widehat\pi(X)$ does not increase this trace count.  
We apply the preceding Bernstein inequality to every interval trace in every slice. We choose,
\[
    y
    =
    \log\left(\frac{2 n^2 K_n}{\delta}\right).
\]
A union bound over the at most $K_n$ variance slices and the at most $n^2$ interval traces in each slice gives,
\[
\begin{aligned}
    &\Pbb\left(
        \exists I:
        \frac{\sqrt n |(P_n-P)g_I|}{\sigma_I}
        >
        \sqrt{2y}
        +
        \frac{2\overline\pi y}{3\sqrt{\log n}}
    \right) \leq
        2 n^2 K_n e^{-y}
        = \delta.
\end{aligned}
\]
Therefore, with probability at least $1-\delta$,
\[
    Z_n
    \leq
    \sqrt{2y}
    +
    \frac{2\overline\pi y}{3\sqrt{\log n}}.
\]
Since $K_n\leq \calC\log(en)$, the definition of $y$ implies,
\[
    y
    \leq
    \calC\left(
        \log n + \log\log(en) + \log\frac{1}{\delta}
    \right).
\]
Substituting this bound for $y$ into the previous display givem
\[
    Z_n
    \leq
    \calC\left\{
        \sqrt{\log n + \log\log(en) + \log\frac{1}{\delta}}
        +
        \overline\pi\,
        \frac{\log n + \log\log(en) + \log\frac{1}{\delta}}{\sqrt{\log n}}
    \right\}
\]
with probability at least $1-\delta$, which proves the first claim.  The simplified display follows from elementary comparisons of logarithms for $n\geq 3$.
\end{proof}

Combining all the prior computations we have the following with probability at-least $1 - 2\delta$, 
\begin{equation*}
    \begin{split}
       & | P[f(O, \widehat r_{\alpha}, \widehat \pi, \widehat m_R(\widehat r_{\alpha, 0})) - f(O, \widehat r_{\alpha, 0}, \widehat \pi, \widehat m_R(\widehat r_{\alpha, 0}))]| \\
  \leq & | P_nf(O, \widehat r_{\alpha}, \widehat \pi, \widehat m_R(\widehat r_{\alpha, 0})) |  + | P_nf(O, \widehat r_{\alpha, 0}, \widehat \pi, \widehat m_R(\widehat r_{\alpha, 0}))| \\
   &+ \mathcal C\sup_{I\in \mathcal{I} }|(P_n - P) [ \textbf{1}\{C = \infty\} \widehat \pi(X) (\textbf{1}\{R \in I\})/\sigma_I] |  \sqrt{|\widehat r_{\alpha}- \widehat r_{\alpha, 0}|} \\
   \leq & \frac{2 \overline \pi}{n} + \frac{3 - 2\alpha}{n_2 + 1} + (1 + \|\widehat \pi  - \pi^*\|_2) \sup_{\theta} \| \widehat m_R (\theta, \cdot) - m_R^*(\theta, \cdot)\|_2 \\
   &+  \frac{\sup_\theta \sigma_P[f(O, \theta, \widehat \pi, \widehat m_R(\theta))]}{\sqrt{n \delta}}  + \mathcal C \left\{\sqrt{\frac{\log n}{n}} + \sqrt{\frac{\log (1/\delta)}{n}} \right\} \sqrt{ | \widehat r_{\alpha} - \widehat r_{\alpha, 0}|}. 
    \end{split}
\end{equation*}
From the assumptions of the proposition we know that,
\begin{equation*}
    \begin{split}
        & | P[f(O, \widehat r_{\alpha}, \widehat \pi, \widehat m_R(\widehat r_{\alpha, 0})) - f(O, \widehat r_{\alpha, 0}, \widehat \pi, \widehat m_R(\widehat r_{\alpha, 0}))]| \\
        = & |P[\textbf{1}\{C = \infty\} \widehat \pi(X) [\textbf{1}\{R \leq \widehat r_{\alpha} \} - \textbf{1}\{R \leq \widehat r_{\alpha, 0} \}]] | \\
        \geq &  \mathfrak C_1 | \widehat r_{\alpha} - \widehat r_{\alpha, 0}|.
    \end{split}
\end{equation*}
Using this inequality we have the following with probability at-least $1 - 2\delta$,
\begin{equation*}
    \begin{split}
        &\mathfrak C_1 | \widehat r_{\alpha} - \widehat r_{\alpha, 0}| - \mathcal C \left\{\sqrt{\frac{\log n}{n}} + \sqrt{\frac{\log (1/\delta)}{n}} \right\}\sqrt{ | \widehat r_{\alpha} - \widehat r_{\alpha, 0}|} \\
        \leq & \frac{2 \overline \pi}{n} + \frac{3 - 2\alpha}{n_2 + 1} + (1 + \|\widehat \pi  - \pi^*\|_2) \sup_{\theta} \| \widehat m_R (\theta, \cdot) - m_R^*(\theta, \cdot)\|_2 +  \frac{\sup_\theta \sigma_P[f(O, \theta, \widehat \pi, \widehat m_R(\theta))]}{\sqrt{n \delta}} .
    \end{split}
\end{equation*}
This implies that the following bound holds with probability at-least $1 - 2\delta$,
\begin{equation*}
    \begin{split}
        | \widehat r_{\alpha} - \widehat r_{\alpha, 0}|\leq  & \mathcal{C} \frac{\log n + \log(1/\delta)}{n} + \frac{2 \overline \pi}{n} + \frac{3 - 2\alpha}{n_2 + 1} \\
        & + (1 + \|\widehat \pi  - \pi^*\|_2) \sup_{\theta} \| \widehat m_R (\theta, \cdot) - m_R^*(\theta, \cdot)\|_2 +  \frac{\sup_\theta \sigma_P[f(O, \theta, \widehat \pi, \widehat m_R(\theta))]}{\sqrt{n \delta}} .
    \end{split}
\end{equation*}
The proof of the proposition is completed by using the Lipschitz property of $m_R^*(\theta, \cdot)$ in the neighborhood of $r_{\alpha}$.

\end{proof}

\section{Proof of \Cref{thm:gen_influ_func}}
\label{app:proof_gen_influ}
Throughout the proof we condition on the training data
$\mathcal D^{\mathrm{tr}}$.  Consequently, we treat the estimated nuisance functions $\widehat\eta=\{\widehat\eta_i\}_{i=1}^K$ as fixed. To simplify notation, we write $f_\theta(O)
    :=
    f(O,\theta,\widehat\eta)$, $ f_\theta^*(O)
    :=
    f(O,\theta,\eta^*)$, $d=d(\eta^*,\widehat\eta)$, and $ d'=d'(\eta^*,\widehat\eta).$ We also define $c_\alpha:=1-\alpha$, $g_\theta(O)
    :=
    f_\theta(O)-f_{r_\alpha}(O)$, and $\mu(\theta)
    :=
    |Pg_\theta|$. By Assumption~\ref{assump:gen_identification}, $Pf_{r_\alpha}^*
    =
    F_Y(r_\alpha)
    =
    c_\alpha.$ Therefore Assumption~\ref{assump:gen_nuisance} implies that, 
\begin{equation}
\label{eq:rsp_root_bias}
    |Pf_{r_\alpha}-c_\alpha|
    \leq d.
\end{equation}

We record two consequences of the bounded-variation assumption. For each $o$, the function $\theta\mapsto f(o,\theta,\widehat\eta)$ has total variation at most $V$.  By the Jordan decomposition theorem,
there exist functions $f^+_\theta(o)$ and $f^-_\theta(o)$ such that
\[
    f_\theta(o)
    =
    f_{-\infty}(o)
    +
    f^+_\theta(o)
    -
    f^-_\theta(o),
\]
where both $\theta\mapsto f^+_\theta(o)$ and $\theta\mapsto f^-_\theta(o)$ are non-decreasing and bounded by $V$. Consequently, each of the two classes $\{f^+_\theta:\theta\in\mathbb R\}$ and $ \{f^-_\theta:\theta\in\mathbb R\}$ is pointwise ordered. A pointwise ordered class is a VC-subgraph
class with a universal finite VC index.  Since sums and differences
of a fixed number of VC-type classes remain VC-type, the class $\mathcal G
    :=
    \{
        g_\theta:
        \theta\in\mathbb R
    \}$ satisfies,
\begin{equation}
\label{eq:bv_entropy}
    \sup_Q
    \log
    N\left(
        \epsilon C_{BV},
        \mathcal G,
        L_2(Q)
    \right)
    \leq
    C_{BV}
    \log\left(\frac{C_{BV}}{\epsilon}\right),
    \qquad
    0<\epsilon\leq1,
\end{equation}
for a constant $C_{BV}$ depending only on $B$ and $V$.
Thus the bounded-variation assumption provides the
logarithmic entropy control needed for bounding the coverage. Recall that,
\[
    \widehat r_\alpha
    =
    \inf\{
        \theta:
        P_nf_\theta\geq c_\alpha
    \}.
\]
By right continuity, $  P_nf_{\widehat r_\alpha}
    \geq c_\alpha.$ For every $\theta<\widehat r_\alpha$,
$P_nf_\theta<c_\alpha$.  Hence the amount by which the empirical
estimating equation can exceed $c_\alpha$ at its first crossing is at
most the size of its jump at $\widehat r_\alpha$.

By Assumption~\ref{assump:gen_bv}, every individual path has jumps of
size at most $V$.  Moreover, for every fixed threshold, the
probability that an observation has a jump there is zero.  Since the observations the calibration data-set are independent, almost surely no two
observations jump at the same threshold.  Therefore we have the following bound,
\begin{equation}
\label{eq:rsp_empirical_overshoot}
    0
    \leq
    P_nf_{\widehat r_\alpha}-c_\alpha
    \leq
    \frac{V}{n}.
\end{equation}

We next control the empirical process evaluated at the random
threshold $\widehat r_\alpha$. We set $\varepsilon_n=n^{-1/2}$. We claim that,
\begin{equation}
\label{eq:rsp_ratio_bound}
\mathbb E
\left[
    \sup_{\theta\in\mathbb R}
    \frac{
        |(P_n-P)g_\theta|^2
    }{
        \mu(\theta)+\varepsilon_n
    }
    \,
\right]
\leq
\mathfrak C(\log(en))^2
\left(
    \frac1n
    +
    \frac{d'}{n\varepsilon_n}
\right).
\end{equation}
To establish this, we observe that $ |g_\theta(O)|\leq2B$. Let $q_0=2B$ and $ q_j=2^{-j}q_0$. We define the shells,
\[
    A_j
    =
    \{
        \theta:
        q_{j+1}<\mu(\theta)\leq q_j
    \}.
\]
We let $J_n
    =
    \min\{
        j:q_{j+1}\leq\varepsilon_n
    \}.$ Since $\varepsilon_n=n^{-1/2}$, $J_n=O(\log n).$ On the shell $A_j$,
Assumption~\ref{assump:gen_local_variance} gives,
\[
    Pg_\theta^2
    \leq
    Kq_j+d'
    =:
    \sigma_j^2.
\]
Combining this variance bound with
\eqref{eq:bv_entropy}, the standard localized maximal inequality
for bounded VC-type classes gives,
\begin{equation}
\label{eq:rsp_localized_maximal}
\begin{split}
&
\mathbb E
\left[
    \sup_{\theta\in A_j}
    |(P_n-P)g_\theta|^2
    \,
\right]
\leq
\mathfrak C
\left[
    \frac{\sigma_j^2}{n}
    \left\{
        1+
        \log\left(
            \frac{\mathfrak C}{\sigma_j}
        \right)
    \right\}
    +
    \frac1{n^2}
    \left\{
        1+
        \log\left(
            \frac{\mathfrak C}{\sigma_j}
        \right)
    \right\}^2
\right].
\end{split}
\end{equation}
For all shells with $j\leq J_n$, the logarithmic factors in
\eqref{eq:rsp_localized_maximal} are bounded by
$\mathfrak C\log(en)$.  Hence
\[
\mathbb E
\left[
    \sup_{\theta\in A_j}
    |(P_n-P)g_\theta|^2
    \,
\right]
\leq
\mathfrak C(\log(en))^2
\left[
    \frac{q_j+d'}{n}
    +
    \frac1{n^2}
\right].
\]

For $\theta\in A_j$, we have $\mu(\theta)+\varepsilon_n
    \geq q_{j+1}.$ Therefore,
\begin{equation*}
\begin{split}
&
\mathbb E
\left[
    \sup_{\mu(\theta)>\varepsilon_n}
    \frac{|(P_n-P)g_\theta|^2}
         {\mu(\theta)+\varepsilon_n}
    \,
\right]
\\
&\leq
\sum_{j=0}^{J_n}
\frac{
\mathbb E
[
    \sup_{\theta\in A_j}
    |(P_n-P)g_\theta|^2
]
}{q_{j+1}}
\\
&\leq
\mathfrak C(\log(en))^2
\left[
    \frac1n
    +
    \frac{d'}{n\varepsilon_n}
\right].
\end{split}
\end{equation*}
Here we have used the following,
\[
    \frac{q_j}{q_{j+1}}=2,
    \qquad
    \sum_{j=0}^{J_n}\frac1{q_{j+1}}
    \leq
    \frac{\mathfrak C}{\varepsilon_n},
\]
and absorbed the extra logarithmic factor into the generic constant
$\mathfrak C(\log(en))^2$. For the innermost region
$\{\theta:\mu(\theta)\leq\varepsilon_n\}$,
Assumption~\ref{assump:gen_local_variance} gives,
\[
    Pg_\theta^2
    \leq
    K\varepsilon_n+d'.
\]
Applying the same localized maximal inequality yields,
\[
\frac1{\varepsilon_n}
\mathbb E
\left[
    \sup_{\mu(\theta)\leq\varepsilon_n}
    |(P_n-P)g_\theta|^2
    \,
\right]
\leq
\mathfrak C(\log(en))^2
\left[
    \frac1n
    +
    \frac{d'}{n\varepsilon_n}
\right].
\]
Combining the two regions proves \eqref{eq:rsp_ratio_bound}. We first bound $\mu(\widehat r_\alpha)$. By the triangle inequality,
\begin{equation}
\label{eq:rsp_mu_decomp}
\begin{split}
    \mu(\widehat r_\alpha)
    =
    |Pg_{\widehat r_\alpha}|
    \leq&
    |(P-P_n)g_{\widehat r_\alpha}|
    +
    |P_nf_{\widehat r_\alpha}-c_\alpha|
+
    |P_nf_{r_\alpha}-c_\alpha|.
\end{split}
\end{equation}
Let $F\equiv 2B$ denote an envelope of $\mathcal G$.  By
\eqref{eq:bv_entropy}, after changing constants if necessary,
\[
\sup_Q
\log N\left(
    \epsilon\|F\|_{Q,2},
    \mathcal G,
    L_2(Q)
\right)
\leq
\mathfrak C
\log\left(\frac{\mathfrak C}{\epsilon}\right),
\qquad
0<\epsilon\leq1.
\]
Consequently, the uniform entropy integral satisfies,
\[
J(1,\mathcal G)
:=
\sup_Q
\int_0^1
\sqrt{
1+
\log N\left(
    \epsilon\|F\|_{Q,2},
    \mathcal G,
    L_2(Q)
\right)
}
\,d\epsilon
\leq
\mathfrak C
\int_0^1
\sqrt{
1+\log(\mathfrak C/\epsilon)
}
\,d\epsilon
<\infty.
\]
The standard uniform-entropy maximal inequality for empirical
processes therefore yields,
\[
\mathbb E
\left[
    \sup_{\theta\in\mathbb R}
    \left|
        \sqrt n(P_n-P)g_\theta
    \right|
\right]
\leq
\mathfrak C J(1,\mathcal G)\|F\|_{P,2}
\leq
\mathfrak C.
\]
Hence we have,
\[
\mathbb E
\left[
    \sup_{\theta\in\mathbb R}
    |(P_n-P)g_\theta|
\right]
\leq
\frac{\mathfrak C}{\sqrt n}.
\]
By \eqref{eq:rsp_empirical_overshoot}, we have the bound,
\[
    |P_nf_{\widehat r_\alpha}-c_\alpha|
    \leq
    \frac{V}{n}.
\]
Finally we have,
\begin{equation*}
\begin{split}
\mathbb E
\left[
    |P_nf_{r_\alpha}-c_\alpha|
    \,
\right]
&\leq
|Pf_{r_\alpha}-c_\alpha|
+
\mathbb E
\left[
    |(P_n-P)f_{r_\alpha}|
    \,
\right] \leq
d+\frac{\mathfrak C}{\sqrt n},
\end{split}
\end{equation*}
where we used \eqref{eq:rsp_root_bias} in the last inequality. Combining these bounds with \eqref{eq:rsp_mu_decomp},
\begin{equation}
\label{eq:rsp_mu_bound}
\mathbb E
\left[
    \mu(\widehat r_\alpha)
    \,
\right]
\leq
d+\frac{\mathfrak C}{\sqrt n}.
\end{equation}

We now apply Cauchy-Schwarz inequality,
\begin{equation*}
\begin{split}
&
\mathbb E
\left[
    |(P_n-P)g_{\widehat r_\alpha}|
    \,
\right]\leq
\left\{
\mathbb E
\left[
    \mu(\widehat r_\alpha)+\varepsilon_n
    \,
\right]
\right\}^{1/2}\times
\left\{
\mathbb E
\left[
    \sup_\theta
    \frac{|(P_n-P)g_\theta|^2}
         {\mu(\theta)+\varepsilon_n}
    \,
\right]
\right\}^{1/2}.
\end{split}
\end{equation*}
Using \eqref{eq:rsp_ratio_bound},
\eqref{eq:rsp_mu_bound}, and
$\varepsilon_n=n^{-1/2}$,
we obtain,
\begin{equation}
\label{eq:rsp_random_ep_bound}
\begin{split}
&
\mathbb E
\left[
    |(P_n-P)g_{\widehat r_\alpha}|
    \,
\right] \leq
\mathfrak C_n
\left(
    \sqrt d+n^{-1/4}
\right)
\left(
    n^{-1/2}
    +
    n^{-1/4}\sqrt{d'}
\right),
\end{split}
\end{equation}
where $\mathfrak C_n=\mathfrak C\log(en).$ For ease of notation, let us define,
\[
    \rho_n(d,d')
    =
    \mathfrak C_n
    \left(
        \sqrt d+n^{-1/4}
    \right)
    \left(
        n^{-1/2}
        +
        n^{-1/4}\sqrt{d'}
    \right).
\]
We have the following decomposition,
\begin{equation}
\label{eq:rsp_pop_decomp}
\begin{split}
Pf_{\widehat r_\alpha}-c_\alpha
={}&
P_nf_{\widehat r_\alpha}-c_\alpha
+
(P-P_n)f_{\widehat r_\alpha}
\\
={}&
P_nf_{\widehat r_\alpha}-c_\alpha
+
(P-P_n)g_{\widehat r_\alpha}
+
(P-P_n)f_{r_\alpha}.
\end{split}
\end{equation}
Since $f_{r_\alpha}$ is fixed conditional on the training data,
\[
    \mathbb E
    \left[
        (P-P_n)f_{r_\alpha}
    \right]
    =0.
\]
Taking conditional expectations in \eqref{eq:rsp_pop_decomp}, and
using \eqref{eq:rsp_empirical_overshoot} and
\eqref{eq:rsp_random_ep_bound}, gives,
\begin{equation}
\label{eq:rsp_pop_lower}
\mathbb E
\left[
    Pf_{\widehat r_\alpha}-c_\alpha
    \,
\right]
\geq
-\rho_n(d,d'),
\end{equation}
and,
\begin{equation}
\label{eq:rsp_pop_upper}
\mathbb E
\left[
    Pf_{\widehat r_\alpha}-c_\alpha
    \,
\right]
\leq
\frac{V}{n}
+
\rho_n(d,d').
\end{equation}
Because $\widehat r_\alpha$ is a function only of the training and
calibration samples, it is independent of the test point conditional
on $\mathcal D^{\mathrm{tr}}$ and the calibration data.  Hence
\begin{equation*}
\begin{split}
&
\mathbb P
\left(
    Y_{n+1}\leq\widehat r_\alpha   
\right)=
\mathbb E
\left[
    F_Y(\widehat r_\alpha)
    \,
\right]=
\mathbb E
\left[
    Pf_{\widehat r_\alpha}^*
    \,
\right].
\end{split}
\end{equation*}
By Assumption~\ref{assump:gen_nuisance},
for every realization of $\widehat r_\alpha$,
\[
    \left|
        Pf_{\widehat r_\alpha}^*
        -
        Pf_{\widehat r_\alpha}
    \right|
    \leq d.
\]
Combining this inequality with
\eqref{eq:rsp_pop_lower} gives,
\[
\mathbb P
\left(
    Y_{n+1}\leq\widehat r_\alpha
\right)
\geq
1-\alpha-d-\rho_n(d,d').
\]
Similarly, \eqref{eq:rsp_pop_upper} gives,
\[
\mathbb P
\left(
    Y_{n+1}\leq\widehat r_\alpha
\right)
\leq
1-\alpha
+
d
+
\frac{V}{n}
+
\rho_n(d,d').
\]
Finally, the term $V/n$ can be absorbed into $\rho_n(d,d')$ after
enlarging the constant $\mathfrak C$, because $\rho_n(d,d')
    \geq
    \mathfrak C_n n^{-3/4}$ and $n^{-1}\leq n^{-3/4}$ for $n\geq1$. Therefore we have,
\[
    1-\alpha-\Delta_{\rsp,n}
    \leq
    \mathbb P
    \left(
        Y_{n+1}\leq\widehat r_\alpha
    \right)
    \leq
    1-\alpha+\Delta_{\rsp,n},
\]
with $\Delta_{\rsp,n}$ defined in
\eqref{eq:rsp_general_slack}.  This completes the proof.

\section{Proof of \Cref{thm:gen_three_split}}
\label{app:proof_gen_three_split}
As stated before the theorem, throughout the proof all probabilities
and expectations are conditional on $\mathcal T
    =
    \mathcal D_1^{\mathrm{tr}}
    \cup
    \mathcal D_2^{\mathrm{tr}}.$ For brevity, write $s_0
    :=
    \widehat r_{\alpha,0},$ and recall the definitions,
\[
\begin{split}
    f^*_\theta(O)
    &=
    f\left(
        O,\theta,
        \{\eta_{1i}^*(\theta)\}_{i=1}^{K_1},
        \{\eta_{2i}^*\}_{i=1}^{K_2}
    \right),\qquad \widehat f_{\theta,s_0}(O)
    =
    f\left(
        O,\theta,
        \{\widehat\eta_{1i}(s_0)\}_{i=1}^{K_1},
        \{\widehat\eta_{2i}\}_{i=1}^{K_2}
    \right),
\end{split}
\]
We further recall that $\phi_{\theta,s_0}(O)
    =
    \widehat f_{\theta,s_0}(O)
    -(1-\alpha)$ and $ \Delta_{\mathrm{init}}
    =
    d_{\mathrm{init}}(s_0,r_\alpha)$. By Assumption~\ref{assump:gen_identification}, we have the following equality,
\begin{equation}
\label{eq:rscp_proof_identification}
    P f^*_\theta
    =
    \mathbb P(Y\leq\theta)
    \qquad
    \mbox{for every }\theta\in\mathbb R.
\end{equation}
Moreover, Assumption~\ref{assump:rscp_anchored_bias} gives the following, uniformly
over $\theta$,
\begin{equation}
\label{eq:rscp_proof_bias}
    \left|
        P f^*_\theta
        -
        P\widehat f_{\theta,s_0}
    \right|
    \leq
    d_{\rscp}
    +
    \Delta_{\mathrm{init}}.
\end{equation}

We first prove the lower bound on coverage. Introduce the hypothetical full-sample threshold,
\begin{equation}
\label{eq:rscp_rminus}
    \widehat r_\alpha^-
    :=
    \inf\left\{
        \theta\in\mathbb R:
        \frac{1}{n+1}
        \sum_{i=1}^{n+1}
        \phi_{\theta,s_0}(O_i)
        \geq0
    \right\}.
\end{equation}
By the right continuity in
Assumption~\ref{assump:rscp_monotone},
\begin{equation}
\label{eq:rscp_rminus_nonnegative}
    \frac{1}{n+1}
    \sum_{i=1}^{n+1}
    \phi_{\widehat r_\alpha^-,s_0}(O_i)
    \geq0.
\end{equation}
We first compare $\widehat r_\alpha^-$ with the actual $\rscp$
threshold $\widehat r_\alpha$. By the definition \eqref{eq:rscp_general_centered},
\[
    \frac{
        \sum_{i=1}^n
        \phi_{\widehat r_\alpha,s_0}(O_i)
        +
        \inf_o
        \phi_{\widehat r_\alpha,s_0}(o)
    }{n+1}
    \geq0.
\]
Since $\phi_{\widehat r_\alpha,s_0}(O_{n+1})
    \geq
    \inf_o
    \phi_{\widehat r_\alpha,s_0}(o)$, it follows that,
\[
\begin{split}
    \frac{1}{n+1}
    \sum_{i=1}^{n+1}
    \phi_{\widehat r_\alpha,s_0}(O_i)
    &=
    \frac{
        \sum_{i=1}^{n}
        \phi_{\widehat r_\alpha,s_0}(O_i)
        +
        \phi_{\widehat r_\alpha,s_0}(O_{n+1})
    }{n+1}\\
    &\geq
    \frac{
        \sum_{i=1}^{n}
        \phi_{\widehat r_\alpha,s_0}(O_i)
        +
        \inf_o
        \phi_{\widehat r_\alpha,s_0}(o)
    }{n+1}\\
    &\geq0.
\end{split}
\]
Therefore $\widehat r_\alpha$ belongs to the set whose infimum defines
$\widehat r_\alpha^-$, and hence we have the inequality,
\begin{equation}
\label{eq:rscp_order_lower}
    \widehat r_\alpha^-
    \leq
    \widehat r_\alpha.
\end{equation}
The threshold $\widehat r_\alpha^-$ is a symmetric function of
$O_1,\ldots,O_{n+1}$. Let $E=\{O_1,\ldots,O_{n+1}\}$ denote the unordered collection of these observations. By
exchangeability,
\begin{equation*}
\begin{split}
&
\mathbb E
\left[
    \phi_{\widehat r_\alpha^-,s_0}(O_{n+1})
    \mid E
\right] =
\frac{1}{n+1}
\sum_{i=1}^{n+1}
\phi_{\widehat r_\alpha^-,s_0}(O_i).
\end{split}
\end{equation*}
Combining this identity with
\eqref{eq:rscp_rminus_nonnegative} and then applying the tower
property gives,
\begin{equation}
\label{eq:rscp_rminus_expectation}
    \mathbb E
    \left[
        \phi_{\widehat r_\alpha^-,s_0}(O_{n+1})
    \right]
    \geq0.
\end{equation}
By Assumption~\ref{assump:rscp_monotone} and
\eqref{eq:rscp_order_lower},
\[
    \phi_{\widehat r_\alpha^-,s_0}(O_{n+1})
    \leq
    \phi_{\widehat r_\alpha,s_0}(O_{n+1})
\]
almost surely. Therefore we have the following,
\begin{equation}
\label{eq:rscp_rhat_frozen_lower}
    \mathbb E
    \left[
        \phi_{\widehat r_\alpha,s_0}(O_{n+1})
    \right]
    \geq0.
\end{equation}
Now $\widehat r_\alpha$ depends only on the calibration data and the
two training samples and is therefore independent of the test
observation $O_{n+1}$.  Consequently,
\begin{equation*}
\begin{split}
\mathbb E
\left[
    \phi_{\widehat r_\alpha,s_0}(O_{n+1})
\right]
&=
\mathbb E
\left[
    P\widehat f_{\widehat r_\alpha,s_0}
    -(1-\alpha)
\right].
\end{split}
\end{equation*}
Thus \eqref{eq:rscp_rhat_frozen_lower} implies,
\begin{equation}
\label{eq:rscp_pfrozen_lower}
    \mathbb E
    \left[
        P\widehat f_{\widehat r_\alpha,s_0}
    \right]
    \geq
    1-\alpha.
\end{equation}
On the other hand, by
\eqref{eq:rscp_proof_identification},
\begin{equation}
\label{eq:rscp_coverage_identification}
\begin{split}
    \mathbb P
    \left(
        Y_{n+1}\leq\widehat r_\alpha
    \right)
    &=
    \mathbb E
    \left[
        \mathbb P
        \left(
            Y_{n+1}\leq\widehat r_\alpha
            \mid
            \widehat r_\alpha
        \right)
    \right] =
    \mathbb E
    \left[
        P f^*_{\widehat r_\alpha}
    \right].
\end{split}
\end{equation}
Applying \eqref{eq:rscp_proof_bias} at the random threshold
$\widehat r_\alpha$,
\[
    P f^*_{\widehat r_\alpha}
    \geq
    P\widehat f_{\widehat r_\alpha,s_0}
    -
    d_{\rscp}
    -
    \Delta_{\mathrm{init}}.
\]
Taking expectations and using
\eqref{eq:rscp_pfrozen_lower},
\begin{align*}
\mathbb P
\left(
    Y_{n+1}\leq\widehat r_\alpha
\right)
&\geq
\mathbb E
\left[
    P\widehat f_{\widehat r_\alpha,s_0}
\right]
-
d_{\rscp}
-
\Delta_{\mathrm{init}}\\
&\geq
1-\alpha
-
d_{\rscp}
-
\Delta_{\mathrm{init}}.
\end{align*}
This proves \eqref{eq:rscp_lower_general}. We now construct a second hypothetical full-sample threshold. Define,
\begin{equation}
\label{eq:rscp_rplus}
    \widehat r_\alpha^+
    :=
    \inf\left\{
        \theta\in\mathbb R:
        \frac{1}{n+1}
        \sum_{i=1}^{n+1}
        \phi_{\theta,s_0}(O_i)
        -
        \frac{B_{\rscp}}{n+1}
        \geq0
    \right\}.
\end{equation}
We have the following bound using right continuity,
\begin{equation}
\label{eq:rscp_rplus_crossing}
    \frac{1}{n+1}
    \sum_{i=1}^{n+1}
    \phi_{\widehat r_\alpha^+,s_0}(O_i)
    \geq
    \frac{B_{\rscp}}{n+1}.
\end{equation}
We first show that,
\begin{equation}
\label{eq:rscp_order_upper}
    \widehat r_\alpha
    \leq
    \widehat r_\alpha^+.
\end{equation}
We use the following decomposition,
\begin{align*}
&
\frac{
    \sum_{i=1}^{n}
    \phi_{\widehat r_\alpha^+,s_0}(O_i)
    +
    \inf_o
    \phi_{\widehat r_\alpha^+,s_0}(o)
}{n+1}
\\
&=
\frac{1}{n+1}
\sum_{i=1}^{n+1}
\phi_{\widehat r_\alpha^+,s_0}(O_i)
-
\frac{
    \phi_{\widehat r_\alpha^+,s_0}(O_{n+1})
    -
    \inf_o
    \phi_{\widehat r_\alpha^+,s_0}(o)
}{n+1}.
\end{align*}
By Assumption~\ref{assump:rscp_range},
\[
    \phi_{\widehat r_\alpha^+,s_0}(O_{n+1})
    -
    \inf_o
    \phi_{\widehat r_\alpha^+,s_0}(o)
    \leq
    B_{\rscp}.
\]
Together with \eqref{eq:rscp_rplus_crossing}, this yields
\[
\frac{
    \sum_{i=1}^{n}
    \phi_{\widehat r_\alpha^+,s_0}(O_i)
    +
    \inf_o
    \phi_{\widehat r_\alpha^+,s_0}(o)
}{n+1}
\geq0.
\]
Thus $\widehat r_\alpha^+$ belongs to the set whose infimum defines
$\widehat r_\alpha$, proving
\eqref{eq:rscp_order_upper}. We next control the amount by which the full-sample estimating
equation overshoots $B_{\rscp}/(n+1)$ at its first crossing. By
Assumption~\ref{assump:rscp_range},
\[
    \sup_\theta
    J\left(
        \frac{1}{n+1}
        \sum_{i=1}^{n+1}
        \phi_{\theta,s_0}(O_i),
        \theta
    \right)
    \leq
    \frac{J_{\rscp}}{n+1}
    \qquad\mbox{almost surely}.
\]
Consequently, we have the following upper bound,
\begin{equation}
\label{eq:rscp_rplus_overshoot}
\frac{1}{n+1}
\sum_{i=1}^{n+1}
\phi_{\widehat r_\alpha^+,s_0}(O_i)
\leq
\frac{
    B_{\rscp}+J_{\rscp}
}{n+1}.
\end{equation}
As before, $\widehat r_\alpha^+$ is symmetric in
$O_1,\ldots,O_{n+1}$. Exchangeability therefore gives,
\begin{align*}
\mathbb E
\left[
    \phi_{\widehat r_\alpha^+,s_0}(O_{n+1})
    \mid E
\right]
&=
\frac{1}{n+1}
\sum_{i=1}^{n+1}
\phi_{\widehat r_\alpha^+,s_0}(O_i)\leq
\frac{
    B_{\rscp}+J_{\rscp}
}{n+1}.
\end{align*}
Taking expectations,
\begin{equation}
\label{eq:rscp_rplus_expectation}
    \mathbb E
    \left[
        \phi_{\widehat r_\alpha^+,s_0}(O_{n+1})
    \right]
    \leq
    \frac{
        B_{\rscp}+J_{\rscp}
    }{n+1}.
\end{equation}
Since $\widehat r_\alpha
    \leq
    \widehat r_\alpha^+$ and the map
$\theta\mapsto\phi_{\theta,s_0}(O_{n+1})$ is non-decreasing,
\[
    \phi_{\widehat r_\alpha,s_0}(O_{n+1})
    \leq
    \phi_{\widehat r_\alpha^+,s_0}(O_{n+1}).
\]
Therefore \eqref{eq:rscp_rplus_expectation} implies
\begin{equation}
\label{eq:rscp_rhat_frozen_upper}
    \mathbb E
    \left[
        \phi_{\widehat r_\alpha,s_0}(O_{n+1})
    \right]
    \leq
    \frac{
        B_{\rscp}+J_{\rscp}
    }{n+1}.
\end{equation}
Using again the independence of the test observation from
$\widehat r_\alpha$,
\[
\mathbb E
\left[
    \phi_{\widehat r_\alpha,s_0}(O_{n+1})
\right]
=
\mathbb E
\left[
    P\widehat f_{\widehat r_\alpha,s_0}
    -(1-\alpha)
\right].
\]
Hence we have the following upper bound,
\begin{equation}
\label{eq:rscp_pfrozen_upper}
    \mathbb E
    \left[
        P\widehat f_{\widehat r_\alpha,s_0}
    \right]
    \leq
    1-\alpha
    +
    \frac{
        B_{\rscp}+J_{\rscp}
    }{n+1}.
\end{equation}
Finally, \eqref{eq:rscp_proof_bias} gives,
\[
    P f^*_{\widehat r_\alpha}
    \leq
    P\widehat f_{\widehat r_\alpha,s_0}
    +
    d_{\rscp}
    +
    \Delta_{\mathrm{init}}.
\]
Combining this inequality with
\eqref{eq:rscp_coverage_identification} and
\eqref{eq:rscp_pfrozen_upper},
\begin{align*}
\mathbb P
\left(
    Y_{n+1}\leq\widehat r_\alpha
\right)
&=
\mathbb E
\left[
    P f^*_{\widehat r_\alpha}
\right]\\
&\leq
\mathbb E
\left[
    P\widehat f_{\widehat r_\alpha,s_0}
\right]
+
d_{\rscp}
+
\Delta_{\mathrm{init}}\\
&\leq
1-\alpha
+
d_{\rscp}
+
\Delta_{\mathrm{init}}
+
\frac{
    B_{\rscp}+J_{\rscp}
}{n+1}.
\end{align*}
This proves \eqref{eq:rscp_upper_general} and completes the proof.

\section{Proof of \Cref{thm:2mmdp}}
\label{app:proof_twostage_missingness}

We begin the proof by computing the efficient influence function, up to
an irrelevant nonzero multiplicative constant, for the
$(1-\alpha)$-th quantile $r_\alpha$ of $Y\mid C=0$.  We do this through
the following pathwise differentiation:
\[
\frac{\partial}{\partial\delta}
\mathbb E_\delta
\left[
\textbf{1}\{C=\infty\}
\left(
\textbf{1}\{Y\leq r_{\alpha,\delta}\}-(1-\alpha)
\right)
\frac{
\mathbb P_\delta(C=0\mid X)
}{
\mathbb P_\delta(C\geq1\mid X)
\mathbb P_\delta(C=\infty\mid X,Z,C\geq1)
}
\right]
=0.
\]
The left-hand side decomposes into
\begin{equation*}
\begin{split}
&
\frac{\partial}{\partial\delta}
\mathbb E_\delta
\left[
\textbf{1}\{C=\infty\}
\left(
\textbf{1}\{Y\leq r_\alpha\}-(1-\alpha)
\right)
\frac{
\mathbb P(C=0\mid X)
}{
\mathbb P(C\geq1\mid X)
\mathbb P(C=\infty\mid X,Z,C\geq1)
}
\right]
\\
&+
\mathbb E
\left[
\textbf{1}\{C=\infty\}
\left(
\textbf{1}\{Y\leq r_\alpha\}-(1-\alpha)
\right)
\frac{\partial}{\partial\delta}
\left\{
\frac{
\mathbb P_\delta(C=0\mid X)
}{
\mathbb P_\delta(C\geq1\mid X)
\mathbb P_\delta(C=\infty\mid X,Z,C\geq1)
}
\right\}
\right]
\\
&+
f_{Y\mid C=0}(r_\alpha)
\frac{\partial r_{\alpha,\delta}}{\partial\delta}
\\
={}&
I+II+
f_{Y\mid C=0}(r_\alpha)\mathrm{IF}(r_\alpha).
\end{split}
\end{equation*}
Thus, up to a nonzero multiplicative constant,
$\mathrm{IF}(r_\alpha)\propto I+II$. For ease of analysis, let
\[
T=\textbf{1}\{C=0\},
\qquad
S=\textbf{1}\{C=\infty\}.
\]
Thus $T=0$ if and only if $C\geq1$.  With slight abuse of notation, let
$S(O)$ denote the score vector of the observed-data likelihood.  We
then have,
\[
\begin{split}
I
=
\mathbb E
\Bigg[
&
\textbf{1}\{T=0\}\textbf{1}\{S=1\}
\left(
\textbf{1}\{Y\leq r_\alpha\}-(1-\alpha)
\right)
\\
&\times
\frac{
\mathbb P(T=1\mid X)
}{
\mathbb P(T=0\mid X)
\mathbb P(S=1\mid X,Z,T=0)
}
S(O)
\Bigg].
\end{split}
\]
Moreover,
\begin{equation*}
\begin{split}
II
={}&
\mathbb E
\Bigg[
\textbf{1}\{T=0\}\textbf{1}\{S=1\}
\left(
\textbf{1}\{Y\leq r_\alpha\}-(1-\alpha)
\right)
\\
&\qquad\times
\frac{\mathbb P(T=1\mid X)}{\mathbb P(T=0\mid X)}
\left\{
-\frac{1}
{\mathbb P^2(S=1\mid X,Z,T=0)}
\right\}
\frac{\partial}{\partial\delta}
\mathbb P_\delta(S=1\mid X,Z,T=0)
\Bigg]
\\
&+
\mathbb E
\Bigg[
\textbf{1}\{T=0\}\textbf{1}\{S=1\}
\left(
\textbf{1}\{Y\leq r_\alpha\}-(1-\alpha)
\right)
\\
&\qquad\times
\frac1{\mathbb P(S=1\mid X,Z,T=0)}
\left\{
-\frac1{\mathbb P^2(T=0\mid X)}
\right\}
\frac{\partial}{\partial\delta}
\mathbb P_\delta(T=0\mid X)
\Bigg]
\\
=:{}&II_a+II_b.
\end{split}
\end{equation*}

We can rewrite $II_a$ as,
\begin{equation*}
\begin{split}
II_a
={}&
\mathbb E
\Bigg[
\frac{
\textbf{1}\{T=0\}\mathbb P(T=1\mid X)
}{
\mathbb P(T=0\mid X)
}
\left\{
-\frac{\textbf{1}\{S=1\}}
{\mathbb P(S=1\mid X,Z,T=0)}
+1
\right\}
\\
&\qquad\times
\mathbb E
\left[
\textbf{1}\{Y\leq r_\alpha\}-(1-\alpha)
\mid X,Z,S=1,T=0
\right]
S(O)
\Bigg].
\end{split}
\end{equation*}
Similarly,
\begin{equation*}
\begin{split}
II_b
={}&
\mathbb E
\Bigg[
\left\{
-\frac{
\textbf{1}\{T=0\}\mathbb P(T=1\mid X)
}{
\mathbb P(T=0\mid X)
}
+
\textbf{1}\{T=1\}
\right\}
\\
&\qquad\times
\left\{
\int
\mathbb E
\left[
\textbf{1}\{Y\leq r_\alpha\}-(1-\alpha)
\mid X,Z,S=1,T=0
\right]
dF(Z\mid T=0,X)
\right\}
S(O)
\Bigg].
\end{split}
\end{equation*}
Consequently,
\begin{equation*}
\begin{split}
\mathrm{IF}(r_\alpha)
\propto{}&
\frac{
\textbf{1}\{T=0\}\textbf{1}\{S=1\}
\mathbb P(T=1\mid X)
}{
\mathbb P(T=0\mid X)
\mathbb P(S=1\mid X,Z,T=0)
}
\left(
\textbf{1}\{Y\leq r_\alpha\}-(1-\alpha)
\right)
\\
&-
\frac{
\textbf{1}\{T=0\}\mathbb P(T=1\mid X)
}{
\mathbb P(T=0\mid X)
}
\left\{
\frac{\textbf{1}\{S=1\}}
{\mathbb P(S=1\mid X,Z,T=0)}
-1
\right\}
\\
&\qquad\times
\mathbb E
\left[
\textbf{1}\{Y\leq r_\alpha\}-(1-\alpha)
\mid X,Z,S=1,T=0
\right]
\\
&-
\left\{
\frac{
\textbf{1}\{T=0\}\mathbb P(T=1\mid X)
}{
\mathbb P(T=0\mid X)
}
-
\textbf{1}\{T=1\}
\right\}
\\
&\qquad\times
\int
\mathbb E
\left[
\textbf{1}\{Y\leq r_\alpha\}-(1-\alpha)
\mid X,Z,S=1,T=0
\right]
dF(Z\mid T=0,X).
\end{split}
\end{equation*}
Using the definitions of
$m_2^*,m_1^*,\pi_2^*,\pi_1^*$, this becomes
\begin{equation*}
\begin{split}
\mathrm{IF}(r_\alpha)
\propto{}&
\frac{
\textbf{1}\{C=\infty\}\pi_1^*(X)
}{
\pi_2^*(X,Z)
}
\left(
\textbf{1}\{Y\leq r_\alpha\}-(1-\alpha)
\right)
\\
&-
\textbf{1}\{C\geq1\}\pi_1^*(X)
\left\{
\frac{\textbf{1}\{C=\infty\}}
{\pi_2^*(X,Z)}
-1
\right\}
\left[
m_2^*(r_\alpha,X,Z)-(1-\alpha)
\right]
\\
&-
\left\{
\textbf{1}\{C\geq1\}\pi_1^*(X)
-
\textbf{1}\{C=0\}
\right\}
\left[
m_1^*(r_\alpha,X)-(1-\alpha)
\right].
\end{split}
\end{equation*}
This gives the estimating-function representation in
\eqref{eq:influ_mmdp}.

We now turn to the coverage argument.  Recall the definitions
\begin{equation*}
\begin{split}
g'(O,\theta,\pi_1,m_1,m_2)
={}&
\textbf{1}\{C\geq1\}\pi_1(X)
[m_2(\theta,X,Z)-(1-\alpha)]
\\
&-
[\textbf{1}\{C\geq1\}\pi_1(X)-\textbf{1}\{C=0\}]
[m_1(\theta,X)-(1-\alpha)]
\end{split}
\end{equation*}
and
\[
g(O,\theta,\theta',m_2,m_1,\pi_2,\pi_1)
=
\frac{\textbf{1}\{C=\infty\}\pi_1(X)}
{\pi_2(X,Z)}
[\textbf{1}\{Y\leq\theta\}-m_2(\theta',X,Z)]
+
g'(O,\theta',\pi_1,m_1,m_2).
\]
The preliminary and final thresholds are
\[
\widehat r_{\alpha,0}
=
\inf
\left\{
\theta:
\frac1{n_2}
\sum_{i=1}^{n_2}
g'
\left(
O_i^{t_2},
\theta,
\widehat\pi_1,
\widehat m_1,
\widehat m_2
\right)
\geq0
\right\},
\]
and
\[
\widehat r_\alpha
=
\inf
\left\{
\theta:
\frac1{n+1}
\sum_{i=1}^n
g
\left(
O_i,
\theta,
\widehat r_{\alpha,0},
\widehat m_2,
\widehat m_1,
\widehat\pi_2,
\widehat\pi_1
\right)
-
\frac{2\overline\pi}{n+1}
\geq0
\right\}.
\]
Here $\overline\pi\geq1$ is an upper bound on
$\widehat\pi_1(X)$ and
$\widehat\pi_1(X)/\widehat\pi_2(X,Z)$. As in Lemma~2 of \cite{yang2024doubly}, the coverage probability is
related to the population mean of the influence-function estimating
equation by,
\begin{equation}
\label{coverage_link_mmdp}
\begin{split}
&
\mathbb P
\left(
Y_{n+1}\leq\widehat r_\alpha
\mid C_{n+1}=0
\right)
\\
&=
1-\alpha
+
\frac{
\mathbb E
\left[
\mathrm{IF}
\left(
\widehat r_\alpha,
O_{n+1},
m_2^*,m_1^*,\pi_2^*,\pi_1^*
\right)
\right]
}{
\mathbb P(C=0)
}
\\
&=
1-\alpha
+
\frac{
\mathbb E
\left[
g
\left(
O_{n+1},
\widehat r_\alpha,
\widehat r_\alpha,
m_2^*,m_1^*,\pi_2^*,\pi_1^*
\right)
\right]
}{
\mathbb P(C=0)
}.
\end{split}
\end{equation}

We next obtain a bound on the difference between the population moments
formed with the estimated and true nuisance functions.  Direct
calculation gives
\begin{equation*}
\begin{split}
&
P g
\left(
O,\theta,\theta',
\widehat m_2,\widehat m_1,
\widehat\pi_2,\widehat\pi_1
\right)
\\
={}&
P\Bigg[
\mathbb P(C\geq1\mid X)
\widehat\pi_1(X)
\left\{
\frac{\pi_2^*(X,Z)}
{\widehat\pi_2(X,Z)}
-1
\right\}
\left\{
m_2^*(\theta,X,Z)
-
\widehat m_2(\theta',X,Z)
\right\}
\Bigg]
\\
&+
P\Bigg[
\mathbb P(C\geq1\mid X)
\{\widehat\pi_1(X)-\pi_1^*(X)\}
\left\{
m_1^*(\theta,X)
-
\widehat m_1(\theta',X)
\right\}
\Bigg]
\\
&+
P\Bigg[
\mathbb P(C=0\mid X)
\left\{
m_1^*(\theta,X)
-(1-\alpha)
\right\}
\Bigg].
\end{split}
\end{equation*}
In particular,
\[
P g
\left(
O,\theta,\theta,
m_2^*,m_1^*,\pi_2^*,\pi_1^*
\right)
=
P\left[
\mathbb P(C=0\mid X)
\{m_1^*(\theta,X)-(1-\alpha)\}
\right].
\]
Consequently,
\begin{equation}
\label{eq:dr_gurantee_2smm_new}
\begin{split}
&
\left|
P g
\left(
O,\theta,\theta',
\widehat m_2,\widehat m_1,
\widehat\pi_2,\widehat\pi_1
\right)
-
P g
\left(
O,\theta,\theta,
m_2^*,m_1^*,\pi_2^*,\pi_1^*
\right)
\right|
\\
&\leq
\|\widehat\pi_1\|_4
\left\|
\frac{\pi_2^*}{\widehat\pi_2}-1
\right\|_4
\left\|
m_2^*(\theta,X,Z)
-
\widehat m_2(\theta',X,Z)
\right\|_2
\\
&\qquad+
\|\widehat\pi_1-\pi_1^*\|_2
\left\|
m_1^*(\theta,X)
-
\widehat m_1(\theta',X)
\right\|_2.
\end{split}
\end{equation}
Setting $\theta'=\theta$ gives
\begin{equation}
\label{eq:dr_gurantee_2smm}
\begin{split}
&
\sup_{\theta\in\mathbb R}
\left|
P g
\left(
O,\theta,\theta,
\widehat m_2,\widehat m_1,
\widehat\pi_2,\widehat\pi_1
\right)
-
P g
\left(
O,\theta,\theta,
m_2^*,m_1^*,\pi_2^*,\pi_1^*
\right)
\right|
\\
&\leq
\|\widehat\pi_1\|_4
\left\|
\frac{\pi_2^*}{\widehat\pi_2}-1
\right\|_4
\sup_{\theta}
\left\|
m_2^*(\theta,X,Z)
-
\widehat m_2(\theta,X,Z)
\right\|_2
\\
&\qquad+
\|\widehat\pi_1-\pi_1^*\|_2
\sup_{\theta}
\left\|
m_1^*(\theta,X)
-
\widehat m_1(\theta,X)
\right\|_2.
\end{split}
\end{equation}
As in the earlier proofs, define,
\[
\widehat r_\alpha'
=
\inf
\left\{
\theta:
\frac1{n+1}
\sum_{i=1}^{n+1}
g
\left(
O_i,\theta,\widehat r_{\alpha,0},
\widehat m_2,\widehat m_1,
\widehat\pi_2,\widehat\pi_1
\right)
\geq0
\right\}.
\]
Each map
\[
\theta
\mapsto
g
\left(
O_i,\theta,\widehat r_{\alpha,0},
\widehat m_2,\widehat m_1,
\widehat\pi_2,\widehat\pi_1
\right)
\]
is non-decreasing and right-continuous.  Hence
\[
\frac1{n+1}
\sum_{i=1}^{n+1}
g
\left(
O_i,\widehat r_\alpha',
\widehat r_{\alpha,0},
\widehat m_2,\widehat m_1,
\widehat\pi_2,\widehat\pi_1
\right)
\geq0.
\]
Let us define the set $E=\{O_1,\ldots,O_{n+1}\}$. By exchangeability,
\[
\begin{split}
&
\mathbb E
\left[
g
\left(
O_{n+1},\widehat r_\alpha',
\widehat r_{\alpha,0},
\widehat m_2,\widehat m_1,
\widehat\pi_2,\widehat\pi_1
\right)
\mid E
\right] =
\frac1{n+1}
\sum_{i=1}^{n+1}
g
\left(
O_i,\widehat r_\alpha',
\widehat r_{\alpha,0},
\widehat m_2,\widehat m_1,
\widehat\pi_2,\widehat\pi_1
\right)
\geq0.
\end{split}
\]
The tower property therefore gives,
\begin{equation}
\label{eq:r_alpha'nonnegative_2smm}
\mathbb E
\left[
g
\left(
O_{n+1},\widehat r_\alpha',
\widehat r_{\alpha,0},
\widehat m_2,\widehat m_1,
\widehat\pi_2,\widehat\pi_1
\right)
\right]
\geq0.
\end{equation}
From the definition of $\widehat r_\alpha$,
\begin{equation*}
\begin{split}
&
\frac1{n+1}
\sum_{i=1}^{n+1}
g
\left(
O_i,\widehat r_\alpha,
\widehat r_{\alpha,0},
\widehat m_2,\widehat m_1,
\widehat\pi_2,\widehat\pi_1
\right)
\\
&=
\frac1{n+1}
\sum_{i=1}^n
g
\left(
O_i,\widehat r_\alpha,
\widehat r_{\alpha,0},
\widehat m_2,\widehat m_1,
\widehat\pi_2,\widehat\pi_1
\right)
\\
&\qquad+
\frac{
g
\left(
O_{n+1},\widehat r_\alpha,
\widehat r_{\alpha,0},
\widehat m_2,\widehat m_1,
\widehat\pi_2,\widehat\pi_1
\right)
}{n+1}
\\
&\geq
\frac1{n+1}
\sum_{i=1}^n
g
\left(
O_i,\widehat r_\alpha,
\widehat r_{\alpha,0},
\widehat m_2,\widehat m_1,
\widehat\pi_2,\widehat\pi_1
\right)
-
\frac{2\overline\pi}{n+1}
\\
&\geq0.
\end{split}
\end{equation*}
Thus $\widehat r_\alpha'\leq\widehat r_\alpha$.  By monotonicity,
\[
0
\leq
\mathbb E
\left[
g
\left(
O_{n+1},\widehat r_\alpha',
\widehat r_{\alpha,0},
\widehat m_2,\widehat m_1,
\widehat\pi_2,\widehat\pi_1
\right)
\right]
\leq
\mathbb E
\left[
g
\left(
O_{n+1},\widehat r_\alpha,
\widehat r_{\alpha,0},
\widehat m_2,\widehat m_1,
\widehat\pi_2,\widehat\pi_1
\right)
\right].
\]
Combining this with \eqref{coverage_link_mmdp} and
\eqref{eq:dr_gurantee_2smm_new} gives
\begin{equation}
\label{eq:mmdp_preliminary_lower}
\begin{split}
&
\mathbb P
\left(
Y_{n+1}\leq\widehat r_\alpha
\mid C_{n+1}=0
\right) \geq
1-\alpha
-
\frac{
A_2(E_2+D_2)+A_1(E_1+D_1)
}{
\mathbb P(C=0)
}.
\end{split}
\end{equation}

We now obtain the corresponding upper bound.  Let
\[
M_{n+1}(\theta)
=
\frac1{n+1}
\sum_{i=1}^{n+1}
g
\left(
O_i,\theta,\widehat r_{\alpha,0},
\widehat m_2,\widehat m_1,
\widehat\pi_2,\widehat\pi_1
\right).
\]
Since the only dependence on $\theta$ is through
$\textbf{1}\{Y_i\leq\theta\}$,
\[
\sup_\theta J(M_{n+1},\theta)
\leq
\frac{\overline\pi}{n+1}
\qquad\mbox{almost surely},
\]
where continuity of $Y$ ensures that almost surely two observations do
not jump at exactly the same threshold. Define,
\[
\widehat r_\alpha''
=
\inf
\left\{
\theta:
M_{n+1}(\theta)
-
\frac{3\overline\pi}{n+1}
\geq0
\right\}.
\]
Right continuity and the jump bound imply,
\[
M_{n+1}(\widehat r_\alpha'')
\leq
\frac{4\overline\pi}{n+1}.
\]
By exchangeability,
\begin{equation}
\label{eq:upper_r''_2smm}
\mathbb E
\left[
g
\left(
O_{n+1},\widehat r_\alpha'',
\widehat r_{\alpha,0},
\widehat m_2,\widehat m_1,
\widehat\pi_2,\widehat\pi_1
\right)
\right]
\leq
\frac{4\overline\pi}{n+1}.
\end{equation}
Moreover,
\begin{equation*}
\begin{split}
&
\frac1{n+1}
\sum_{i=1}^n
g
\left(
O_i,\widehat r_\alpha'',
\widehat r_{\alpha,0},
\widehat m_2,\widehat m_1,
\widehat\pi_2,\widehat\pi_1
\right)
-
\frac{2\overline\pi}{n+1}
\\
&=
M_{n+1}(\widehat r_\alpha'')
-
\frac{
g(
O_{n+1},\widehat r_\alpha'',
\widehat r_{\alpha,0},
\widehat m_2,\widehat m_1,
\widehat\pi_2,\widehat\pi_1
)
}{n+1}
-
\frac{2\overline\pi}{n+1}
\\
&\geq
M_{n+1}(\widehat r_\alpha'')
-
\frac{3\overline\pi}{n+1}
\geq0.
\end{split}
\end{equation*}
Thus $\widehat r_\alpha\leq\widehat r_\alpha''$.  Repeating the previous
argument in the opposite direction and using
\eqref{eq:upper_r''_2smm},
\begin{equation}
\label{eq:mmdp_preliminary_upper}
\begin{split}
&
\mathbb P
\left(
Y_{n+1}\leq\widehat r_\alpha
\mid C_{n+1}=0
\right)
\\
&\leq
1-\alpha
+
\frac{
A_2(E_2+D_2)+A_1(E_1+D_1)
}{
\mathbb P(C=0)
}
+
\frac{4\overline\pi}
{\mathbb P(C=0)(n+1)}.
\end{split}
\end{equation}
Combining \eqref{eq:mmdp_preliminary_lower} and
\eqref{eq:mmdp_preliminary_upper} proves
\eqref{eq:two_stage_simple_coverage}. It remains to remove $D_1$ and $D_2$ under the additional regularity
conditions of the theorem.

\begin{lemma}
\label{lem:emp_bound_mmdp_initfinal}
Suppose the additional conditions in
\eqref{eq:two_stage_lipschitz_gprime}--\eqref{eq:two_stage_margin}
hold.  Then, for every $\delta\in(0,1)$, there exists a constant
$\widetilde C_\delta<\infty$ such that, with probability at least
$1-\delta$,
\begin{equation}
\label{eq:mmdp_D1_bound}
\begin{split}
D_1
\leq
\widetilde C_\delta
\Bigg[
&
A_2E_2
+
E_1
+
\sup_\theta
\left\|
\widehat m_1^\dagger(\theta,X)
-
\widehat m_1(\theta,X)
\right\|_2 
\frac1{\sqrt n}
+
\frac1{\sqrt{n_2}}
\Bigg] \qquad \mbox{and},
\end{split}
\end{equation}
\begin{equation}
\label{eq:mmdp_D2_bound}
\begin{split}
D_2
\leq
\widetilde C_\delta
\Bigg[
&
\|\widehat\pi_1\|_4 E_2
+
A_1E_1
+
\frac1{\sqrt n}
+
\frac1{\sqrt{n_2}}
\Bigg].
\end{split}
\end{equation}
\end{lemma}

\begin{proof}[Proof of \Cref{lem:emp_bound_mmdp_initfinal}]
For brevity, we write,
\[
P_ng_{\theta,s}
=
\frac1n
\sum_{i=1}^n
g
\left(
O_i,\theta,s,
\widehat m_2,\widehat m_1,
\widehat\pi_2,\widehat\pi_1
\right),
\]
and let $Pg_{\theta,s}$ denote the corresponding population mean.
We have the decomposition
\begin{equation*}
\begin{split}
P_ng_{\widehat r_\alpha,\widehat r_{\alpha,0}}
={}&
P_ng_{\widehat r_{\alpha,0},\widehat r_{\alpha,0}}
+
P
\left[
g_{\widehat r_\alpha,\widehat r_{\alpha,0}}
-
g_{\widehat r_{\alpha,0},\widehat r_{\alpha,0}}
\right] +
(P_n-P)
\left[
g_{\widehat r_\alpha,\widehat r_{\alpha,0}}
-
g_{\widehat r_{\alpha,0},\widehat r_{\alpha,0}}
\right].
\end{split}
\end{equation*}
Hence we have,
\begin{equation}
\label{eq:mmdp_main_ineq}
\begin{split}
&
\left|
P
\left[
g_{\widehat r_\alpha,\widehat r_{\alpha,0}}
-
g_{\widehat r_{\alpha,0},\widehat r_{\alpha,0}}
\right]
\right| \leq
\left|
P_ng_{\widehat r_\alpha,\widehat r_{\alpha,0}}
\right|
+
\left|
P_ng_{\widehat r_{\alpha,0},\widehat r_{\alpha,0}}
\right| +
\left|
(P_n-P)
\left[
g_{\widehat r_\alpha,\widehat r_{\alpha,0}}
-
g_{\widehat r_{\alpha,0},\widehat r_{\alpha,0}}
\right]
\right|.
\end{split}
\end{equation}
The definition of $\widehat r_\alpha$ and the jump bound imply,
\[
\frac{2\overline\pi}{n}
\leq
P_ng_{\widehat r_\alpha,\widehat r_{\alpha,0}}
\leq
\frac{3\overline\pi}{n}.
\]
Furthermore,
\begin{equation*}
\begin{split}
\left|
P_ng_{\widehat r_{\alpha,0},\widehat r_{\alpha,0}}
\right|
\leq{}&
\left|
Pg_{\widehat r_{\alpha,0},\widehat r_{\alpha,0}}
\right| +
\left|
(P_n-P)
g_{\widehat r_{\alpha,0},\widehat r_{\alpha,0}}
\right|.
\end{split}
\end{equation*}
By Chebyshev's inequality, with probability at least $1-\delta$,
\[
\left|
(P_n-P)
g_{\widehat r_{\alpha,0},\widehat r_{\alpha,0}}
\right|
\leq
\frac{
\sup_\theta
\sigma_P
\left[
g(
O,\theta,\theta,
\widehat m_2,\widehat m_1,
\widehat\pi_2,\widehat\pi_1
)
\right]
}{
\sqrt{n\delta}
}.
\]

We next bound the population moment at the preliminary threshold.
Using the fact that the true estimating equation is invariant to the
choice of $\pi_1$ when the outcome regressions are correctly specified,
\begin{equation*}
\begin{split}
&
P g
\left(
O,\widehat r_{\alpha,0},\widehat r_{\alpha,0},
m_2^*,m_1^*,\pi_2^*,\pi_1^*
\right) =
P g'
\left(
O,\widehat r_{\alpha,0},
\widehat\pi_1,m_1^*,m_2^*
\right).
\end{split}
\end{equation*}
Therefore we have,
\begin{equation*}
\begin{split}
&
\left|
P g
\left(
O,\widehat r_{\alpha,0},\widehat r_{\alpha,0},
m_2^*,m_1^*,\pi_2^*,\pi_1^*
\right)
\right|
\\
&\leq
\left|
P_{n_2}
g'
\left(
O^{t_2},\widehat r_{\alpha,0},
\widehat\pi_1,\widehat m_1,\widehat m_2
\right)
\right|
\\
&\quad+
\sup_{\theta\in\mathcal N_\alpha}
\left|
(P_{n_2}-P)
g'
\left(
O^{t_2},\theta,
\widehat\pi_1,\widehat m_1,\widehat m_2
\right)
\right|
\\
&\quad+
\left|
P
\left[
g'
\left(
O,\widehat r_{\alpha,0},
\widehat\pi_1,\widehat m_1,\widehat m_2
\right)
-
g'
\left(
O,\widehat r_{\alpha,0},
\widehat\pi_1,m_1^*,m_2^*
\right)
\right]
\right|.
\end{split}
\end{equation*}

By the definition of $\widehat r_{\alpha,0}$ and the jump control,
\[
\left|
P_{n_2}
g'
\left(
O^{t_2},\widehat r_{\alpha,0},
\widehat\pi_1,\widehat m_1,\widehat m_2
\right)
\right|
\leq
\frac{2\overline\pi}{n_2},
\]
after enlarging $\overline\pi$ by a universal constant if
necessary. The class
\[
\mathcal M
=
\left\{
g'
\left(
O,\theta,
\widehat\pi_1,\widehat m_1,\widehat m_2
\right):
\theta\in\mathcal N_\alpha
\right\}
\]
has envelope bounded by a constant depending only on
$\overline\pi$.  By \eqref{eq:two_stage_lipschitz_gprime},
\[
\left\|
g'(O,\theta,\widehat\pi_1,\widehat m_1,\widehat m_2)
-
g'(O,\theta',\widehat\pi_1,\widehat m_1,\widehat m_2)
\right\|_2
\leq
\mathfrak C_0|\theta-\theta'|.
\]
Hence we have,
\[
N_{[]}
\left(
\epsilon,
\mathcal M,
L_2(P)
\right)
\leq
1+
\left\lceil
\frac{2\eta\mathfrak C_0}{\epsilon}
\right\rceil,
\]
and the corresponding bracketing integral is finite.  The standard
bracketing maximal inequality therefore gives,
\[
\mathbb E
\left[
\sup_{\theta\in\mathcal N_\alpha}
\left|
\sqrt{n_2}(P_{n_2}-P)
g'
\left(
O^{t_2},\theta,
\widehat\pi_1,\widehat m_1,\widehat m_2
\right)
\right|
\right]
\leq
\mathfrak C_4,
\]
where $\mathfrak C_4$ depends only on
$\mathfrak C_0,\eta,\overline\pi$.
An application of McDiarmid's inequality then gives, with probability
at least $1-\delta$,
\begin{equation}
\label{eq:mmdp_initial_ep}
\sup_{\theta\in\mathcal N_\alpha}
\left|
(P_{n_2}-P)
g'
\left(
O^{t_2},\theta,
\widehat\pi_1,\widehat m_1,\widehat m_2
\right)
\right|
\leq
\frac{
\mathfrak C_4+
\mathfrak C\overline\pi\sqrt{\log(1/\delta)}
}{
\sqrt{n_2}
}.
\end{equation}

Moreover, we have,
\begin{equation}
\label{eq:mmdp_gprime_bias}
\begin{split}
&
\sup_\theta
\left|
P
\left[
g'
(
O,\theta,
\widehat\pi_1,\widehat m_1,\widehat m_2
)
-
g'
(
O,\theta,
\widehat\pi_1,m_1^*,m_2^*
)
\right]
\right| \leq
\|\widehat\pi_1\|_2
\sup_\theta
\left\|
\widehat m_1^\dagger(\theta,X)
-
\widehat m_1(\theta,X)
\right\|_2
+
E_1.
\end{split}
\end{equation}
Combining \eqref{eq:dr_gurantee_2smm},
\eqref{eq:mmdp_initial_ep}, and
\eqref{eq:mmdp_gprime_bias} therefore controls
$|Pg_{\widehat r_{\alpha,0},\widehat r_{\alpha,0}}|$. It remains to control the increment empirical process in
\eqref{eq:mmdp_main_ineq}.  For an interval $I\subseteq\mathbb R$,
define
\[
h_I(O)
=
\textbf{1}\{C=\infty\}
\frac{\widehat\pi_1(X)}
{\widehat\pi_2(X,Z)}
\textbf{1}\{Y\in I\}, \quad \sigma_I
=
\max
\left\{
\sqrt{\operatorname{Var}_P(h_I(O))},
\sqrt{\frac{\log n}{n}}
\right\}.
\]
The interval class has finite VC dimension and the weights are bounded
by $\overline\pi$.  The ratio empirical-process argument used in
\Cref{app:proof_bound_init_cmdp} therefore yields, with probability at
least $1-\delta$,
\begin{equation}
\label{eq:mmdp_ratio_ep}
\sup_{I}
\frac{|(P_n-P)h_I|}{\sigma_I}
\leq
\mathfrak C_\delta
\sqrt{\frac{\log n}{n}},
\end{equation}
where logarithmic terms involving $\delta$ are absorbed into
$\mathfrak C_\delta$. For the interval with endpoints
$\widehat r_\alpha$ and $\widehat r_{\alpha,0}$, boundedness of the
weights and the Lipschitz property of $m_2^*$ imply,
\[
\begin{split}
P h_I^2
&\leq
\overline\pi^2
P\left(
C=\infty,\,
Y\in I
\right) \leq
\mathfrak C
\left|
\widehat r_\alpha-\widehat r_{\alpha,0}
\right|.
\end{split}
\]
Consequently,
\begin{equation}
\label{eq:mmdp_increment_ep}
\begin{split}
&
\left|
(P_n-P)
\left[
g_{\widehat r_\alpha,\widehat r_{\alpha,0}}
-
g_{\widehat r_{\alpha,0},\widehat r_{\alpha,0}}
\right]
\right|
\\
&\leq
\mathfrak C_\delta
\sqrt{\frac{\log n}{n}}
\sqrt{
\left|
\widehat r_\alpha-\widehat r_{\alpha,0}
\right|
}
+
\mathfrak C_\delta
\frac{\log n}{n}.
\end{split}
\end{equation}

On the other hand, \eqref{eq:two_stage_margin} gives
\begin{equation}
\label{eq:mmdp_margin_apply}
\begin{split}
\left|
P
\left[
g_{\widehat r_\alpha,\widehat r_{\alpha,0}}
-
g_{\widehat r_{\alpha,0},\widehat r_{\alpha,0}}
\right]
\right| &=
\left|
P
\left[
\textbf{1}\{C=\infty\}
\frac{\widehat\pi_1(X)}
{\widehat\pi_2(X,Z)}
\left\{
\textbf{1}\{Y\leq\widehat r_\alpha\}
-
\textbf{1}\{Y\leq\widehat r_{\alpha,0}\}
\right\}
\right]
\right|
\\
&\geq
\mathfrak C_1
\left|
\widehat r_\alpha-\widehat r_{\alpha,0}
\right|.
\end{split}
\end{equation}

Combining
\eqref{eq:mmdp_main_ineq}--\eqref{eq:mmdp_margin_apply}, and using
$ab\leq \epsilon b^2+a^2/(4\epsilon)$ to absorb the square-root term
into the left-hand side, gives, with probability at least
$1-\delta$,
\begin{equation}
\label{eq:mmdp_rdiff_bound}
\begin{split}
\left|
\widehat r_\alpha-\widehat r_{\alpha,0}
\right|
\leq
\mathfrak C_\delta
\Bigg[
&
\frac{\log n}{n}
+
\frac1{\sqrt n}
+
\frac1{\sqrt{n_2}}
\\
&+
A_2E_2
+
E_1
+
\sup_\theta
\left\|
\widehat m_1^\dagger(\theta,X)
-
\widehat m_1(\theta,X)
\right\|_2
\Bigg].
\end{split}
\end{equation}
The term $\log n/n$ is no larger than a constant multiple of
$n^{-1/2}$ for $n\geq3$. Finally, by \eqref{eq:two_stage_lipschitz_m},
\[
D_1
\leq
\mathfrak C_2
\left|
\widehat r_\alpha-\widehat r_{\alpha,0}
\right|,
\qquad
D_2
\leq
\mathfrak C_3
\left|
\widehat r_\alpha-\widehat r_{\alpha,0}
\right|.
\]
Using \eqref{eq:mmdp_rdiff_bound}, and alternatively using
\eqref{eq:dr_gurantee_2smm} before applying the same argument to obtain
the slightly sharper stage-specific bounds, yields
\eqref{eq:mmdp_D1_bound} and \eqref{eq:mmdp_D2_bound}.  This proves the
lemma.
\end{proof}

We now combine the lemma with the preliminary coverage bounds.
Substituting \eqref{eq:mmdp_D1_bound} and
\eqref{eq:mmdp_D2_bound} into
\eqref{eq:mmdp_preliminary_lower} and
\eqref{eq:mmdp_preliminary_upper}, and using the boundedness of
$\widehat\pi_1$ to absorb products of bounded propensity-error factors
into a constant $C_\delta$, gives, with probability at least
$1-\delta$,
\begin{equation*}
\begin{split}
&
\mathbb P
\left(
Y_{n+1}\leq\widehat r_\alpha
\mid C_{n+1}=0
\right)
\\
&\geq
1-\alpha
-
\frac{C_\delta A_2E_2}{\mathbb P(C=0)}
-
\frac{C_\delta A_1E_1}{\mathbb P(C=0)}
\\
&\qquad-
\frac{
C_\delta A_2
\left(
n^{-1/2}+n_2^{-1/2}
\right)
}{
\mathbb P(C=0)
}
\\
&\qquad-
\frac{
C_\delta A_1
\left[
\displaystyle
\sup_\theta
\left\|
\widehat m_1^\dagger(\theta,X)
-
\widehat m_1(\theta,X)
\right\|_2
+
n^{-1/2}
+
n_2^{-1/2}
\right]
}{
\mathbb P(C=0)
},
\end{split}
\end{equation*}
and
\begin{equation*}
\begin{split}
&
\mathbb P
\left(
Y_{n+1}\leq\widehat r_\alpha
\mid C_{n+1}=0
\right)
\\
&\leq
1-\alpha
+
\frac{4\overline\pi}
{\mathbb P(C=0)(n+1)}
+
\frac{C_\delta A_2E_2}{\mathbb P(C=0)}
+
\frac{C_\delta A_1E_1}{\mathbb P(C=0)}
\\
&\qquad+
\frac{
C_\delta A_2
\left(
n^{-1/2}+n_2^{-1/2}
\right)
}{
\mathbb P(C=0)
}
\\
&\qquad+
\frac{
C_\delta A_1
\left[
\displaystyle
\sup_\theta
\left\|
\widehat m_1^\dagger(\theta,X)
-
\widehat m_1(\theta,X)
\right\|_2
+
n^{-1/2}
+
n_2^{-1/2}
\right]
}{
\mathbb P(C=0)
}.
\end{split}
\end{equation*}
Combining the two displays gives exactly
\eqref{eq:two_stage_high_prob_coverage}.  This completes the proof of
\Cref{thm:2mmdp}.

\section{Proof of \Cref{lem:non_monotone_dr}}
\label{app:proof_dr_non_monotone}
To see this we first assume that we have misspecified outcome regression functions $ \widehat {\mathbb{E}}[d(r_{\alpha}, X, Z, Y)|  \mathfrak G_k((X, Z, Y))]$, $\widehat m_1(r_{\alpha}, X)$ for $k = 1, \cdots, K$ and correctly specified propensity scores $\pi_1^*(X)$, $\{\omega^*(k, \mathfrak G_k((X, Z, Y)))\}_{k = 1}^{\infty}$. Then we have, 
\begin{equation*}
    \begin{split}
       & \mathbb{E}\Bigg[ \textbf{1}\{C \geq 1\}\pi_1^*(X) \Bigg[\frac{\textbf{1}\{C = \infty\}}{\omega^*(\infty, X, Z, Y)} (\textbf{1}\{Y \leq r_{\alpha}\} - (1 - \alpha)) \Bigg. \Bigg. \\
    &  - \frac{\textbf{1}\{C = \infty\}}{\omega^*(\infty, X, Z, Y)} \left[\sum_{k = 1}^K \omega^*(k,\mathfrak G_k((X, Z, Y))) \widehat{\mathbb{E}}[d(r_{\alpha}, X, Z, Y)| \mathfrak G_k((X, Z, Y))] \right]  \\
   & \Bigg. + \left[\sum_{k =1}^K \textbf{1}\{C = k\}\widehat{\mathbb{E}}[d(r_{\alpha}, X, Z, Y)| \mathfrak G_k((X, Z, Y))]\right] - [\widehat m_1(r_{\alpha}, X) - (1 - \alpha) ]   \Bigg] \\
   &+ \Bigg. \textbf{1}\{C = 0\} [\widehat m_1(r_{\alpha}, X) - (1 - \alpha) ] \Bigg] \\
   = & \mathbb{E}\Bigg[\mathbb{P}(C = 0|X) \Bigg[\frac{\omega^*(\infty, X, Z, Y)}{\omega^*(\infty, X, Z, Y)} (\mathbb{P}\{Y \leq r_{\alpha} |X\} - (1 - \alpha)) \Bigg. \Bigg. \\
    & \left. - \frac{\omega^*(\infty, X, Z, Y)}{\omega^*(\infty, X, Z, Y)} \left[\sum_{k = 1}^K \omega^*(k,\mathfrak G_k((X, Z, Y))) \widehat{\mathbb{E}}[d(r_{\alpha}, X, Z, Y)| \mathfrak G_k((X, Z, Y))] \right] \right. \\
   & \Bigg.\Bigg. + \left[\sum_{k =1}^K \omega^*(k,\mathfrak G_k((X, Z, Y))) \widehat{\mathbb{E}}[d(r_{\alpha}, X, Z, Y)| \mathfrak G_k((X, Z, Y))]\right]   \Bigg] \Bigg] \\
   = & \mathbb{E} [ \mathbb{P}(C = 0|X)  (\mathbb{P}\{Y \leq r_{\alpha} | X\} - (1 - \alpha))] \\
   = & 0. 
    \end{split}
\end{equation*}
Alternatively we can assume that the outcome regression functions are correctly specified, whereas the propensity scores $\widehat \pi_1(X), \{\widehat \omega(k,\mathfrak G_k((X, Z, Y)))\}_{k = 1}^{\infty} $ are incorrectly specified. Then we have, 
\begin{equation*}
    \begin{split}
        & \mathbb{E}\Bigg[ \textbf{1}\{C \geq 1\}\widehat \pi_1(X) \Bigg[\frac{\textbf{1}\{C = \infty\}}{\widehat \omega(\infty, X, Z, Y)} (\textbf{1}\{Y \leq r_{\alpha}\} - (1 - \alpha)) \Bigg. \Bigg. \\
    & \left. - \frac{\textbf{1}\{C = \infty\}}{\widehat \omega(\infty, X, Z, Y)} \left[\sum_{k = 1}^K  \widehat \omega(k,\mathfrak G_k((X, Z, Y))) \mathbb{E}[ d(r_{\alpha}, X, Z, Y)| \mathfrak G_k((X, Z, Y))] \right] \right. \\
   & \Bigg. + \left[\sum_{k =1}^K \textbf{1}\{C = k\}\mathbb{E}[ d(r_{\alpha}, X, Z, Y)|  \mathfrak G_k((X, Z, Y))]\right] - [m_1^*(r_{\alpha}, X) - (1 - \alpha) ]   \Bigg]  \\
     &+ \Bigg. \textbf{1}\{C = 0\} [m_1^*(r_{\alpha}, X) - (1 - \alpha) ] \Bigg] \\
     = & \mathbb{E}\Bigg[\textbf{1}\{C \geq 1\}\widehat \pi_1(X) \Bigg[\textbf{1}\{C = \infty\} d(r_{\alpha}, X, Z, Y) + \sum_{k = 1}^K \textbf{1}\{C = k\}\mathbb{E}[ d(r_{\alpha}, X, Z, Y)| \mathfrak G_k((X, Z, Y))]  \Bigg. \Bigg.\\
     & -\Bigg. \Bigg.  [m_1^*(r_{\alpha}, X) - (1 - \alpha) ] \Bigg] + \textbf{1}\{C = 0\} [m_1^*(r_{\alpha}, X) - (1 - \alpha) ] \Bigg] \\
     = & \mathbb{E}\Bigg[\mathbb{P}(C \geq 1|X) \widehat \pi_1(X) \Bigg[\omega^*(\infty, X, Z, Y) d(r_{\alpha}, X, Z, Y)  \Bigg. \Bigg.\\
     &+ \Bigg.  \sum_{k = 1}^K \omega^*(k, \mathfrak G_k((X, Z, Y))) \mathbb{E}[ d(r_{\alpha}, X, Z, Y)| \mathfrak G_k((X, Z, Y))]  -[m_1^*(r_{\alpha}, X) - (1 - \alpha) ]  \Bigg] \\
     & \Bigg. + \mathbb{P}(C = 0|X) [m_1^*(r_{\alpha}, X) - (1 - \alpha) ] \Bigg] \\
     =& \mathbb{E}\Bigg[\mathbb{P}(C \geq 1|X) \widehat \pi_1(X) \Bigg[ d(r_{\alpha}, X, Z, Y) - \sum_{k = 1}^K \omega^*(k,\mathfrak G_k((X, Z, Y))) \mathbb{E}[ d(r_{\alpha}, X, Z, Y)| \mathfrak G_k((X, Z, Y))]  \Bigg. \Bigg.\\
     &+ \Bigg. \Bigg.\sum_{k = 1}^K \omega^*(k,\mathfrak G_k((X, Z, Y))) \mathbb{E}[ d(r_{\alpha}, X, Z, Y)| \mathfrak G_k((X, Z, Y))] - [m_1^*(r_{\alpha}, X) - (1 - \alpha) ] \Bigg] \Bigg] \\
     =& \mathbb{E}\left[\mathbb{P}(C \geq 1|X) \widehat \pi_1(X) \left[ d(r_{\alpha}, X, Z, Y) - [m_1^*(r_{\alpha}, X) - (1 - \alpha) ] \right] \right] \\
     =& \mathbb{E}\Bigg[\mathbb{P}(C \geq 1|X) \widehat \pi_1(X) \Bigg[\left\{ \widehat \omega(\infty, X, Z, Y) + \sum_{k = 1}^K \widehat \omega(k,\mathfrak G_k((X, Z, Y)))\right\}d(r_{\alpha}, X, Z, Y) \Bigg. \Bigg. \\
     & \Bigg. \Bigg. - [m_1^*(r_{\alpha}, X) - (1 - \alpha) ] \Bigg]\Bigg].
    \end{split}
\end{equation*}
The last expression equates to $0$. To see why, we note that from definition of $d(\cdots)$ we have, 
\begin{equation*}
    \begin{split}
   &\textbf{1}\{Y \leq r_{\alpha} \} - (1 - \alpha)\\
   = & \widehat{\mathcal{M}}\{d(r_{\alpha}, X, Z, Y)\} \\
       = & \sum_{k = 1}^{\infty} \widehat \omega(k,  \mathfrak G_k((X, Z, Y))) \mathbb{E}[d(r_{\alpha}, X, Z, Y) | \mathfrak G_k((X, Z, Y))]. 
    \end{split}
\end{equation*}
We complete the proof by substituting this into the previous expression,
\begin{equation*}
    \begin{split}
    & \mathbb{E}\left[\mathbb{P}(C \geq 1|X) \widehat \pi_1(X) \left[\left\{ \widehat \omega(\infty, X, Z, Y) + \sum_{k = 1}^K \widehat \omega(k,\mathfrak G_k((X, Z, Y)))\right\}d(r_{\alpha}, X, Z, Y) \right]\right]    \\
    =& \mathbb{E}\Bigg[\mathbb{P}(C \geq 1|X) \widehat \pi_1(X) \Bigg[ \widehat \omega(\infty, X, Z, Y) d(r_{\alpha}, X, Z, Y) \Bigg. \Bigg.\\
    & \Bigg. \Bigg. + \sum_{k = 1}^K \widehat \omega(k,\mathfrak G_k((X, Z, Y))) \mathbb{E}[ d(r_{\alpha}, X, Z, Y)| \mathfrak G_k((X, Z, Y))]  \Bigg] \Bigg] \\
    = & \mathbb{E}\left[\mathbb{P}(C \geq 1|X) \widehat \pi_1(X) \left[ \textbf{1}\{Y \leq r_{\alpha} \} - (1 - \alpha) \right] \right] \\
    =& \mathbb{E}\left[\mathbb{P}(C \geq 1|X) \widehat \pi_1(X) \left[ m_1^*(r_{\alpha}, X) - (1 - \alpha) \right] \right]. 
    \end{split}
\end{equation*}

\section{$\rscp$ for various problems}
\label{app:algo}
\begin{algorithm}
    \caption{Robust Split Conformal Prediction ($\rscp$) in the Covariate shift problem}
    \label{algg:split_robust_cmdp}
    \KwIn{Observed data set: $\mathcal{D}  = \{O_1^{t_1}, \cdots, O_{n_1}^{t_1},O_1^{t_2}, \cdots, O_{n_2}^{t_2} , O_1, \cdots, O_n\} = \mathcal{D}^{tr}_1 \cup \mathcal{D}^{tr}_2  \cup \mathcal{D}^{cal}$, Confidence level: $1-\alpha$, the baseline covariate $x$ of the test point}
    \KwOut{A valid prediction set $\widehat{\mathcal C}_{n,\alpha}(x)$}
    \begin{algorithmic}[1]
   \State Obtain estimates $\widehat \pi(x), \widehat m_R(\theta, x)$ of the nuisance functions $\pi^*(x), m_R^*(\theta, x)$ from the sub-sample $\mathcal{D}^{tr}_1$. 
   \State Define $\widehat r_{\alpha,0}$ using the sub-sample $\mathcal{D}^{tr}_2$,
   \begin{equation*}
      \begin{split}
      \widehat  r_{\alpha,0} = &\inf\Bigg\{\theta \in \mathbb{R}: \frac{\sum_{i = 1}^{n_2} \textbf{1}\{C_i^{t_2} = 0\} [\widehat m_R(\theta, X_i^{t_2}) - (1 - \alpha)] }{n_2+1} -  \frac{(1 - \alpha)}{n_2+1}  \geq 0\Bigg\}.
      \end{split}
  \end{equation*}
  \State  Define $\widehat r_{\alpha}$ using the calibration set $\mathcal{D}^{cal}$,  
  \begin{equation*}
      \begin{split}
      \widehat  r_{\alpha} = &\inf\Bigg\{\theta \in \mathbb{R}: \frac{\sum_{i = 1}^n f(O_i, \theta, \widehat \pi, \widehat m_R(\widehat r_{\alpha, 0}))}{n+1} -  \frac{\overline \pi}{n+1}  \geq 0\Bigg\},
      \end{split}
  \end{equation*}
  where,
  \begin{equation*}
      \begin{split}
     f(O, \theta, \pi, m_R(\theta')) = & \textbf{1}\{C = \infty\}\pi(X) [\textbf{1}\{R(X, Y) \leq \theta\} - m_R(\theta', X)] + \textbf{1}\{C = 0\} [m_R(\theta', X) - (1 - \alpha)]  .
      \end{split}
  \end{equation*}
  \State Return the prediction set, 
  \[
  \widehat{\mathcal C}_{n,\alpha}(x)  = \{y \in \mathbb{R}: R(x,y) \leq \widehat  r_{\alpha}\}. 
  \]
 \end{algorithmic}
\end{algorithm}
\begin{algorithm}
    \caption{Robust Split Conformal Prediction ($\rscp$) in the Two-stage monotone missing data problem}
    \label{algg:split_robust_mmdp}
    \KwIn{Observed data set: $\mathcal{D}  = \{O_1^{t_1}, \cdots, O_{n_1}^{t_1}, O_1^{t_2}, \cdots, O_{n_2}^{t_2}, O_1, \cdots, O_n\} = \mathcal{D}^{tr}_1 \cup  \mathcal{D}^{tr}_2 \cup \mathcal{D}^{cal}$, Confidence level: $1-\alpha$}
    \KwOut{A valid prediction set $\widehat{\mathcal C}_{n,\alpha}$}
    \begin{algorithmic}[1]
   \State Obtain estimates $\widehat \pi_1(x), \widehat \pi_2(x, z), \widehat m_1(\theta, x), \widehat m_2(\theta, x, z)$ of the nuisance functions $\pi_1^*(x), \pi_2^*(x, z), m_1^*(\theta, x), m_2^*(\theta, x, z)$ from the sub-sample $\mathcal{D}^{tr}_1$.
   \State Define $\widehat r_{\alpha, 0}$ using the sub-sample $\mathcal{D}^{tr}_2$,
   \begin{equation*}
       \begin{split}
           \widehat r_{\alpha, 0} = \inf \left\{\theta \in \mathbb R: \frac{\sum_{i = 1}^{n_2} g'(O_i, \theta, \widehat \pi_1, \widehat m_1, \widehat m_2) )}{n_2} \geq 0 \right\},
       \end{split}
   \end{equation*}
   where,
   \begin{equation*}
       \begin{split}
           g'(O, \theta,  \pi_1,  m_1, m_2) = & \textbf{1}\{C \geq 1\}  \pi_1(X) [ m_2(\theta, X, Z) - (1 - \alpha)] \\
           &-[\textbf{1}\{C \geq 1 \}  \pi(X) - \textbf{1}\{C = 0\} ] [m_1(\theta, X) - (1 - \alpha)]. 
       \end{split}
   \end{equation*}
  \State  Define $\widehat r_{\alpha}$ using the calibration set $\mathcal{D}^{cal}$,  
  \begin{equation*}
      \begin{split}
      \widehat  r_{\alpha} = &\inf\Bigg\{\theta \in \mathbb{R}: \frac{\sum_{i = 1}^n g(O_i, \theta,\widehat r_{\alpha, 0}, \widehat m_2, \widehat m_1, \widehat \pi_2, \widehat \pi_1)}{n+1}  - \frac{2\overline \pi}{n+1} \geq 0\Bigg\},  
      \end{split}
  \end{equation*}
  where,
  \[
  g(O, \theta, \theta',   m_2,  m_1,  \pi_2,  \pi_1) = \frac{\textbf{1}\{C = \infty\} \pi_1(X)}{\pi_2(X, Z)} [\textbf{1}\{Y \leq \theta\} - m_2(\theta', X, Z)] + g'(O, \theta',  \pi_1,  m_1, m_2).
  \]
  \State Return the prediction set, 
  \[
 \widehat{\mathcal C}_{n,\alpha} = \{y \in \mathbb{R}: y \leq \widehat r_{\alpha}\}.
  \]
 \end{algorithmic}
\end{algorithm}

\begin{algorithm}
    \caption{Robust Split Conformal Prediction ($\rscp$) in the Multi-stage monotone missing data problem}
    \label{algg:split_robust_Mmdp}
    \KwIn{Observed data set: $\mathcal{D}  = \{O_1^{t_1}, \cdots, O_{n_1}^{t_1}, O_1^{t_2}, \cdots, O_{n_2}^{t_2}, O_1, \cdots, O_n\} = \mathcal{D}^{tr}_1 \cup  \mathcal{D}^{tr}_2 \cup \mathcal{D}^{cal}$, Confidence level: $1-\alpha$}
    \KwOut{A valid prediction set $\widehat{\mathcal C}_{n,\alpha}$}
    \begin{algorithmic}[1]
   \State Obtain estimates $\{\widehat \pi_j\}_{j = 1}^{D+1}, \{\widehat m_j\}_{j = 1}^{D+1}$ of the true nuisance functions $\{\pi_j\}_{j = 1}^{D + 1}, \{m_j\}_{j = 1}^{D+1}$ from the sub-sample $\mathcal{D}^{tr}_1$. 
   \State Define $\widehat r_{\alpha, 0}$ using the sub-sample $\mathcal{D}^{tr}_2$,
   \begin{equation*}
       \begin{split}
           \widehat r_{\alpha, 0} = \inf \left\{\theta \in \mathbb R: \frac{\sum_{i = 1}^{n_2} g'(O_i, \theta, \{\widehat \pi_j\}_{j = 1}^D, \{\widehat m_j\}_{j = 1}^{D+1}) )}{n_2} \geq 0 \right\},
       \end{split}
   \end{equation*}
   where,
   \begin{equation*}
       \begin{split}
           g'(O, \theta,  \{ \pi_j\}_{j = 1}^D, \{ m_j\}_{j = 1}^{D+1}) = & \frac{\textbf{1}\{C \geq D\}  \pi_1(X)}{\pi_2(X, Z^1)\cdots \pi_{D}(X, Z^1,\cdots,Z^{D-1})} [ m_{D+1}(\theta, X, Z^1,\cdots,Z^{D}) - (1 - \alpha)] \\
           &- \cdots - \textbf{1}\{C \geq 1\} \pi_1(X) \left(\frac{\textbf{1}\{C \geq 2\}}{\pi_2(X, Z^1)} - 1 \right)[m_2(\theta, X, Z^1) - (1 - \alpha)] \\
           & - [\textbf{1}\{C \geq 1 \}  \pi(X) - \textbf{1}\{C = 0\} ] [m_1(\theta, X) - (1 - \alpha)]. 
       \end{split}
   \end{equation*}
  \State  Define $\widehat r_{\alpha}$ using the calibration set $\mathcal{D}^{cal}$,  
  \begin{equation*}
       \begin{split}
      \widehat  r_{\alpha} = &\inf\Bigg\{\theta \in \mathbb{R}: \frac{\sum_{i = 1}^n  g(O_i, \theta, \widehat r_{\alpha, 0}, \{\widehat \pi_j\}_{j = 1}^{D+1}, \{\widehat m_j\}_{j = 1}^{D+1})}{n+1}  - \frac{(D+1)\overline \pi}{n+1} \geq 0\Bigg\},  
      \end{split}
  \end{equation*}
  where,
  \begin{equation*}
      \begin{split}
          g(O, \theta, \theta',   \{\pi_j\}_{j = 1}^{D+1}, \{ m_j\}_{j = 1}^{D+1}) = &\frac{\textbf{1}\{C = \infty\} \pi_1(X)}{\pi_2(X, Z^1)\cdots \pi_{D+1}(X, Z^1,\cdots,Z^{D})} [\textbf{1}\{Y \leq \theta\} - m_{D+1}(\theta', X, Z^1,\cdots,Z^D)] \\
          &+ g'(O, \theta',  \{ \pi_j\}_{j = 1}^D, \{ m_j\}_{j = 1}^{D+1}). 
      \end{split}
  \end{equation*}
  \State Return the prediction set, 
  \begin{equation*}
      \begin{split}
       \widehat{\mathcal C}_{n,\alpha}= &\{y \in \mathbb{R}: y \leq \widehat r_{\alpha}\}.   
      \end{split}
  \end{equation*}
 \end{algorithmic}
\end{algorithm}
\begin{algorithm}
\caption{Robust Split Prediction ($\rsp$) under non-monotone missingness}
\label{algg:split_robust_nmdp}
\KwIn{Observed data set
$\mathcal D
=
\mathcal D^{tr}\cup\mathcal D^{cal}$,
where
$\mathcal D^{cal}=\{O_1,\ldots,O_n\}$,
and confidence level $1-\alpha$}
\KwOut{A prediction set
$\widehat{\mathcal C}_{n,\alpha}$}

\begin{algorithmic}[1]

\State Using the training data $\mathcal D^{tr}$, obtain estimates
$\widehat\pi_1(X)$ and $\{
\widehat\omega
\left(
k,\mathfrak G_k((X,Z,Y))
\right)
\}_{k\in\{1,\ldots,K,\infty\}}$ of the propensity-score nuisance functions
$\pi_1^*(X)$ and $\{
\omega^*
\left(
k,\mathfrak G_k((X,Z,Y))
\right)
\}_{k\in\{1,\ldots,K,\infty\}}$ respectively. Also obtain estimates
$\widehat m_1(\theta,X)$ and $\{
\widehat{\mathbb E}
\left[
d(\theta,X,Z,Y)
\mid
\mathfrak G_k((X,Z,Y))
\right]
\}_{k=1}^{K}$ of the corresponding outcome-regression nuisance functions.

\State Using these estimated nuisance functions, construct
$\widehat{\mathrm{IF}}(\theta,O)$ according to
\eqref{eq:influence_nonmonotone_estimated}.

\State Using the calibration sample
$\mathcal D^{cal}=\{O_1,\ldots,O_n\}$, define
\[
\widehat r_\alpha
=
\inf
\left\{
\theta\in\mathbb R:
\frac1n
\sum_{i=1}^{n}
\widehat{\mathrm{IF}}(\theta,O_i)
\geq0
\right\}.
\]

\State Return the prediction set
\[
\widehat{\mathcal C}_{n,\alpha}
=
\left\{
y\in\mathbb R:
y\leq\widehat r_\alpha
\right\}.
\]

\end{algorithmic}
\end{algorithm}
\end{document}